\documentclass[11pt]{article}

\usepackage[margin=1in]{geometry}
\usepackage{amsmath,amssymb,amsthm,mathtools,bm}
\usepackage{booktabs,array,microtype}
\usepackage[inline]{enumitem}
\usepackage{natbib}
\usepackage{graphicx}
\usepackage{xcolor}
\usepackage[colorlinks=true,linkcolor=blue!60!black,
            citecolor=blue!60!black,urlcolor=blue!60!black]{hyperref}

\allowdisplaybreaks
\newcommand{\covLamThousand}{96.0}
\newcommand{\biasLamThousand}{+0.0293}
\newcommand{\covEtaThousand}{94.8}
\newcommand{\biasEtaThousand}{-0.0163}
\newcommand{\covPsiThousand}{93.5}
\newcommand{\biasPsiThousand}{-0.0080}
\newcommand{\covLamTenk}{96.0}
\newcommand{\biasLamTenk}{+0.0032}
\newcommand{\covEtaTenk}{96.0}
\newcommand{\biasEtaTenk}{-0.0001}
\newcommand{\covPsiTenk}{96.8}
\newcommand{\biasPsiTenk}{+0.0001}
\newcommand{\covLamHundk}{98.0}
\newcommand{\biasLamHundk}{-0.0007}
\newcommand{\covEtaHundk}{97.0}
\newcommand{\biasEtaHundk}{+0.0026}
\newcommand{\covPsiHundk}{98.0}
\newcommand{\biasPsiHundk}{+0.0005}
\newcommand{\repSmall}{400}
\newcommand{\repLarge}{100}
\newcommand{\covPsiSmallPct}{93.5\%}
\newcommand{\biasRthousand}{-18}
\newcommand{\biasRtenk}{-21}
\newcommand{\feasCov}{93.0 / 94.0 / 94.0}
\newcommand{\nullMassTenkPct}{52.2\%}
\newcommand{\sizeTenkPct}{3.3\%}
\newcommand{\Ieff}{0.716}
\newcommand{\powTwoPct}{44.6\%}
\newcommand{\powFourPct}{89.2\%}
\newcommand{\powFullPct}{100.0\%}
\newcommand{\missKernelSummary}{$-5\%$ (triangular) and $-1\%$ (truncated Gaussian)}
\newcommand{\missTwoharmPsi}{$\widehat\psi=0.0907$}
\newcommand{\dmN}{3,000}
\newcommand{\dmR}{0.482}
\newcommand{\dmLamHat}{2.01}
\newcommand{\dmRHat}{2.05}
\newcommand{\dmPsiGraph}{0.028}
\newcommand{\dmPsiCI}{[0.001,0.150]}
\newcommand{\dmPsiDirect}{0.0527}
\newcommand{\dmPsiDirectSE}{0.0081}
\newcommand{\dmRgeo}{0.461}
\newcommand{\dmShare}{96\%}
\newcommand{\dmTransferPct}{5\%}
\newcommand{\bkC}{0.135}
\newcommand{\bkR}{-0.057}
\newcommand{\bkBDirect}{-0.045}
\newcommand{\bkBDirectSE}{0.026}
\newcommand{\bkRgeo}{0.470}
\newcommand{\bkGap}{0.53}
\newcommand{\bkMedKm}{13.2}
\newcommand{\bkPctOneKm}{9.5\%}
\newcommand{\bkMaxDeg}{109}
\newcommand{\bkHubShare}{11\%}
\newcommand{\ellipseCovThousand}{94.0\%}
\newcommand{\jointNaiveThousand}{86.2\%}
\newcommand{\ellipseCovTenk}{96.8\%}
\newcommand{\jointNaiveTenk}{91.0\%}
\newcommand{\ellipseCovHundk}{98.0\%}
\newcommand{\jointNaiveHundk}{94.0\%}
\newcommand{\chanSortNull}{-0.025}
\newcommand{\chanSortCenter}{0.129}
\newcommand{\chanOverlapCenter}{0.252}
\newcommand{\gridTable}{\begin{tabular}{@{}lccccccc@{}}\toprule $(\lambda,\eta,\psi)$ & $s_0$ & RMSE$_\psi$ & cov$_\psi$ & ellipse & len$_\psi$ & cov$_{s}$ & solve\% \\ \midrule $(1.5,6,0.020)$ & +0.08 & 0.0185 & 94.7 & 94.0 & 0.069 & 97.0 & 100\\ $(1.5,14,0.120)$ & +0.59 & 0.0179 & 98.3 & 95.7 & 0.075 & 97.0 & 99\\ $(3.0,9,0.084)$ & +0.34 & 0.0204 & 91.7 & 93.0 & 0.073 & 92.7 & 100\\ $(4.5,6,0.050)$ & +0.12 & 0.0225 & 94.3 & 92.7 & 0.089 & 96.0 & 100\\ $(4.5,14,0.140)$ & +0.51 & 0.0185 & 95.3 & 91.7 & 0.083 & 97.7 & 89\\ $(2.0,12,0.005)$ & -0.04 & 0.0090 & 96.0 & 96.7 & 0.036 & 95.7 & 100\\ $(3.0,9,0.150)$ & +0.47 & 0.0198 & 95.0 & 92.7 & 0.089 & 98.3 & 72\\ $(2.5,7,0.040)$ & +0.16 & 0.0193 & 93.3 & 94.0 & 0.070 & 92.7 & 100\\ \bottomrule\end{tabular}}
\newcommand{\rmseThreePsi}{0.0202}
\newcommand{\rmseFivePsi}{0.0133}
\newcommand{\covFivePct}{95.0/95.8/94.5\%}
\newcommand{\jTestPct}{5.8\%}
\newcommand{\pairShareA}{0.49}
\newcommand{\pairShareB}{-0.02}
\newcommand{\pairShareHatA}{0.49}
\newcommand{\pairShareHatB}{-0.01}
\newcommand{\sizeNuisA}{5.2\%}
\newcommand{\sizeNuisB}{5.5\%}
\newcommand{\sizeNuisC}{5.3\%}
\newcommand{\massNuisC}{0.51}
\newcommand{\cityTable}{\begin{tabular}{@{}lccccccccccc@{}}\toprule Region & $n$ & $\widehat\ell$ & $\widehat C$ & $\widehat r$ & $\widehat r_{\mathrm{rank}}$ & $r_{\mathrm{cfg}}$ & $C_{\mathrm{cfg}}$ & $\widehat\psi_{\mathrm{dir}}$ & $\widehat\psi_{\mathrm{gr}}$ CI & $r_{\mathrm{geo}}$ & $\sigma_{\min}$\\ \midrule SF Bay & 1,411 & 4.1 & 0.135 & $-0.057$ & $+0.105$ & $-0.170$ & 0.008 & 0.0020 & [0,0.021] & $+0.47$ & 0.054\\ LA & 1,304 & 4.7 & 0.133 & $-0.117$ & $+0.053$ & $-0.259$ & 0.009 & 0.0014 & [0,0.067] & $+0.46$ & 0.037\\ NYC & 1,040 & 2.4 & 0.123 & $-0.017$ & $+0.123$ & $-0.193$ & 0.003 & 0.0007 & [0,0.303] & $+0.26$ & 0.002\\ Denver & 1,025 & 8.8 & 0.186 & $-0.072$ & $+0.108$ & $-0.331$ & 0.020 & 0.0004 & [0,0.113] & $+0.33$ & 0.065\\ London & 704 & 2.3 & 0.198 & $+0.119$ & $+0.136$ & $-0.123$ & 0.003 & 0.0001 & [0,0.154] & $+0.51$ & 0.071\\ \bottomrule\end{tabular}}
\newcommand{\bkSpearman}{+0.105}
\newcommand{\bkCfgC}{0.008}
\newcommand{\bkCfgR}{-0.170}
\newcommand{\bkChSort}{0.463}
\newcommand{\bkChOv}{0.018}
\newcommand{\bkJacSmin}{0.054}
\newcommand{\bkJacCond}{201}
\newcommand{\bkRLo}{-0.12}
\newcommand{\bkRHi}{+0.12}
\newcommand{\bkRgeoLo}{+0.26}
\newcommand{\bkRgeoHi}{+0.51}
\newcommand{\bkPsiDirLo}{0.0001}
\newcommand{\bkPsiDirHi}{0.0020}
\newcommand{\dmTwelveCI}{[0.051,0.123]}
\newcommand{\hausSize}{3.2\%}
\newcommand{\hausLRtwo}{3.4\%}
\newcommand{\hausLRfive}{2.8\%}
\newcommand{\hausHeavy}{92.2\%}
\newcommand{\lrTwoPsi}{0.0816}
\newcommand{\lrFivePsi}{0.0819}
\newcommand{\heavyPsi}{0.0234}
\newcommand{\bkDistC}{0.005}
\newcommand{\bkDistR}{+0.011}
\newcommand{\ldBeta}{0.030}
\newcommand{\ldC}{0.0036}
\newcommand{\ldR}{+0.053}
\newcommand{\btSFp}{0.32}
\newcommand{\btSize}{4.5\%}
\newcommand{\gwTable}{\begin{tabular}{@{}lccccccccccc@{}}\toprule Region & $n$ & $\widehat\ell$ & $\widehat C$ & $\widehat r$ & $\widehat r_{\mathrm{rank}}$ & $r_{\mathrm{cfg}}$ & $C_{\mathrm{cfg}}$ & $\widehat\psi_{\mathrm{dir}}$ & $\widehat\psi_{\mathrm{gr}}$ CI & $r_{\mathrm{geo}}$ & $\sigma_{\min}$\\ \midrule Austin & 3,925 & 8.1 & 0.162 & $-0.004$ & $+0.307$ & $-0.273$ & 0.007 & 0.0040 & [0,0.004] & $+0.48$ & 0.075\\ Stockholm & 5,768 & 6.6 & 0.181 & $+0.275$ & $+0.331$ & $-0.211$ & 0.001 & 0.0001 & [0,0.013] & $+0.53$ & 0.077\\ SF Bay & 3,356 & 5.8 & 0.174 & $+0.018$ & $+0.329$ & $-0.206$ & 0.004 & 0.0002 & [0,0.020] & $+0.57$ & 0.083\\ Dallas & 2,164 & 5.4 & 0.166 & $-0.037$ & $+0.296$ & $-0.193$ & 0.005 & 0.0010 & [0,0.036] & $+0.58$ & 0.069\\ NYC & 1,678 & 2.3 & 0.168 & $-0.001$ & $+0.302$ & $-0.247$ & 0.003 & 0.0013 & [0,0.303] & $+0.38$ & 0.044\\ London & 1,331 & 2.0 & 0.178 & $+0.148$ & $+0.267$ & $-0.229$ & 0.002 & 0.0001 & [0,0.079] & $+0.55$ & 0.063\\ \bottomrule\end{tabular}}
\newcommand{\gwAusBDir}{-0.063}
\newcommand{\gwAusBDirSE}{0.016}
\newcommand{\gwAusPsiD}{0.0040}
\newcommand{\gwAusPsiDSE}{0.0020}
\newcommand{\gwAusR}{-0.004}
\newcommand{\gwAusSpear}{$+0.307$}
\newcommand{\gwRLo}{$-0.04$}
\newcommand{\gwRHi}{$+0.27$}
\newcommand{\gwRgeoLo}{$+0.38$}
\newcommand{\gwRgeoHi}{$+0.58$}
\newcommand{\bootDraws}{100--150}
\newcommand{\gofTable}{\begin{tabular}{@{}lccccccc@{}}\toprule City & $r_{\mathrm{mod}}$ & $\widehat r$ & resid & $z$ & rail & $\%{>}\widehat R$ & $\Delta_{\mathrm{sort}}(\widehat b_{\mathrm{dir}})$\\ \midrule BK SF Bay & $+0.470$ & $-0.057$ & $-0.528$ & $-21$ & $b{=}0$ & 91 & $-0.0106$\\ BK LA & $+0.461$ & $-0.117$ & $-0.578$ & $-14$ & $b{=}0$ & 84 & $+0.0098$\\ BK NYC & $+0.138$ & $-0.017$ & $-0.155$ & $-5$ & $b{=}0.55$, $R{=}0.1$ & 97 & --\\ BK Denver & $+0.328$ & $-0.072$ & $-0.400$ & $-10$ & $b{=}0$, $R{=}R_{\max}$ & 84 & --\\ BK London & $+0.506$ & $+0.119$ & $-0.387$ & $-8$ & $b{=}0$ & 93 & $-0.0009$\\ GW Austin & $+0.484$ & $-0.004$ & $-0.488$ & $-37$ & $b{=}0$ & 96 & $-0.0143$\\ GW Stockholm & $+0.534$ & $+0.275$ & $-0.259$ & $-16$ & $b{=}0$ & 98 & $+0.0027$\\ GW SF Bay & $+0.572$ & $+0.018$ & $-0.554$ & $-32$ & $b{=}0$ & 95 & $-0.0032$\\ GW Dallas & $+0.576$ & $-0.037$ & $-0.614$ & $-29$ & $b{=}0$ & 91 & $+0.0070$\\ GW NYC & $+0.381$ & $-0.001$ & $-0.382$ & $-11$ & $b{=}0$ & 94 & $-0.0041$\\ GW London & $+0.550$ & $+0.148$ & $-0.402$ & $-11$ & $b{=}0$ & 95 & $-0.0020$\\ \bottomrule\end{tabular}}
\newcommand{\gapLoAll}{+0.26}
\newcommand{\gapHiAll}{+0.61}
\newcommand{\rgeoLoAll}{+0.26}
\newcommand{\rgeoHiAll}{+0.58}
\newcommand{\robsLoAll}{-0.12}
\newcommand{\robsHiAll}{+0.27}
\newcommand{\fracOutLoAll}{84}
\newcommand{\fracOutHiAll}{98}
\newcommand{\dsortBoundAll}{0.02}
\newcommand{\bkFracLo}{84}
\newcommand{\bkFracHi}{97}
\newcommand{\gwFracLo}{91}
\newcommand{\gwFracHi}{98}
\newcommand{\bkRhatSF}{0.97}
\newcommand{\bkDsortMax}{0.02}
\newcommand{\gwAusDsort}{-0.0143}
\newcommand{\gwAusDsortFlip}{+0.0156}
\newcommand{\bkNycSmin}{0.002}
\newcommand{\bkNycCond}{5,349}
\newcommand{\bkSminLo}{0.037}
\newcommand{\bkSminHi}{0.071}
\newcommand{\gwStoR}{+0.275}
\newcommand{\gwStoCfgR}{-0.211}
\newcommand{\gwStoRgeo}{+0.534}
\newcommand{\pairDeltaA}{+0.155}
\newcommand{\pairDeltaB}{+0.012}
\newcommand{\ieffLo}{0.39}
\newcommand{\ieffHi}{2.92}
\newcommand{\gofZWorstOld}{37}
\newcommand{\richTable}{\begin{tabular}{@{}llcccccccc@{}}\toprule City & & $\widehat\ell$ & $\widehat C$ & $\widehat r$ & $q_{25}$ & $q_{50}$ & $q_{75}$ & cv & $p_{99}$\\ \midrule SF Bay & obs & 4.10 & 0.13 & $-0.057$ & 2.62 & 13.17 & 30.87 & 1.66 & 31.80\\  & model & 5.63 (0.94) & 0.12 (0.00) & $+0.017$ (0.02) & 4.83 (1.14) & 17.35 (1.45) & 38.01 (2.13) & 1.68 (0.07) & 44.46 (8.82)\\  & $z$ & $-1.6$ & $+3.6$ & $-3.6$ & $-1.9$ & $-2.9$ & $-3.4$ & $-0.3$ & $-1.4$\\ \midrule Austin & obs & 8.11 & 0.16 & $-0.004$ & 1.27 & 3.65 & 9.99 & 2.23 & 78.00\\  & model & 14.31 (3.11) & 0.16 (0.02) & $-0.025$ (0.02) & 1.17 (0.11) & 3.41 (0.39) & 8.97 (0.69) & 2.13 (0.07) & 144.26 (31.76)\\  & $z$ & $-2.0$ & $+0.2$ & $+0.9$ & $+0.9$ & $+0.6$ & $+1.5$ & $+1.4$ & $-2.1$\\ \bottomrule\end{tabular}}
\newcommand{\richSFMaxZ}{3.6}
\newcommand{\richAusMaxZ}{2.1}
\newcommand{\richSFHub}{0.11}
\newcommand{\richSFMaxDeg}{108}
\newcommand{\cfHubSF}{+0.097}
\newcommand{\cfCloSF}{-0.052}
\newcommand{\cfLrSF}{+0.048}
\newcommand{\cfBSF}{+0.018}
\newcommand{\cfHubAus}{+0.144}
\newcommand{\cfCloAus}{-0.027}
\newcommand{\cfLrAus}{+0.004}
\newcommand{\cfBAus}{+0.013}
\newcommand{\richAusEllZ}{-2.0}
\newcommand{\richAusPZ}{-2.1}
\newcommand{\richAusGoodZ}{1.5}
\newcommand{\richAusEllOver}{76}
\newcommand{\cfCompTable}{\begin{tabular}{@{}lcc@{}}\toprule counterfactual (off) & SF Bay & Austin\\ \midrule propensities $\sigma$ & +0.097 & +0.144\\ closure $\tau$ & -0.052 & -0.027\\ long range $w$ & +0.048 & +0.004\\ dependence set to $\widehat b_{\mathrm{dir}}$ & +0.018 & +0.013\\ \bottomrule\end{tabular}\\[7pt] \begin{tabular}{@{}lcccccc@{}}\toprule model & AUC & $C$ & $r$ & $q_{50}$ & cv & $p_{99}$\\ \midrule geometry ($\psi{=}0$) & 0.54 & 0.135 & $+0.473$ & 0.4 & 1.96 & 35\\ DC-SBM ($K{=}6$) & 0.82 & 0.050 & $-0.051$ & 23.6 & 1.64 & 22\\ spatial logistic $+$ node effects & 0.86 & 0.025 & $-0.056$ & 13.6 & 1.56 & 29\\ richer class & 0.52/0.76 & 0.121 & $+0.022$ & 16.9 & 1.66 & 45\\ \textbf{observed} & & \textbf{0.135} & $\mathbf{-0.057}$ & \textbf{13.2} & \textbf{1.66} & \textbf{32}\\ \bottomrule\end{tabular}}
\newcommand{\aucLogit}{0.86}
\newcommand{\aucDcsbm}{0.82}
\newcommand{\aucRichCond}{0.76}
\newcommand{\aucGeo}{0.54}
\newcommand{\compLogC}{0.025}
\newcommand{\compDcsbmC}{0.050}
\newcommand{\oaAuthors}{73{,}342}
\newcommand{\oaEdgesAll}{403{,}117}
\newcommand{\oaNLn}{27{,}861}
\newcommand{\oaNLE}{118{,}212}
\newcommand{\oaNLell}{8.49}
\newcommand{\oaNLC}{0.361}
\newcommand{\oaNLr}{+0.166}
\newcommand{\oaNLrs}{+0.367}
\newcommand{\oaBdir}{+0.0932}
\newcommand{\oaBdirSE}{0.0056}
\newcommand{\oaBdirSEcl}{0.0244}
\newcommand{\oaBdirZcl}{3.8}
\newcommand{\oaPsiDir}{0.0087}
\newcommand{\oaNLatoms}{420}
\newcommand{\oaNLtop}{15}
\newcommand{\oaGroBdir}{+0.0018}
\newcommand{\oaGroTop}{100}
\newcommand{\oaSubN}{6{,}500}
\newcommand{\oaSubE}{6{,}787}
\newcommand{\oaSubEll}{2.09}
\newcommand{\oaSubC}{0.353}
\newcommand{\oaSubR}{+0.230}
\newcommand{\oaSubRs}{+0.400}
\newcommand{\oaSubBdir}{+0.0853}
\newcommand{\oaSubBdirSEcl}{0.0244}
\newcommand{\oaSubCfgR}{-0.072}
\newcommand{\oaSubCfgRsd}{0.009}
\newcommand{\oaSubCfgC}{0.0003}
\newcommand{\oaSubQtf}{0.45}
\newcommand{\oaSubQf}{0.81}
\newcommand{\oaSubQsf}{57.4}
\newcommand{\oaSubQnn}{148}
\newcommand{\oaSubFracLocal}{56}
\newcommand{\oaSubFracFar}{39}
\newcommand{\oaSubDbar}{+0.22}
\newcommand{\oaSubDbarZ}{12}
\newcommand{\oaRegionTable}{\begin{tabular}{@{}lrrrrrrrrr@{}}\toprule region & $n$ & edges & $\widehat\ell$ & $\widehat C$ & $\widehat r$ & $\widehat b_{\mathrm{dir}}$ & cl.~se & atoms & top\\ \midrule Netherlands & 27{,}861 & 118{,}212 & 8.49 & 0.361 & $+0.166$ & $+0.0932$ & 0.0244 & 420 & 15\\ Randstad & 14{,}584 & 54{,}668 & 7.50 & 0.376 & $+0.160$ & $+0.0839$ & 0.0378 & 166 & 29\\ South Holland & 6{,}868 & 22{,}248 & 6.48 & 0.402 & $+0.165$ & $+0.1270$ & 0.0518 & 71 & 44\\ Amsterdam & 4{,}380 & 12{,}903 & 5.89 & 0.431 & $+0.181$ & $+0.0470$ & 0.0438 & 35 & 95\\ Utrecht & 2{,}886 & 8{,}832 & 6.12 & 0.409 & $+0.152$ & $+0.0668$ & 0.0615 & 33 & 83\\ Eindhoven & 1{,}959 & 6{,}486 & 6.62 & 0.395 & $+0.087$ & $+0.0572$ & 0.0567 & 17 & 91\\ Groningen & 2{,}209 & 6{,}161 & 5.58 & 0.433 & $+0.200$ & $+0.0018$ & 0.0131 & 4 & 100\\ \bottomrule\end{tabular}}
\newcommand{\oaBaseEll}{5.33}
\newcommand{\oaBaseR}{+0.661}
\newcommand{\oaBaseGofZ}{-24}
\newcommand{\oaBaseFracOut}{98.6}
\newcommand{\oaBaseRmax}{0.22}
\newcommand{\oaBaseDsort}{+0.018}
\newcommand{\oaBaseDsortFlip}{-0.018}
\newcommand{\oaFreeB}{-0.55}
\newcommand{\oaFreeMaxZ}{10.2}
\newcommand{\oaFreeZell}{+10.2}
\newcommand{\oaFreeZr}{-5.9}
\newcommand{\oaFreeZp}{+6.9}
\newcommand{\oaFreeEll}{1.54}
\newcommand{\oaFixTau}{0.64}
\newcommand{\oaFixL}{200}
\newcommand{\oaFixMaxZ}{7.8}
\newcommand{\oaFixZr}{-7.8}
\newcommand{\oaFixZcv}{-5.2}
\newcommand{\oaFixZqsf}{+3.8}
\newcommand{\oaFixRhat}{+0.368}
\newcommand{\oaFixCfClo}{-0.212}
\newcommand{\oaFixSortContrib}{+0.0067}
\newcommand{\oaFixSignedGap}{0.0132}
\newcommand{\oaProfSlopeFix}{+0.050}
\newcommand{\oaProfSlopeFree}{+0.085}
\newcommand{\oaFixSpear}{+0.507}
\newcommand{\oaTable}{\begin{tabular}{@{}lllll@{}}\toprule  & model $r$ & model $C$ & fitted moments & flags\\ \midrule observed (subsample{,} $n=6{,}500$) & $+0.230$ & 0.353 & & \\ configuration null & $-0.072\,(0.009)$ & 0.0003 & degrees & ---\\ hard-disc class & $+0.661$ & 0.355 & $(\ell,C,r)$; misses $\ell$ & $R$, $b$ railed; $99\%$ outside $\widehat R$\\ richer, $b$ free & $+0.331$ & 0.360 & eight; misses $\ell,r,\mathrm{cv},p_{99}$ & $b$ railed $-0.55$; $\max|z|{=}10.2$\\ richer, $b=\widehat b_{\mathrm{dir}}$ & $+0.360$ & 0.360 & eight; misses $r,\mathrm{cv}$ & $\max|z|{=}7.8$ ($z_r{=}-7.8$)\\ \bottomrule\end{tabular}}
\newcommand{\oaFreeTwoB}{-0.55}
\newcommand{\oaFreeTwoMaxZ}{19.4}
\newcommand{\oaSensLo}{+0.082}
\newcommand{\oaSensHi}{+0.100}
\newcommand{\oaSensZlo}{3.6}
\newcommand{\oaSensZhi}{4.0}

\newcommand{\T}{\mathbb T}
\newcommand{\R}{\mathbb R}
\newcommand{\E}{\mathbb E}
\newcommand{\Pp}{\mathbb P}
\newcommand{\Cov}{\operatorname{Cov}}
\newcommand{\Var}{\operatorname{Var}}
\newcommand{\TV}{\mathrm{TV}}
\newcommand{\Pois}{\operatorname{Pois}}
\newcommand{\one}{\mathbf 1}
\newcommand{\trans}{\mathsf T}
\newcommand{\eps}{\varepsilon}

\theoremstyle{plain}
\newtheorem{theorem}{Theorem}
\newtheorem{lemma}{Lemma}

\theoremstyle{remark}
\newtheorem*{remark}{Remark}

\usepackage{mathrsfs}

\theoremstyle{plain}
\newtheorem{sitheorem}{Theorem}[section]
\newtheorem{siproposition}[sitheorem]{Proposition}
\newtheorem{silemma}[sitheorem]{Lemma}
\newtheorem{sicorollary}[sitheorem]{Corollary}

\title{Decomposing Degree Assortativity in Sparse Spatial Networks}
\author{
Marios Papamichalis\thanks{\raggedright Human Nature Lab, Yale
University, New Haven, CT 06511, USA.
\href{mailto:marios.papamichalis@yale.edu}{\nolinkurl{marios.papamichalis@yale.edu}}}
\and Regina Ruane\thanks{\raggedright Department of Statistics and
Data Science, The Wharton School, University of Pennsylvania,
Philadelphia, PA, USA.
\href{mailto:ruanej@wharton.upenn.edu}{\nolinkurl{ruanej@wharton.upenn.edu}}}
}
\date{}

\begin{document}
\maketitle

\begin{abstract}
Spatial networks are typically assortative: well-connected nodes
link to other well-connected nodes, and the usual reading is
sorting, popular nodes seeking each other out.  In space there is a
rival explanation: nearby nodes draw on the same pool of potential
neighbors, so their degrees move together even when popularity and
location are unrelated.  This paper asks when the two can be told
apart from network data, and gives a three-part answer.  First, in a
sparse random-connection model in which latent popularity and
position are coupled by a copula, we prove that the limiting
assortativity splits exactly into an intensity-covariance channel
and a shared-neighbor channel; because the first mixes sorting with
density effects, we define the sorting contribution as the
increment produced by the dependence when mean degree and
transitivity are held fixed, exactly zero at independence.  Second, we prove that three observable summaries
(mean degree, transitivity, assortativity) identify the model
globally, including on the boundary of no dependence, by a
computer-assisted proof in exact rational arithmetic that also
verifies the conditions for valid inference; simulations confirm
the asymptotics.  Third, on data
the machinery acts as a gate: no sorting estimate is reported unless
the model first fits.  In eleven metropolitan location-based
networks the model class is rejected and direct estimates find
essentially no dependence; the observed assortativity is evidence of
neither sorting nor shared opportunity.  In a national co-authorship
network the verdict reverses: productive authors
concentrate where researcher density is high, yet graph-only
attribution is again refused, the misfit pointing to the team
structure of multi-author papers.
\end{abstract}


\section{Introduction}
\label{sec:intro}

In many networks the degrees at the two ends of an edge are
correlated, a property measured by the assortativity coefficient
\citep{newman2002}.  When nodes carry locations, positive
assortativity is commonly interpreted as sorting: popular nodes
concentrate in the same places, so edges tend to join two locally
popular endpoints.  In a friendship network from a location-based
service, for example, positive assortativity suggests that highly
active users live in the same neighborhoods.  A second mechanism,
however, produces the same signature without any sorting.  The
endpoints of an edge are spatially close; close nodes share potential
neighbors; and the shared opportunities correlate the two degrees.
That geometry alone generates positive assortativity through shared
neighbors is known: \citet{kaufmann2026} decompose a neighbor's degree
into Poisson contributions tied to the overlap of influence regions
and study real networks including Brightkite, and \citet{patularu2026}
derive explicit limits of the Pearson assortativity coefficient for
geometric graphs under local convergence.

Single-graph separation of other mechanism pairs exists in other
model classes: \citet{peixoto2022} disentangles homophily from
triadic closure and community structure within a Bayesian
block-model class, and \citet{fosdick2015} test and model dependence
between a network and observed nodal attributes.  What the literature
does not provide, and what this paper supplies for one specified
model class, is a population definition of the intensity and overlap
contributions to assortativity, a proof that observable graph
summaries determine them globally, and inference that remains valid
at the independence boundary.  Any measured assortativity mixes the
two mechanisms, and an estimate of popularity--location dependence
that ignores the opportunity channel attributes to preference what is
due to geometry.  Three questions therefore arise.  Can the contribution
of each mechanism be defined as a population quantity?  Are the two
contributions identified from observable graph summaries?  Can the
resulting estimators be equipped with valid inference, including at
the boundary case of no dependence, which is the natural null
hypothesis?  The questions matter beyond the model studied here
because assortativity is a standard summary in applied network
analysis and is routinely given a mechanistic interpretation.

This paper answers the three questions in a specific model class:
finite-range random-connection graphs in which a node's latent
popularity and position are coupled by a copula that preserves both
marginals.  In the sparse finite-range limit, the edge-endpoint degree
correlation is
\begin{equation}
 r=
 \frac{\overbrace{\Cov_{\Pi}\{\Lambda(W_1,X),\Lambda(W_2,X)\}}^{\text{intensity covariance}}
       +\overbrace{\E_{\Pi}\,\Gamma(W_1,W_2,X;U)}^{\text{shared neighbors}}}
      {\Var_{\Pi}\{\Lambda(W_1,X)\}
       +\E_{\Pi}\Lambda(W_1,X)},
 \label{eq:r-decomp-intro}
\end{equation}
where $\Pi$ is the edge-Palm law of the endpoint marks, position, and
tangent displacement $U$; $\Lambda$ is the endpoint connection
intensity; and $\Gamma$ is the shared-neighbor intensity, which depends
on $U$ at the connection-range scale
(Theorem~\ref{thm:realized-main}).  We write $r_{\mathrm{int}}$ and
$r_{\mathrm{overlap}}$ for the two summands.  The first term is
deliberately not called sorting: it aggregates every source of
endpoint-intensity dependence, including inhomogeneity of positions,
and it can be large at $\psi=0$.  The sorting contribution of the
dependence parameter is instead defined as a counterfactual
increment, $\Delta_{\mathrm{sort}}(\psi)=r(\psi)-r(0)$ along the
profile that holds mean degree and transitivity fixed; it is zero at
independence by construction, and
Theorem~\ref{thm:identification-main}(iii) proves it strictly
increasing.  The paper reports channel estimates and delta-method
intervals for $r_{\mathrm{int}}$ and $r_{\mathrm{overlap}}$, and
$\Delta_{\mathrm{sort}}$ evaluated at the estimated dependence, in
simulation and on data.  Collapsing the two endpoints to one
macroscopic location deletes $r_{\mathrm{overlap}}$ and therefore
overstates the intensity channel.

\begin{figure}[t]
\centering
\includegraphics[width=.62\linewidth]{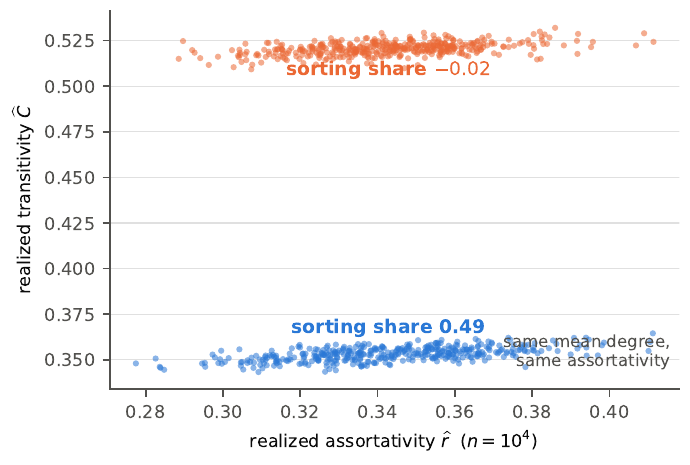}
\caption{Two models matched on population mean degree and
assortativity, one with intensity share $\pairShareA$ (blue) and one
with intensity share $\pairShareB$ (orange); $400$ realized networks
each at $n=10^4$.  The realized assortativity distributions coincide,
while the realized transitivity distributions do not overlap; the
three-moment method recovers both shares
(Section~\ref{sec:simulations}).}
\label{fig:thesis}
\end{figure}

Two facts organize the analysis.  First, mean degree and assortativity
alone cannot separate the channels.  Figure~\ref{fig:thesis} shows two
parameter points inside the certified region with identical $(\ell,r)$
and intensity shares of $\pairShareA$ and $\pairShareB$ (sorting
contributions $\Delta_{\mathrm{sort}}=\pairDeltaA$ and $\pairDeltaB$):
the realized assortativity distributions coincide, while the realized
transitivity distributions do not overlap.  This motivates
identification from the triple $(\ell,C,r)$ of mean degree,
transitivity, and assortativity.
Second, the channels are functionals of the parameters
$(\lambda,\eta,\psi)$, so channel attribution requires inverting the
map from parameters to observables, and the inversion must be shown
to be well posed.

The model couples a node's latent popularity rank and position by a
periodic one-harmonic conditional law that preserves both marginals
exactly.  Exact finite-sample symmetries (translation, kernel
symmetries, a sign flip) make the harmonic's phase and signed amplitude
unidentifiable; the orbit invariant is the squared amplitude
$\psi$, an $L^2$ dependence functional of the copula
(Lemma~\ref{lem:equivalence-main}).  We treat $\psi$ as the estimand
throughout.

The contributions are as follows.  First, an exact two-channel
decomposition of the limiting assortativity, carried by the tangent
displacement coordinate of the edge-Palm law
(Theorem~\ref{thm:realized-main}), together with a sorting
contribution $\Delta_{\mathrm{sort}}$ defined so that it vanishes
exactly at independence.  Second, certified global
identification: an exact-rational interval computation proves that
$(\lambda,\eta,\psi)\mapsto(\ell,C,r)$ is injective on the region
$B=[1,5]\times[5,15]\times[0,0.15]$, including the boundary $\psi=0$,
with quantitative bounds on the Jacobian
(Theorem~\ref{thm:identification-main}); a companion result covers the
planar hard-disc case.  Global injectivity is stronger than the local
rank conditions standard in moment-based inference, and it rules out
spurious solution branches in the inversion.  Third, verified
inference: the same computation discharges the uniform-convergence,
rank, and unique-root conditions of the motif-moment central limit
theorem and generalized method of moments
(Theorem~\ref{thm:inference-main}), with the numerical evaluator,
simulation, and optimization accuracies stated as explicit
implementation conditions, and yields a boundary theory at $\psi=0$
(censored-Gaussian $\sqrt n$ law for $\widehat\psi$, $n^{1/4}$
amplitude rate, and $\tfrac12\chi^2_0+\tfrac12\chi^2_1$ test,
Theorem~\ref{thm:boundary-main}) whose model-side conditions are all
verified; the boundary laws themselves are the standard consequences
of a squared boundary parameter
\citep{chernoff1954,selfliang1987,andrews1999,andrews2001}.  Fourth, a
simulation study spanning the certified region: eight parameter
configurations with marginal and simultaneous coverage, a comparison
of the exactly identified three-moment estimator with an
overidentified five-motif GMM and its $\chi^2_2$ specification test,
boundary test size across nuisance configurations, and the
mechanism-separation experiment of Figure~\ref{fig:thesis}
(Section~\ref{sec:simulations}).  Fifth, a diagnostic application to
eleven prespecified metropolitan networks from two location-based
services, five from Brightkite and six from Gowalla, in which the
fitted model class is rejected by the very moment it is meant to
explain, and the pattern of the rejection, together with null
comparisons, localizes the mechanisms the class omits
(Section~\ref{sec:application}).  Sixth, a constructive response to
the rejection: a closure-augmented two-scale class with degree
propensities fitted under the same goodness-of-fit discipline, a
competitor study across model families on held-out edges and
posterior-predictive motifs, a boundary theory made uniform over the
certified null face by a second exact-rational certificate, and
first-order identification of the signed amplitude under irregular
designs (Sections~\ref{sec:inference}
and~\ref{sec:application}; Supplement, Sections~S4 and~S5).
Seventh, an out-of-domain application of the two-view design to a
national co-authorship network, institution coordinates as
positions, productivity as the mark, in which the direct view
measures a clearly positive interior dependence and graph-only
attribution is refused by the goodness-of-fit gate at every level of
the model ladder, the residual naming the mechanism the class lacks
(Section~\ref{sec:collab}).

Relative to the latent-space, graphon-copula, and spatial random-graph
literatures
\citep{raschke2014,idowu2025,kaufmann2026,vdh2023local,vdh2023cluster,
patularu2026,peixoto2022,krivitsky2009,fosdick2015,rastelli2016,
shalizi2024}, the contribution is not the observation that spatial
graphs are assortative, which is prior work, but the combination of a
marginal-preserving dependence family, an exact two-channel
decomposition carried by the tangent displacement coordinate, and
global identification with verified inference for the channels from a
single graph.

The application uses check-in data, which provide each user with a
home location and an observed activity rank.  In these data the
latent-dependence estimand has a directly computable counterpart:
$\psi$ can be estimated once from the location--rank pairs and once
from the graph moments with ranks hidden.  The empirical finding is a
rejection, and the paper presents it as one.  Fitted to each city's
mean degree and transitivity, the model class over-predicts that
city's assortativity by $\gapLoAll$ to $\gapHiAll$ (predictions of
$\rgeoLoAll$ to $\rgeoHiAll$ against observed values of $\robsLoAll$
to $\robsHiAll$); because assortativity is increasing in
$\psi$ after compensation
(Theorem~\ref{thm:identification-main}(iii)), no admissible
dependence strength closes the gap, and the constrained graph-only
fits sit on the $\psi=0$ boundary in ten of eleven cities as a
consequence of the misfit, not as a measurement of independence.  The
fitted connection ranges exclude $\fracOutLoAll$ to $\fracOutHiAll$
percent of the observed ties.  The short-range
geometry-plus-dependence class is therefore decisively insufficient
for these networks.  The direct estimator, which does not use the
graph, gives the sorting evidence: it finds no dependence of activity
rank on the spatial score in ten cities, and in Austin a small
negative dependence
($\widehat\psi_{\mathrm{dir}}=\gwAusPsiD\pm\gwAusPsiDSE$) whose sign
is invisible to the homogeneous graph law
(Lemma~\ref{lem:equivalence-main}) and whose size, propagated through
the fitted geometry, changes model assortativity by less than
$\dsortBoundAll$.  The observed assortativity is consistent with
degree heterogeneity as the dominant channel (degree-preserving
configuration nulls produce disassortativity of the observed sign,
overshooting its size, and collapse transitivity; the median
within-region tie length is $\bkMedKm$ km on Brightkite,
incompatible with a short connection range), and on both platforms
the observed value sits above the configuration null by a positive
residual that none of the originally fitted mechanisms explains; a
closure-augmented two-scale class with degree propensities, fitted
in Section~\ref{sec:application}, carries most of this structure and
localizes what remains.  Section~\ref{sec:collab} then repeats the
two-view protocol on a co-authorship network built for the purpose,
where the direct verdict reverses, an interior positive dependence,
while the goodness-of-fit gate refuses graph-only attribution under
every class fitted, including the augmented one.

Notation and terminology are fixed as follows.  $\ell$ denotes the
population mean degree, $C$ the transitivity (the global clustering
coefficient), and $r$ the Pearson assortativity coefficient; hats
denote sample or estimated versions.  $\widehat\psi_{\mathrm{dir}}$
and $\widehat\psi_{\mathrm{gr}}$ denote the direct and graph-only
estimators of $\psi$.  ``Configuration null'' means the
degree-preserving rewiring null, and ``geometry-only baseline'' means
the fitted model with $\psi=0$ and the remaining parameters
recalibrated to the observed mean degree and transitivity, whose
assortativity is written $r_{\mathrm{geo}}$.  The decomposition
concerns the Pearson coefficient, because the Pearson numerator is the
covariance functional that the edge-Palm limit controls and that
admits the closed forms of Section~\ref{sec:realized}; the Spearman
(rank) coefficient is reported alongside it throughout the
application, the two disagree in sign in several cities, and every
substantive claim in this paper is a claim about the Pearson
coefficient.  A rank-based analogue of the decomposition is an open
problem.
Section~\ref{sec:model} defines the model and the invariant estimand;
Section~\ref{sec:realized} derives the decomposition;
Section~\ref{sec:identification} states the identification
certificate; Section~\ref{sec:inference} develops inference;
Sections~\ref{sec:simulations} and~\ref{sec:application} report the
simulation study and the application; Section~\ref{sec:discussion}
discusses limitations and extensions.

\section{Model and invariant estimand}
\label{sec:model}

Let $\T^d=\R^d/\mathbb Z^d$ carry normalized Haar measure, and write
$\delta_{\T^d}(x,y)\in[-1/2,1/2)^d$ for the centered representative of
$y-x$.  Let $\eps_n>0$ with $\rho_n=n\eps_n^d\to\rho\in(0,\infty)$.
Each node draws $(X_i,W_i)$ independently from $dx\,H_x(dw)$ on
$\T^d\times[0,\infty)$; given the latent variables, unordered edges are
independent with
\begin{equation}
 p_{ij,n}=1-\exp\!\left\{-\frac{\lambda}{\rho_n}W_iW_j
 k\!\left(\frac{\delta_{\T^d}(X_i,X_j)}{\eps_n}\right)\right\},
 \label{eq:finite-link-main}
\end{equation}
with $\lambda>0$ and $k$ nonnegative, bounded, even, compactly
supported, under the local total-variation continuity condition
$\sup_{d_{\T^d}(x,y)\le a}\|H_x-H_y\|_{\TV}\to0$ as $a\downarrow0$.

The marginal-preserving submodel puts $W_i=\omega(V_i)$ and couples
$V_i\in[0,1]$ to $X_i$ by the one-harmonic conditional density
\begin{equation}
 h_{b,\varphi,m}(v\mid x)
 =1+b\,a(v)\,\phi_{m,\varphi}(x),
 \qquad
 a(v)=\sqrt3(2v-1),\quad
 \phi_{m,\varphi}(x)=\sqrt2\cos(2\pi m^{\trans}x-\varphi),
 \label{eq:harmonic-density-main}
\end{equation}
for $m\in\mathbb Z^d\setminus\{0\}$ and $|b|\le1/\sqrt6$.  Both
marginals are exactly uniform.  The dependence-strength target is the
squared $L^2$ norm of the copula perturbation,
\begin{equation}
 \psi=\int_0^1\!\!\int_{\T^d}
 \{h_{b,\varphi,m}(v\mid x)-1\}^2\,dx\,dv=b^2 .
 \label{eq:psi-main}
\end{equation}
The family \eqref{eq:harmonic-density-main} is a
Farlie--Gumbel--Morgenstern-type first-order coupling
\citep{nelsen2006}; it is a minimal family that varies
popularity--position dependence while holding both marginals fixed.

\begin{lemma}[Exact symmetries and the invariant strength]
\label{lem:equivalence-main}
For every $n\ge2$ the labeled adjacency law is invariant under torus
translations and kernel symmetries acting on $(\varphi,m)$ as
$\varphi\mapsto\varphi+2\pi(Qm)^{\trans}\widetilde\tau$,
$m\mapsto Qm$, and under
$(b,\varphi)\mapsto(-b,\varphi+\pi)$.  Phase and signed amplitude are
therefore not identified from the adjacency law, while $\psi=b^2$ is
invariant under every such transformation.  Rescalings of popularity,
kernel height, and bandwidth generate the further exact aliases
recorded in the supplement and are removed by the normalizations
adopted below.
\end{lemma}

The lemma (proved in Section~S1 of the Supplement, together with the full orbit
description) fixes the estimand; it does not establish that $\psi$ is
identified.  Identification is established in
Theorem~\ref{thm:identification-main}.

\section{Edge-Palm and realized limits: the two channels of
assortativity}
\label{sec:realized}

Rescale displacements by $U_{ij,n}=\delta_{\T^d}(X_i,X_j)/\eps_n$ and
define the limiting link
$q(w,w';u)=1-\exp\{-\tfrac{\lambda}{\rho}ww'k(u)\}$, the mean degree
$\ell$, the endpoint intensity $\Lambda(w,x)$, and the shared-neighbor
intensity $\Gamma(w_1,w_2,x;u)$ as in \eqref{eq:ell-main-app}
(Section~S2 of the Supplement restates all definitions).  Let $\Pi$ be the edge-Palm
law of $(W_1,W_2,X,U)$ and let $\mathsf Q$ attach to it two endpoint
degrees $D_i^\star=1+\Pois$ mixtures with intensities
$\Lambda(W_i,X)-\Gamma$ and common part $\Pois(\Gamma)$.

\begin{theorem}[Joint edge-Palm, local-count, and realized-assortativity
limits]
\label{thm:realized-main}
Under the stated conditions: \textup{(i)} $np_n\to\ell$ and the
conditional law of $(W_1,W_2,X_1,U_{12,n},D_1,D_2)$ given $A_{12}=1$
converges in total variation to $\mathsf Q$; \textup{(ii)} the
edge-rooted empirical measure converges in $L^2$ against continuous
scores of quadratic degree growth; \textup{(iii)} edge, triangle, and
two-star densities converge jointly in $L^2$ to $(\ell,t,s_2)$;
\textup{(iv)} realized transitivity and the realized assortativity
coefficient converge in probability to $C=3t/s_2$ and to $r$ of
\eqref{eq:r-decomp-intro}, whose denominator is strictly positive.
\end{theorem}

Equation \eqref{eq:r-decomp-intro} is the population form of the
decomposition stated in the introduction, and the second channel
depends on the tangent displacement $U$: $\Gamma$ decreases as the
displacement approaches the connection range, because distant
endpoints share few potential neighbors.  The proof (Supplement, Section~S2)
combines a multivariate Poisson approximation for the joint
endpoint-degree law near a realized edge, an Efron--Stein
concentration argument for the realized counts, and a uniform
quadratic domination for the score class.

Three remarks concern the estimands.  First, in the certified submodel
of Section~\ref{sec:identification} the channels have closed forms,
\begin{equation}
 r_{\mathrm{int}}=\frac{2\eta(PE-S^2)}{2\eta(QE-S^2)+SE},
 \qquad
 r_{\mathrm{overlap}}=\frac{\tfrac34\,TE}{2\eta(QE-S^2)+SE},
 \qquad r=r_{\mathrm{int}}+r_{\mathrm{overlap}},
 \label{eq:channels}
\end{equation}
so channel estimates and delta-method intervals follow from any
estimate of $(\lambda,\eta,\psi)$.  Second, $r_{\mathrm{int}}$ is the
intensity-covariance channel from \emph{all} sources, not from
$\psi$ alone: even at $\psi=0$, where positions are homogeneous, the
exponential link induces a small negative endpoint-intensity
covariance under edge sampling ($r_{\mathrm{int}}=\chanSortNull$ at
the design center, against $+\chanSortCenter$ at $\psi=0.0835$, where
$r_{\mathrm{overlap}}=\chanOverlapCenter$); under the inhomogeneous
fixed-position designs of Section~\ref{sec:application} the same
channel is large and positive at $\psi=0$, driven entirely by density.
For this reason $r_{\mathrm{int}}$ is never interpreted as sorting.
Third, the sorting contribution of the dependence parameter is the
compensated increment
\begin{equation}
 \Delta_{\mathrm{sort}}(\psi_1)
 = r\{\lambda(\psi_1),\eta(\psi_1),\psi_1\}
 - r\{\lambda(0),\eta(0),0\},
 \label{eq:dsort}
\end{equation}
where $\psi\mapsto(\lambda(\psi),\eta(\psi))$ is the unique profile
holding $(\ell,C)$ at their observed values
(unique by Theorem~\ref{thm:identification-main}(i)).  By
construction $\Delta_{\mathrm{sort}}(0)=0$, and by
Theorem~\ref{thm:identification-main}(iii) it is strictly increasing,
so it is a well-defined, sign-unambiguous measure of what the
dependence adds to assortativity once intensity and range are
recalibrated.  The geometry-only baseline of
Section~\ref{sec:application} is exactly the second term of
\eqref{eq:dsort} evaluated under the conditional design.

\section{Certified global identification}
\label{sec:identification}

Specialize to $d=1$, $\rho=1$, $W=V\sim\mathrm{Unif}(0,1)$, $m=1$,
$b=\sqrt\psi$, hard-range kernel $k_\eta(z)=\one\{|z|\le\eta\}$, and
exact link $H_\lambda(u,v)=1-e^{-\lambda uv}$.  The three observables
reduce exactly to
\begin{equation}
 \ell=2\eta E,\qquad
 C=\frac{3T}{4S},\qquad
 r=\frac{2\eta(PE-S^2)+\tfrac34TE}{2\eta(QE-S^2)+SE},
 \label{eq:M3-main}
\end{equation}
with $E,S,T,P,Q$ the five mark-motif integrals (polynomial in $\psi$,
analytic in $\lambda$; Supplement, Section~S3).  Write
$M_3(\lambda,\eta,\psi)=(\ell,C,r)$.

\begin{theorem}[Certified global identification]
\label{thm:identification-main}
Let $B=[1,5]\times[5,15]\times[0,0.15]$.  Then:
\begin{enumerate}[label=\textup{(\roman*)}]
\item $M_3$ is \emph{injective on $B$}, including the boundary
$\psi=0$: the population mean degree, transitivity, and assortativity
determine $(\lambda,\eta,\psi)$ uniquely on $B$.  The $\ell$-image of
$B$ spans $[2.034,17.042]$.
\item Uniformly on $B$, $\det DM_3<-1/64$ and
$\sigma_{\min}(DM_3)\ge1/100$; $M_3$ is a real-analytic diffeomorphism
of $\operatorname{int}B$ onto its image with locally Lipschitz inverse.
\item Along any profile holding $(\ell,C)$ fixed, $dr/d\psi>0$
throughout $B$: after compensating intensity and range, assortativity
is strictly increasing in the dependence strength.
\end{enumerate}
\end{theorem}

\begin{figure}[t]
\centering
\includegraphics[width=.82\linewidth]{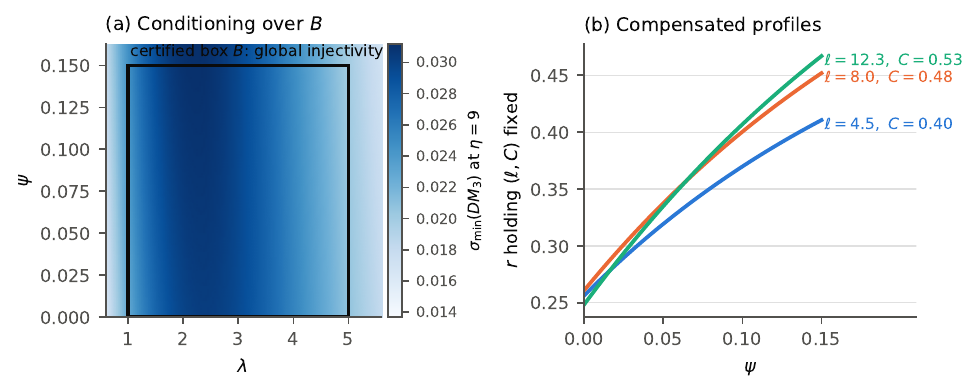}
\caption{Identification.  Left: $\sigma_{\min}(DM_3)$ over
$(\lambda,\psi)$ at $\eta=9$ with the certified region $B$.  Right:
$r$ against $\psi$ along profiles holding $(\ell,C)$ fixed; the
positive slope is certified in
Theorem~\ref{thm:identification-main}(iii).}
\label{fig:identification}
\end{figure}

The proof has two parts.  The analytical part is a
monotone-triangular argument (Supplement, Section~S3): $C$ is free of $\eta$;
given $\ell$, the range solves $2\eta E=\ell$; and on the
$(\lambda,\psi)$ rectangle certified sign conditions force the
$C$-level curves to be graphs along which the reduced assortativity is
strictly monotone.  The computational part is an exact-rational,
outward-rounded interval certificate; its implementation (the
coefficient-level composite arithmetic that eliminates interval
dependency, the adaptive partition, and the precision and truncation
choices) is documented in Section~S3 of the Supplement together with a frozen
verifier output, following standard practice in validated numerics
\citep{tucker2002,moore2009}.  The statistical consequence is that the
inverse problem underlying channel estimation has exactly one
solution on $B$, with a quantitative lower bound on conditioning, so
the channel estimates of Sections~\ref{sec:simulations}
and~\ref{sec:application} do not depend on a choice among multiple
roots.  Part (iii) additionally makes the sorting contribution
\eqref{eq:dsort} well defined: the compensating profile exists and is
unique, and $\Delta_{\mathrm{sort}}$ is strictly increasing with
$\Delta_{\mathrm{sort}}(0)=0$, so a reported sorting contribution is
never an artifact of a solution branch or of the parametrization.  Figure~\ref{fig:identification} maps $\sigma_{\min}(DM_3)$
over the region and shows the compensated profiles of part (iii).

\section{Verified inference and the independence boundary}
\label{sec:inference}

Let $Y_n=\tfrac1n(2E_n,T_n,S_{2,n},H_{3,n},H_{11,n})^{\trans}$ collect
the five connected-motif moments.  Three ingredients allow every
model-side regularity condition for inference in the certified
submodel to be verified (Supplement, Section~S4): closed forms for
the limit means,
\[
 \mu(\theta)=\Bigl(2\eta E,\ \tfrac12\eta^2T,\ 2\eta^2S,\
 8\eta^3Q+12\eta^2S+2\eta E,\
 8\eta^3P+3\eta^2T+8\eta^2S+2\eta E\Bigr)^{\trans},
\]
a uniform bias bound
$\sup_{\theta\in B}\|\E_\theta Y_n-\mu(\theta)\|\le C_B/n$ (and the
same for first derivatives), and the certificate of
Theorem~\ref{thm:identification-main}, which supplies both the unique
population root and the full rank of $G=D\mu(\theta_0)$.

\begin{theorem}[CLT and feasible GMM with verified conditions]
\label{thm:inference-main}
Let $\theta_0$ be interior to a compact
$\Theta_b\subset\operatorname{int}B$ with $\psi_0>0$, and let
$\widehat\theta_n$ be the GMM estimator computed with a certified
evaluator of $\mu$ whose uniform error is $o(n^{-1/2})$, a simulated
weight matrix based on $R_n\to\infty$ replicates at a consistent
pilot, and approximate minimization with optimization error
$o_p(n^{-1})$.  Then
$\sqrt n\{Y_n-\mu_n(\theta_0)\}\Rightarrow N\{0,\Omega(\theta_0)\}$
with $\Omega(\theta_0)\succ0$; the limit center $\mu(\theta_0)$ may
replace $\mu_n(\theta_0)$; the weight matrix is consistent; and
$\sqrt n(\widehat\theta_n-\theta_0)\Rightarrow
N\{0,(G^{\trans}\Omega^{-1}G)^{-1}\}$.  Every model-side regularity
condition (bias, covariance nondegeneracy, moment convergence,
identification, and rank) is verified for the submodel; the
evaluator, simulation, and optimization requirements are explicit
conditions on the implementation (Theorem~S4.3 of the Supplement).
\end{theorem}

\begin{theorem}[Boundary inference at $\psi=0$]
\label{thm:boundary-main}
Let $\vartheta_b=(\lambda_0,\eta_0,0)$ with
$(\lambda_0,\eta_0)\in(1,5)\times(5,15)$, let $\widehat\psi$ be the
GMM estimator constrained to $\psi\ge0$ and computed under the
implementation conditions of Theorem~\ref{thm:inference-main}, and
let $D_n$ be the profiled distance statistic for $H_0:\psi=0$.  Then,
with $I_{\mathrm{eff}}>0$ the efficient information for $\psi$ in the
moment system used,
\[
 \sqrt n\,\widehat\psi\Rightarrow
 \max\{0,N(0,I_{\mathrm{eff}}^{-1})\},
 \qquad
 n^{1/4}\,\widehat b=n^{1/4}\sqrt{\widehat\psi}\Rightarrow
 \bigl[\max\{0,N(0,I_{\mathrm{eff}}^{-1})\}\bigr]^{1/2},
 \qquad
 D_n\Rightarrow\tfrac12\chi^2_0+\tfrac12\chi^2_1 .
\]
\end{theorem}

Section~S5 of the Supplement contains five extensions that connect the certified
submodel to the application.  First, the homogeneous two-dimensional
hard-disc submodel reduces exactly to the one-dimensional moment map
under $2\eta\mapsto\pi\eta^2$ and
$\tfrac34\mapsto c_2=1-\tfrac{3\sqrt3}{4\pi}$ (the geometric-graph
clustering constant), so the same interval method certifies global
injectivity and $\sigma_{\min}\ge1/100$ for $d=2$ on the
corresponding region, with $\psi=0$ included (Theorem~S5.6); the
certified identification theory therefore covers the planar geometry
of the application, and the remaining gap is positional
inhomogeneity.  Second, a conditional-design theorem: for arbitrary
fixed position configurations with bounded local geometry, the motif
counts obey a dependency-graph CLT, and the simulated-moment
estimator used in Section~\ref{sec:application} is consistent and
asymptotically normal under explicit conditions; the identification
and nondegeneracy conditions are checked per dataset by the reported
Jacobian diagnostic.  Third, a multi-harmonic identification result:
if the mark law carries any finite set of non-resonant harmonics with
total energy $\psi_{\mathrm{tot}}$, then mean degree and transitivity
equal the one-harmonic maps at $\psi_{\mathrm{tot}}$ exactly, and the
one-harmonic inversion recovers $\psi_{\mathrm{tot}}$ to within $3\%$
at the worst admissible harmonic profile; this extends the
two-harmonic simulation result of Section~\ref{sec:simulations} to a
theorem.  Fourth, because marks are observable in the intended
applications, Section~S5 derives the joint law of the graph-only and direct
estimators and a Hausman-type specification test
($T_n\Rightarrow\chi^2_1$ in the interior); its size and power are
reported in Section~\ref{sec:simulations}.  Fifth, a
signed-amplitude result for irregular designs: for fixed positions,
the expected motif counts are exactly quadratic in the signed
amplitude $b$, the linear coefficients are explicit design
functionals weighted by a strictly positive constant
$S_a(\lambda)$, and whenever the mean score elevation of
within-range pairs is nonzero, which holds overwhelmingly in all
eleven fitted city designs, the sign of $b$ is identified at first
order, in exact contrast to the homogeneous model
(Proposition~S5.4 of the Supplement).

The boundary theorem uses the fact that the certificate includes
$\psi=0$: the null is a member of the certified region, so
identification under the null, the rank of the nuisance Jacobian, and
$I_{\mathrm{eff}}>0$ are consequences of the certificate and are not
separately assumed (the verification is Proposition~S4.7 of the
Supplement).  The amplitude rate $n^{1/4}$ and the chi-bar-square law
are the standard boundary phenomena
\citep{chernoff1954,selfliang1987,andrews1999,andrews2001}.  The
theorem as stated is pointwise, and the Supplement upgrades it:
under explicit uniform weight and optimization conditions, the
boundary laws hold \emph{uniformly over the entire certified null
face} $[1,5]\times[5,15]\times\{0\}$, with the required uniform
information bound
$\inf I_{\mathrm{eff}}>0$ discharged by a further exact-rational
interval certificate instead of being assumed
(Theorem~S4.9 and the accompanying certificate; the certified
constant is conservative by construction, and numerically the
efficient information along the face lies between $\ieffLo$ and
$\ieffHi$).  The sign
of $b$ is not identified in the homogeneous model
(Lemma~\ref{lem:equivalence-main}), so no test of direction is
attempted there; under the irregular designs of the application the
sign enters at first order (Proposition~S5.4).  The theorem covers
the homogeneous random-design submodel; the conditional
fixed-position designs of Section~\ref{sec:application} are outside
its scope, and the boundary test applied there is therefore
calibrated by a parametric bootstrap whose size is verified by
simulation in Section~\ref{sec:simulations}.

\section{Simulation study}
\label{sec:simulations}

All simulations use the certified submodel and the deterministic
certified evaluator, whose truncation error (below $10^{-15}$) is
negligible relative to Monte Carlo error at the simulated sample
sizes; replication code and frozen results are packaged with the
paper.  The study covers the certified region: eight parameter
configurations spanning weak-information corners, near-boundary and
boundary dependence, and the edges of $B$; per-design intensity
shares, RMSE, coverage, interval lengths, and solver success rates
are tabulated in Section~S5 of the Supplement, and coverage of the
intensity share is nominal across the grid.

\paragraph{Coverage: marginal and simultaneous.}
Table~\ref{tab:coverage} reports \emph{marginal} coverage of nominal
$95\%$ delta-method intervals and bias at the center design
$\theta_0=(3,9,0.0835)$ for $n=10^3,10^4,10^5$ (\repSmall,
\repSmall, and \repLarge\ replications).  Marginal coverage is
within Monte Carlo error of nominal from $n=10^4$; at $n=10^3$ the
$\psi$ interval undercovers mildly (\covPsiSmallPct), and the Newton
solver converges in $91.5\%$ of replicates (failures are retained in
all reported rates).  Simultaneous inference uses the Wald ellipse at
the $\chi^2_3$ level: its coverage is \ellipseCovThousand,
\ellipseCovTenk, and \ellipseCovHundk\ across the three sample sizes;
the naive ``all three marginal intervals cover'' event, which targets
no nominal level, occurs with probability \jointNaiveThousand,
\jointNaiveTenk, \jointNaiveHundk, consistent with three unadjusted
intervals.  A fully feasible arm (own pilot, simulated covariance,
$R=150$) gives marginal coverage \feasCov\ at $n=10^4$ against the
oracle arm's \covLamTenk/\covEtaTenk/\covPsiTenk; with $100$
replicates the difference is not statistically resolvable, and both
numbers are reported.  Measured moment biases
scale as $c/n$ ($n\times$bias of $\widehat r_n$: \biasRthousand\ at
$n=10^3$, \biasRtenk\ at $n=10^4$), consistent with the uniform bias
lemma.

\begin{table}[t]
\centering
\caption{Center-design \emph{marginal} coverage of nominal $95\%$
intervals and bias.  Monte Carlo standard error for coverage is about
$1.1\%$ at $n\le10^4$ and $2.2\%$ at $n=10^5$.}
\label{tab:coverage}
\begin{tabular}{@{}lccccccc@{}}
\toprule
 & \multicolumn{3}{c}{Marginal coverage (\%)} & \multicolumn{3}{c}{Bias} \\
\cmidrule(lr){2-4}\cmidrule(lr){5-7}
$n$ & $\lambda$ & $\eta$ & $\psi$ & $\lambda$ & $\eta$ & $\psi$ \\
\midrule
$10^3$ & \covLamThousand & \covEtaThousand & \covPsiThousand
       & \biasLamThousand & \biasEtaThousand & \biasPsiThousand \\
$10^4$ & \covLamTenk & \covEtaTenk & \covPsiTenk
       & \biasLamTenk & \biasEtaTenk & \biasPsiTenk \\
$10^5$ & \covLamHundk & \covEtaHundk & \covPsiHundk
       & \biasLamHundk & \biasEtaHundk & \biasPsiHundk \\
\bottomrule
\end{tabular}
\end{table}

\paragraph{Exactly identified versus overidentified GMM.}
At the center design, the five-motif optimally weighted GMM of
Theorem~\ref{thm:inference-main} (deterministic evaluator, oracle
weighting from a $1{,}200$-replicate bank) improves on the
three-moment exactly identified system: RMSE falls from
\rmseThreePsi\ to \rmseFivePsi\ for $\psi$ (with similar gains for
$\lambda$ and $\eta$), coverage stays nominal
(\covFivePct), and the $\chi^2_2$ overidentification test has size
\jTestPct\ at the $5\%$ level.

\paragraph{Same $(\ell,r)$, different mechanisms.}
Figure~\ref{fig:thesis} reports the mechanism-separation experiment:
two parameter points inside $B$ matched to the same population mean
degree and assortativity, with intensity shares $\pairShareA$ and
$\pairShareB$ and sorting contributions
$\Delta_{\mathrm{sort}}=\pairDeltaA$ and $\pairDeltaB$.  Across $400$
replicates at $n=10^4$ the realized $\widehat r$ distributions are
indistinguishable while the realized $\widehat C$ distributions are
disjoint, and the full three-moment method recovers the two intensity
shares with means $\pairShareHatA$ and $\pairShareHatB$.  Mean degree
and assortativity alone do not separate the mechanisms; adding
transitivity does, which motivates identification from $(\ell,C,r)$.

\begin{figure}[t]
\centering
\includegraphics[width=.85\linewidth]{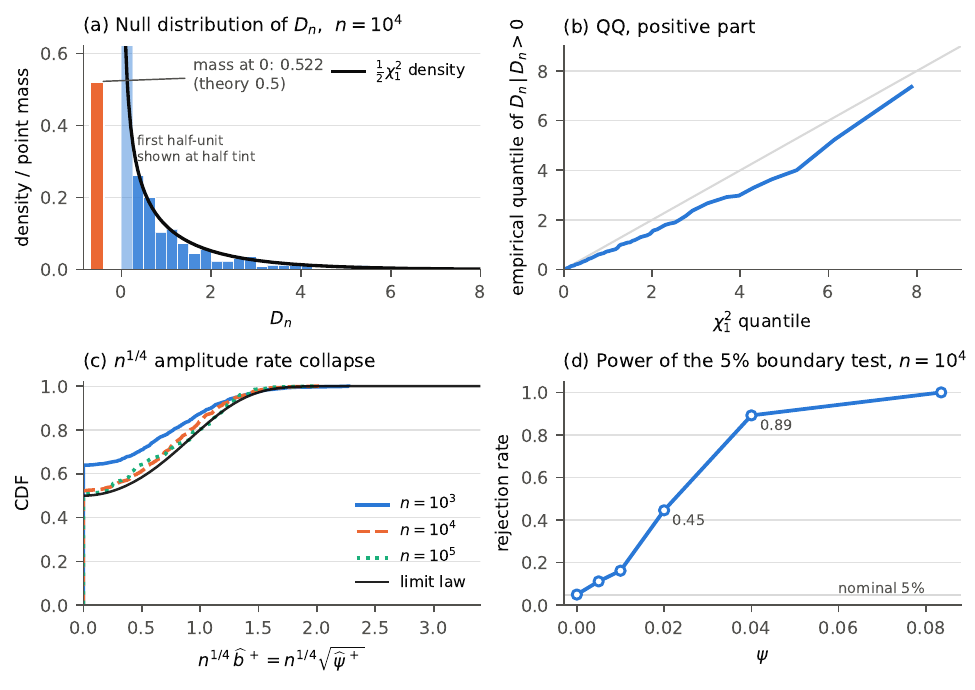}
\caption{Boundary behavior at $\psi=0$ and power.  Left: distribution
of $D_n$ under the null at $n=10^4$ against
$\tfrac12\chi^2_0+\tfrac12\chi^2_1$.  Middle: QQ plot of the positive
part of $D_n$ against $\chi^2_1$.  Right: power of the $5\%$ test at
$n=10^4$.}
\label{fig:boundary}
\end{figure}

\paragraph{Boundary size across nuisances.}
The $5\%$ boundary test has empirical size \sizeNuisA, \sizeNuisB, and
\sizeNuisC\ at nuisance configurations $(\lambda_0,\eta_0)=(2,7)$,
$(3,9)$, $(4,12)$, with null mass at zero $0.54$, $0.54$, and
\massNuisC\ converging to the theoretical $\tfrac12$ (null mass \nullMassTenkPct\ and size \sizeTenkPct\ in the larger center-nuisance run; $I_{\mathrm{eff}}=\Ieff$;
Figure~\ref{fig:boundary}); power at the center nuisance rises to
\powFourPct\ at $\psi=0.04$.

\paragraph{Boundary law and power.}
Under $\psi=0$ the fraction of replications with $D_n=0$ is
\nullMassTenkPct\ at $n=10^4$ (theory: $50\%$), the rejection rate of
the $5\%$ test is \sizeTenkPct, and the positive part of $D_n$ matches
the $\chi^2_1$ quantiles (Figure~\ref{fig:boundary}); the distribution
of $n^{1/4}\sqrt{\widehat\psi^{\,+}}$ is stable in $n$ across
$10^3$--$10^5$, consistent with the $n^{1/4}$ amplitude rate.
Power at $n=10^4$ rises from size \sizeTenkPct\ at $\psi=0$ to
\powTwoPct\ at $\psi=0.02$, \powFourPct\ at $\psi=0.04$, and
\powFullPct\ at $\psi=0.0835$.

\paragraph{Specification test and heavier misspecification.}
The graph-versus-direct test of Section~S5 of the Supplement has empirical size
\hausSize\ at the $5\%$ level at the center design.  Under a convex
(heavier-tailed) transform of the popularity marks, which biases the
graph-only estimate to \heavyPsi\ against the truth $0.0835$ while
leaving the direct estimate unbiased (the transform does not change
the observed ranks), the test rejects with power \hausHeavy.  Under
$2\%$ and $5\%$ admixed uniform long-range ties the graph-only
estimate changes little (\lrTwoPsi\ and \lrFivePsi) and the test
stays at or below nominal size (\hausLRtwo, \hausLRfive): the test
rejects under contaminations that bias the estimand and does not
reject under contaminations that leave it unchanged.  Because the
data of Section~\ref{sec:application} sit at the $\psi=0$ boundary,
where the $\chi^2_1$ calibration does not apply, the test is also
implemented with a parametric-bootstrap null calibrated at the
boundary fit: its empirical size at $\psi_0=0$ is \btSize\ (nominal
$5\%$), and this version is applied to the data.  A
collapsed-displacement baseline (inverting the moment map with the
overlap channel deleted, as a latent-position analysis that ignores
microscopic displacement implicitly does) misattributes the entire
overlap term to the intensity channel and biases $\widehat\psi$
upward by at least $95\%$ at the center design, moving it outside the
admissible range; the tangent coordinate is therefore necessary for
correct channel estimation.

\paragraph{Misspecification.}
Generating with a triangular or truncated-Gaussian kernel of the same
integrated range and fitting the hard-window model moves
$\widehat\psi$ by \missKernelSummary; generating with a two-harmonic
mark law of total strength $\psi_{\mathrm{tot}}=0.0835$ and fitting the
one-harmonic model recovers \missTwoharmPsi\ on average: the estimand
behaves as the total $L^2$ dependence, as predicted by the moment
structure, because cross-harmonic terms cancel in all motifs up to the
four-vertex level.

\section{Application to location-based social networks}
\label{sec:application}

The application analyzes two location-based social networks,
Brightkite and Gowalla; their public snapshots \citep{cho2011}
contain friendship graphs and $4.5$ and $6.4$ million check-ins
respectively, which give each user both a home location
(median check-in coordinate) and an observed activity level, whose
within-sample rank serves as the popularity mark $V$.  Because both
coordinates of the latent pair are observed, the dependence strength
can be estimated in two independent ways: directly from the
$(X_i,V_i)$ pairs (the direct estimator
$\widehat\psi_{\mathrm{dir}}$), and from the graph moments alone with
marks hidden (the graph-only estimator $\widehat\psi_{\mathrm{gr}}$).
Agreement of the two estimates is a check of the graph-only method
that purely latent settings cannot provide.  The design conditions on
the observed home configuration, since uniform positions are
inappropriate for cities: the spatial score is standardized log home
density, the dependence family is
\eqref{eq:harmonic-density-main} in that score, edges follow the exact
link within a hard radius $R$, and $(\lambda,R,\psi)$ is fit by
simulated moment matching of $(\widehat\ell,\widehat C,\widehat r)$
with positions held fixed.  The decomposition question is addressed
by the geometry-only baseline: the model is refit with $\psi=0$,
recalibrating $(\lambda,R)$ to the observed mean degree and
transitivity, and the resulting assortativity $r_{\mathrm{geo}}$ is
the value that geometric overlap and density produce alone.

Because attribution is meaningful only under a model that fits,
goodness of fit is reported before anything else.  The fitting
protocol matches $(\widehat\ell,\widehat C)$ within tolerance where
the constraints $R\ge0.1$ km and $R\le R_{\max}$ permit (radius rails
are flagged, and mean degree is left unmatched there), and then
searches the constrained dependence amplitude $b\in[0,0.55]$ for the
assortativity moment; the residual
$\widehat r-r_{\mathrm{model}}(\widehat\theta)$ is left unmatched
whenever the constrained search fails, and that residual, standardized
by the bootstrap standard deviation of realized model assortativity at
$\widehat\theta$ (a goodness-of-fit $z$), is the specification
diagnostic.  Table~\ref{tab:gof} reports, for every city, the fitted
and observed assortativity, the residual and its $z$, which parameter
bounds the optimizer reached, the share of observed ties longer than
the fitted connection radius, and the model-implied sorting
contribution $\Delta_{\mathrm{sort}}$ at the signed direct estimate
$\widehat b_{\mathrm{dir}}$.  Channel and baseline numbers are
interpreted in the light of this table throughout: in these data they
diagnose the model class; they do not describe the cities.

\begin{table}[t]
\centering
\caption{Goodness of fit before attribution, all eleven cities.
$r_{\mathrm{mod}}$: realized assortativity of the fitted model at
$\widehat\theta$.  resid: $\widehat r-r_{\mathrm{mod}}$, with $z$ its
ratio to the bootstrap standard deviation of realized model
assortativity (\bootDraws\ draws).  rail: parameters at optimizer
bounds ($b\in[0,0.55]$, $R\ge0.1$ km).  $\%{>}\widehat R$: share of
observed ties longer than the fitted radius, which the fitted model
assigns probability zero.  $\Delta_{\mathrm{sort}}$: model-implied
change in assortativity from dependence of the signed directly
estimated size, after recalibrating $(\lambda,R)$ to
$(\widehat\ell,\widehat C)$.  Every fit fails in the assortativity
moment; the class is rejected in every city.}
\label{tab:gof}
\scriptsize
\gofTable
\end{table}

\paragraph{Known-truth validation.}
The pipeline is first applied to a benchmark generated from the
conditional model itself at Brightkite-like scale (clustered
multi-core home configurations,
$(\lambda,R,\psi)=(2,\,2\text{ km},\,0.0625)$, $n=\dmN$).  The
graph-only fit matches all three moments and recovers
$\widehat\lambda=\dmLamHat$, $\widehat R=\dmRHat$ km (within about two
percent) and $\widehat\psi_{\mathrm{gr}}=\dmPsiGraph$ with bootstrap
$95\%$ band \dmPsiCI\ covering the truth; the direct estimate
($\dmPsiDirect\pm\dmPsiDirectSE$) is more precise and consistent;
the geometry-only baseline reproduces \dmShare\ of the benchmark's
assortativity ($r_{\mathrm{geo}}=\dmRgeo$ against \dmR); parameters
transfer to a held-out region within \dmTransferPct; the
band shrinks to \dmTwelveCI\ at $n=12{,}000$, consistent with the
root-$n$ rate; and the finite-difference Jacobian diagnostic is
well conditioned ($\sigma_{\min}=\bkJacSmin$, condition number
\bkJacCond).  The benchmark fixes the reading of the city fits below:
when the class is true, the fit matches all three moments and the
interior truth is recovered.  In the data the assortativity moment is
missed by half a unit in every city; a boundary estimate produced
under such a misfit is a consequence of the rejection, not a
measurement of independence, and it is never presented as one below.

\paragraph{Brightkite: five cities.}
The analysis starts from the full public friendship list
($214{,}078$ edges) and a per-user homes summary of the $4.5$ million
check-ins (median coordinate; $\ge5$ valid check-ins), and studies the
subgraphs induced by five prespecified metropolitan regions ($2{,}894$
edges in SF Bay); all conclusions concern these metropolitan
subgraphs.  Table~\ref{tab:cities} reports, for each region, the
observed moments, Spearman (rank) assortativity, the configuration
null, both estimates of $\psi$, and the geometry-only baseline
$r_{\mathrm{geo}}$; Table~\ref{tab:gof} carries the corresponding
goodness-of-fit rows.

\begin{table}[t]
\centering
\caption{Brightkite cities.  $\widehat r_{\mathrm{rank}}$: Spearman
assortativity.  $r_{\mathrm{cfg}}, C_{\mathrm{cfg}}$:
degree-preserving configuration null (double-edge swaps).
$\widehat\psi_{\mathrm{dir}}$: direct estimate from
$(\text{location},\text{rank})$ pairs.
$\widehat\psi_{\mathrm{gr}}$: graph-only estimate with bootstrap
$95\%$ band (\bootDraws\ draws); boundary values arise under the
misfit documented in Table~\ref{tab:gof} and are diagnostic, not
measurements of independence.  $r_{\mathrm{geo}}$: assortativity of
the geometry-only baseline recalibrated to the observed
$(\widehat\ell,\widehat C)$.  $\sigma_{\min}$: smallest singular
value of the conditional-design moment-map Jacobian at the fit, the
dataset-specific identification diagnostic; NYC is weakly identified
(see text).}
\label{tab:cities}
\scriptsize
\cityTable
\end{table}

Three findings hold in every Brightkite city.  First, \emph{the model
class is rejected, and the rejection is the central empirical fact}.
Every fit matches transitivity, matches mean degree except where the
radius itself rails (Denver at $R_{\max}$, NYC at the lower bound),
and misses assortativity everywhere: $r_{\mathrm{geo}}$ lies between
\bkRgeoLo\ and \bkRgeoHi\ against observed values between \bkRLo\ and
\bkRHi, the residuals and their goodness-of-fit $z$ values are in
Table~\ref{tab:gof}, and no $\psi\ge0$ closes the gap of \bkGap\ at
the SF fit, because assortativity is increasing in $\psi$ after
compensation (Theorem~\ref{thm:identification-main}(iii)).  The
constrained amplitude therefore rails at $b=0$ in four cities: these
boundary values are consequences of the misfit, not measurements of
independence.  The fitted radii tell the same story from the edge
lengths: \bkFracLo\ to \bkFracHi\ percent of observed ties are longer
than $\widehat R$ (median within-region tie \bkMedKm\ km; only
\bkPctOneKm\ within $1$ km; $\widehat R=\bkRhatSF$ km in SF), so the
fitted model assigns probability zero to most of the edges it is
meant to explain.  Clustered homes combined with any short-range
kernel produce strongly positive assortativity at the observed
clustering levels, through intensity covariance (model channel
$r_{\mathrm{int}}=\bkChSort$ at the SF fit, driven by density) and
overlap ($\bkChOv$); the data refuse the prediction.  The fifth city,
NYC, adds an identification failure on top of the misfit: the
finite-difference Jacobian of the conditional moment map is nearly
singular there ($\sigma_{\min}=\bkNycSmin$, condition number
$\bkNycCond$, against $\bkSminLo$--$\bkSminHi$ elsewhere), the
optimizer terminates with
both the amplitude at its upper bound and the radius at its lower
bound, and the bootstrap band spans essentially the admissible range;
the graph moments are uninformative in NYC and the diagnostic
indicates this without reference to the estimate.

Second, \emph{the direct estimator, which does not use the graph,
finds no evidence of the specified dependence}:
$\widehat\psi_{\mathrm{dir}}$ ranges from \bkPsiDirLo\ to \bkPsiDirHi\
with standard errors of the same size (SF Bay:
$\widehat b_{\mathrm{dir}}=\bkBDirect\pm\bkBDirectSE$).  The claim is
deliberately scoped: what is tested is the one-harmonic dependence of
activity rank on the chosen spatial score, not independence of
activity from location in general.  Under the fitted geometry, taken
at face value, dependence of the directly estimated size would change
model assortativity by at most $\bkDsortMax$ in absolute value
(Table~\ref{tab:gof}), so even if the class fit, sorting of the
detected size would be negligible for mixing.

Third, \emph{the configuration null identifies which observed
features require explanation}: degree-preserving rewiring reproduces
assortativity of the observed sign and magnitude (SF:
$r_{\mathrm{cfg}}=\bkCfgR$ against observed $\bkR$), so the mild
disassortativity is consistent with a degree-sequence phenomenon, though the null
overshoots it, leaving a positive residual of
$+0.11$ in SF (maximum degree \bkMaxDeg; the top $1\%$ of users
carry \bkHubShare\ of edge ends), but the rewiring collapses
transitivity
($C_{\mathrm{cfg}}=\bkCfgC$ against observed \bkC); a complementary
null that matches the observed edge-length distribution with random
pairs also produces near-zero transitivity
($C_{\mathrm{dist}}=\bkDistC$ in SF) and near-zero assortativity
($r=\bkDistR$).  The observed triangles are therefore explained
neither by the degree sequence nor, by the first finding, by spatial
overlap; a mechanism outside the fitted class, such as non-spatial
triadic closure, is required, and the present analysis does not
identify it.  Pearson and Spearman assortativity disagree in sign
(SF: $\bkR$ versus \bkSpearman), consistent with the sensitivity of
the Pearson coefficient to hubs \citep{kaufmann2026}; under either
version the gap to $r_{\mathrm{geo}}$ remains, and per the scope
fixed in Section~\ref{sec:intro}, the decomposition claims concern
the Pearson coefficient only.

\begin{figure}[t]
\centering
\includegraphics[width=.46\linewidth]{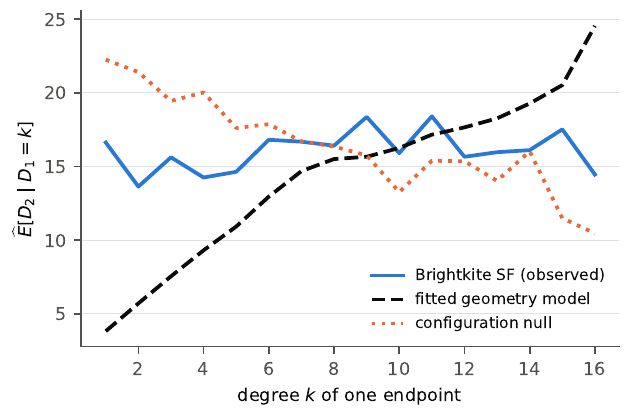}
\caption{Conditional-degree curves for SF Bay: observed versus the
fitted geometry model and the degree-preserving configuration null.}
\label{fig:annd}
\end{figure}

Table~\ref{tab:battery} summarizes the null comparisons; three
further checks follow.  A \emph{fitted} latent-distance model
(logistic distance-decay likelihood estimated on all SF dyads,
$\widehat\beta=\ldBeta$ per km) reproduces mean degree by
construction but produces transitivity \ldC\ against the observed
\bkC: a fitted spatial likelihood, like the randomization nulls,
fails to reproduce the observed transitivity.  The conditional-degree
curves $\widehat E[D_2\mid D_1=k]$ (Figure~\ref{fig:annd}) show that
the comparison does not depend on a scalar summary: the geometry
model implies a steeply increasing profile (from about $4$ at $k=1$
to about $25$ at $k=16$), the configuration null implies a decreasing
profile (the hub effect), and the observed curve is nearly flat and
lies between the two, differing from the geometric profile at every
degree.  The boundary-calibrated two-view test comparing
$\widehat\psi_{\mathrm{dir}}$ and $\widehat\psi_{\mathrm{gr}}$ does
not reject in SF Bay ($p=\btSFp$); the test addresses only the
dependence parameter, and since the class is already rejected through
the unmatched assortativity moment, the non-rejection is reported for
completeness and not as validation of the model.

\begin{table}[t]
\centering
\caption{Implied transitivity and assortativity of each candidate
mechanism, SF Bay.}
\label{tab:battery}
\begin{tabular}{@{}lcc@{}}
\toprule
 & transitivity $C$ & assortativity $r$ \\
\midrule
degree sequence (config.\ null) & \bkCfgC & $\bkCfgR$ \\
edge-length profile (distance null) & \bkDistC & $\bkDistR$ \\
fitted latent-distance likelihood & \ldC & $\ldR$ \\
geometry model matched to $(\ell,C)$ & \bkC & $+\bkRgeo$ \\
\textbf{observed} & \textbf{\bkC} & $\mathbf{\bkR}$ \\
\bottomrule
\end{tabular}
\end{table}

\paragraph{Gowalla: six cities.}
Replication on the second platform ($950{,}327$ friendships, $6.4$
million check-ins, six prespecified metropolitan subgraphs) repeats
the rejection and adds two cases.  Table~\ref{tab:gowalla} reports
the six cities; their goodness-of-fit rows are in
Table~\ref{tab:gof}.  Every fit again matches
$(\widehat\ell,\widehat C)$, misses assortativity ($r_{\mathrm{geo}}$
of \gwRgeoLo\ to \gwRgeoHi\ against observed \gwRLo\ to \gwRHi),
rails at $b=0$, and places \gwFracLo\ to \gwFracHi\ percent of
observed ties beyond the fitted radius; configuration nulls again
collapse transitivity.  Where Brightkite's observed assortativity was
close to its configuration null, Gowalla's lies systematically above
it in all six cities, most clearly in Stockholm (observed $\gwStoR$
against a null of $\gwStoCfgR$ and a geometry-only prediction of
$\gwStoRgeo$): a positive mesoscale residual, consistent with community
structure, that is explained neither by the degree sequence, nor by
the specified dependence ($\widehat\psi_{\mathrm{dir}}\approx0$
there), nor by spatial overlap, and that remains unexplained by every
mechanism fitted in this paper.  Pearson and Spearman assortativity
again diverge (Austin: $\gwAusR$ versus \gwAusSpear).  In Austin the
direct estimator detects small but statistically significant
popularity--location dependence,
$\widehat b_{\mathrm{dir}}=\gwAusBDir\ (\pm\gwAusBDirSE)$ and
$\widehat\psi_{\mathrm{dir}}=\gwAusPsiD\ (\pm\gwAusPsiDSE)$, the only
nonzero detection among the eleven cities (a two-sided $z$ of about
$4.0$; with a Bonferroni correction across the eleven direct tests
the detection stands at $p\approx10^{-3}$).  The estimate is
negative: the most active users concentrate away from the densest
home areas.  The sign is invisible to the homogeneous adjacency law
(Lemma~\ref{lem:equivalence-main}); under the fitted irregular
design it enters the graph law at first order
(Proposition~S5.4 of the Supplement), with information far below
detection at this sample size, and the magnitude lies below the
graph-only detection threshold at $n\approx4{,}000$
(Section~\ref{sec:simulations}).  Propagated through the fitted
geometry, dependence of the detected signed size changes model
assortativity by $\gwAusDsort$ (flipping the sign of
$\widehat b_{\mathrm{dir}}$ gives $\gwAusDsortFlip$; the asymmetry
between the two is the empirical size of the linear-in-$b$ terms that
the conditional design can support, Section~S5 of the Supplement).
The direct estimator measures the dependence; the model-implied
$\Delta_{\mathrm{sort}}$ bounds its contribution to mixing.

\begin{table}[t]
\centering
\caption{Gowalla cities; columns as in Table~\ref{tab:cities}.}
\label{tab:gowalla}
\scriptsize
\gwTable
\end{table}

\paragraph{A richer class, and what it buys.}
The rejection pattern dictates the augmentation: ties an order of
magnitude beyond the fitted radius call for a second connection
scale, the configuration null's overshoot of the observed
disassortativity calls for explicit degree propensities, and
triangles unexplained by geometry call for a closure mechanism.  The
augmented class (Section~S5.5 of the Supplement) couples the original
sorting family with lognormal node propensities $\sigma$, a
two-scale kernel $\one\{d\le R\}+w e^{-d/L}$, and a one-pass triadic
closure $\tau$, and is fit by simulated moment matching to eight
targets: the original triple, the three edge-length quartiles, and
the degree coefficient of variation and $99$th percentile.  Fitted
to SF Bay and to Austin, the class removes the catastrophic
single-moment failure: the worst parametric-bootstrap $z$ falls from
$\gofZWorstOld$ (Table~\ref{tab:gof}) to $\richSFMaxZ$ (SF Bay) and
$\richAusMaxZ$ (Austin), and the fitted models reproduce the held-out
hub economy that the original class cannot represent at all (SF: top
$1\%$ edge-end share $\richSFHub$ against observed \bkHubShare;
maximum degree $\richSFMaxDeg$ against $\bkMaxDeg$).
Table~\ref{tab:rich} reports the per-target detail.  The
reconciliation is nevertheless partial, and the paper reports it as
such.  At SF Bay the fitted member overshoots mean degree and the
short-range edge-length mass at the mean level even where the $z$
gates pass, because the realization variance of the heavy-tailed
class is large; a diffuse model can pass a $z$ gate by variance
inflation, so Table~\ref{tab:rich} reports the model-realization
standard deviations alongside the $z$ values.  At Austin, with the
dependence amplitude pinned to the signed direct estimate, which
closes both optimizer escape routes, the class matches transitivity,
assortativity, all three length quartiles, and the degree
coefficient of variation ($|z|\le\richAusGoodZ$ at the mean level as
well), while mean degree remains inflated by heavy-tailed propensity
realizations ($+\richAusEllOver\%$ at the mean, $z=\richAusEllZ$)
together with the extreme degree percentile ($z=\richAusPZ$).  The
escape routes just mentioned, documented and not hidden, are
themselves diagnostics: unconstrained searches drift toward variance
inflation or toward signed amplitudes far beyond anything the direct
estimates support, which says that the residual misfit is not in the
directions the class spans.  What still binds, in both cities, is
the joint demand for short-range triangles together with long ties
and a moderate degree tail: global one-pass closure creates
triangles wherever two-paths exist, while the data concentrate them
at short range.  A distance-limited closure or an explicit community
layer is the indicated next enrichment.

Within the fitted class, the compensated counterfactuals of
Table~\ref{tab:cf} decompose the model assortativity: switching off
the propensities moves $r$ by $\cfHubSF$ (SF) and $\cfHubAus$
(Austin), the closure pass by $\cfCloSF$ and $\cfCloAus$, the
long-range scale by $\cfLrSF$ and $\cfLrAus$, and dependence of the
directly estimated signed size by $\cfBSF$ and $\cfBAus$.  The
ordering is stable across the basins explored: degree structure
first, closure and range second, sorting an order of magnitude
smaller.  Under a class rich enough to carry the main mechanisms,
the answer to the attribution question is that the
popularity--location dependence remains a minor channel.

\begin{table}[t]
\centering
\caption{The richer class at SF Bay and Austin: parametric-bootstrap
$z$ per target (obs minus model mean over model-realization sd) and
model means with sd in parentheses.  $100$ bootstrap draws at the
fitted parameters.}
\label{tab:rich}
\scriptsize
\richTable
\end{table}

\begin{table}[t]
\centering
\caption{Top: within-class compensated counterfactuals, the change
in model assortativity from switching the mechanism off (or, for the
dependence row, setting $b$ to the signed direct estimate), with
$\lambda$ recalibrated to $\widehat\ell$; simulation se
$\le0.006$.  Bottom: competing families at SF Bay, held-out edge
prediction (AUC over a $20\%$ edge holdout against sampled
non-edges) and posterior-predictive means.  DC-SBM: degree-corrected
blockmodel ($K$ by held-out AUC).  Logistic: spatial logistic with
per-node effects.  Richer: the augmented class (AUC unconditional /
conditional on observed activity marks).  Geometry: the fitted
hard-disc baseline.}
\label{tab:cf}
\scriptsize
\cfCompTable
\end{table}

\paragraph{Competing families.}
Table~\ref{tab:cf} (bottom) confronts four fitted families with the
same held-out-edges protocol and posterior-predictive checks.  The
node-conditional models dominate individual edge prediction: the
spatial logistic with node effects attains AUC $\aucLogit$ and the
degree-corrected blockmodel $\aucDcsbm$, against $\aucRichCond$ for
the richer class conditioned on observed activity and $\aucGeo$ for
geometry; this ordering is partly structural, since exchangeable
generative classes do not learn node identities from the training
graph.  The posterior-predictive columns reverse the ordering.  The
logistic model reproduces mean degree, assortativity, all three
length quartiles, and the degree profile, and collapses transitivity
($C=\compLogC$ against \bkC); the blockmodel matches assortativity
and degrees but produces $C=\compDcsbmC$ and edge lengths twice too
long; geometry matches $(\ell,C)$ by construction and fails
assortativity and lengths; the richer class is the only family
within reach of all eight targets simultaneously, at the cost
documented in Table~\ref{tab:rich}.  Clustering is reproduced only
by families with explicit local structure, and destroyed by purely
node-conditional ones; mixing, lengths, and degrees are carried by
node heterogeneity.  No fitted family passes every check.  The
demonstration is the point: these networks jointly demand local
closure, two connection scales, and degree effects, and the
decomposition of which family fails which moment is the empirical
specification of the missing model.

In summary, the short-range geometry-plus-dependence class is
decisively insufficient for these eleven metropolitan networks: it
cannot jointly reproduce their clustering, edge lengths, and
assortativity, and its constrained fits reach the $\psi=0$ boundary
as a symptom of that insufficiency.  The rejection is nevertheless
informative in a specific, certified sense.  Because assortativity is
increasing in the dependence after compensation, the observed values,
lying far below every geometric prediction, provide no evidence of
positive sorting under any member of the class; the direct estimates,
which do not use the graph, independently find none (Austin's small
negative detection changes model assortativity by less than
$\dsortBoundAll$).  What the observed pattern is consistent with is
degree heterogeneity (Brightkite, where configuration nulls carry
$r$'s sign but overshoot it and destroy $C$), plus a positive
mesoscale residual on both platforms and non-spatial triadic
closure, mechanisms the augmented class of the preceding paragraphs
begins to carry: the augmentation removes the catastrophic misfit,
confirms the mechanism ordering, and leaves a residual pattern that
itself specifies what a fitting model must add.

\section{The two-view design on a collaboration network}
\label{sec:collab}

The metropolitan analysis ends in a rejection and a null, and
neither half of that verdict is forced by the design.  This section
applies the identical protocol to a network chosen in advance as the
candidate positive case: the co-authorship network of physics
articles published 2018--2022 with a Netherlands-affiliated author,
built by the replication package's fetch-and-format script from the
OpenAlex catalog, with institution coordinates as positions and the
corpus work count as the observable productivity mark (\oaAuthors\
authors, \oaEdgesAll\ edges; construction details, the full regional
battery, and the complete fit protocols are in Section~S5.6 of the
Supplement).  Positions are institution point masses, so a
fixed-seed jitter of $0.3$\,km breaks the coordinate atoms, the
density bandwidth of the spatial score is floored at $0.5$\,km, and
the direct estimate carries standard errors clustered at the atom
level.

The two views now separate.  At the country scale ($n=\oaNLn$,
\oaNLE\ edges, $\widehat\ell=\oaNLell$, $\widehat C=\oaNLC$,
$\widehat r=\oaNLr$, Spearman $\oaNLrs$) the direct estimator gives
$\widehat b_{\mathrm{dir}}=\oaBdir$ with clustered standard error
\oaBdirSEcl, a detection at $\oaBdirZcl$ standard errors
($\widehat\psi_{\mathrm{dir}}=\oaPsiDir$): productive authors
concentrate where researcher density is high, an interior positive
signal where the metropolitan estimates gave ten nulls and one
small negative detection, and one robust to the preprocessing
constants ($\widehat b_{\mathrm{dir}}\in[\oaSensLo,\oaSensHi]$
across jitter scales, bandwidth floors, and seeds; Supplement,
Section~S5.6).  The signal vanishes
exactly where the design degenerates: in the Groningen box
\oaGroTop\% of authors share a single coordinate atom, the score
$\phi$ has no variation left, and the direct estimate collapses to
$\oaGroBdir$, an empirical instance of the effective-information
condition of Section~\ref{sec:inference}.

The graph-only route is refused at every level of the ladder
(Table~\ref{tab:oa}), each refusal with a legible residual.  The
configuration null preserves degrees and explains nothing else
($r=\oaSubCfgR$ against $\oaSubR$ observed on the fit subsample;
$C=\oaSubCfgC$ against \oaSubC).  The hard-disc class cannot
reconcile mean degree with clustering under the two-scale edge
lengths (\oaSubFracLocal\% of edges within $1$\,km,
\oaSubFracFar\% beyond $6$\,km), rails both parameters, and leaves
\oaBaseFracOut\% of observed ties outside its fitted radius.  The
augmented class with $b$ free escapes to the negative rail
($b=\oaFreeB$) and is rejected at $\max_k|z_k|=\oaFreeMaxZ$; with
$b$ pinned to the direct estimate it matches mean degree,
clustering, both short length quartiles, and the extreme degree
percentile, and is rejected on assortativity and the degree
coefficient of variation alone ($z_r=\oaFixZr$,
$z_{\mathrm{cv}}=\oaFixZcv$).  The compensated counterfactual
locates the obstruction: reaching $\widehat C=\oaSubC$ requires
closure $\tau=\oaFixTau$, and switching that closure off moves model
assortativity by $\oaFixCfClo$, so the clustering the class can only
buy with iterated closure drags $r$ far above its observed value.
Multi-author papers create triangles in single events, without the
degree-mixing side effect of iterated closure: the co-authorship
graph is a clique projection of the paper--author bipartite
structure \citep{newman2001}, and an explicit team layer is the
mechanism the residual names.  The signed channel, meanwhile,
behaves as Proposition~S5.4 predicts: model assortativity is
monotone increasing in $b$ along compensated profiles (slope
$\oaProfSlopeFix$ per unit $b$), and the sign flip
$b\mapsto-\widehat b_{\mathrm{dir}}$ lowers it by $\oaFixSignedGap$.

\begin{table}[t]
\centering
\caption{The collaboration network's graph-only ladder on the
Netherlands fit subsample.  Model $r$ and $C$ are
parametric-bootstrap means at the fitted parameters (configuration
null: rewiring mean and sd); the last columns list the moments each
fit targets and the failure flags.  Direct view, for comparison:
$\widehat b_{\mathrm{dir}}=\oaSubBdir$ (clustered se
\oaSubBdirSEcl) on the same subsample.}
\label{tab:oa}
\scriptsize
\oaTable
\end{table}

The case supplies what the metropolitan networks could not.  The
estimand is directly measurable, clearly positive, and interior,
so the two-view comparison has a nonzero target on the direct side;
and the goodness-of-fit gate demonstrates its value on real data by
refusing graph-only attribution under every fitted class, including
the augmented one, rather than laundering a misspecified fit into a
sorting estimate.  A graph-only analyst without the marks would
conclude, correctly, that no fitted class supports attribution here
and that the residual points to team structure; the marks, when
observed, answer the question the graph could not.

\section{Discussion}
\label{sec:discussion}

This paper defines the intensity-covariance and shared-neighbor
contributions to assortativity as population quantities in a class of
sparse spatial network models, defines the sorting contribution as a
compensated increment that vanishes exactly at independence, proves
that mean degree, transitivity, and assortativity identify the
generating parameters globally on a certified region, verifies the
regularity conditions for asymptotic inference including the
independence boundary, and applies the resulting methods to eleven
metropolitan networks on two platforms.  The application is a
diagnosis, not a fit: the class is rejected in every city through the
one moment it exists to explain, and the certified monotonicity turns
the rejection into a one-sided conclusion, namely that no admissible
dependence strength reconciles the model with the data and that the
observed assortativity is not evidence of sorting.  What the data are
consistent with, degree heterogeneity and, on one platform, an
unexplained mesoscale residual, is established by null comparisons,
and the mechanisms actually generating the observed transitivity
remain outside the fitted class.

Two elements of the approach apply beyond this model.  The
exact-rational interval method certifies injectivity and rank of any
moment map with polynomial or analytic parameter dependence on a box;
mixture moment problems, exponential-family network statistics, and
spatial point-process summaries have this structure, and the verifier
is packaged for reuse.  The two-view validation design, in which the
same latent quantity is estimated once from auxiliary observables and
once from the graph alone, with a joint law and a specification test
connecting the two, applies whenever a candidate latent attribute is
partially observable, for example activity, wealth, or seniority; the
Austin case shows that the comparison is informative about detection
thresholds as well as about the consistency of the two views, and the
present application shows that it retains value under rejection,
where the direct view supplies the evidence the graph view cannot.

The results have the following limitations.  The certified statements
concern the stated submodels (dimension one with hard range, and the
planar hard-disc case, each with uniform popularity marginal, one
harmonic, and fixed normalizations) on the region $B$; other kernels,
marginals, or normalizations require new certificates.  For the
inhomogeneous fixed-position designs of the application, the
identification, sign-invariance, and boundary properties proved for
the homogeneous model are not established: Section~S5 of the
Supplement states the conditions under which the conditional
estimator is consistent and asymptotically normal, those conditions
are checked per dataset only through the finite-difference Jacobian
diagnostic, which is local, and terms linear in the signed amplitude
can survive under irregular positions, so the graph law need not
factor through $\psi=b^2$ there (the Austin sign-flip comparison in
Table~\ref{tab:gof} measures these terms empirically, and
Proposition~S5.4 turns their presence into a first-order
identification result for the sign).  The realized
limit theory holds for the general bounded finite-range model of
Section~\ref{sec:model}, but conditional edge independence excludes
contagion and strategic formation, and heavy-tailed degree
distributions do not arise under a fixed compact range.  The
decomposition concerns Pearson assortativity; rank assortativity
disagrees with it in sign in several cities, and no rank analogue of
the decomposition is offered.  The empirical analysis conditions on
observed home positions and a chosen spatial score, and its
conclusions concern the eleven metropolitan subgraphs under that
design, not location-based platforms in general.  The collaboration
analysis conditions likewise on institution coordinates, which are
point masses resolved only by a fixed jitter, on a corpus-internal
work count as the productivity proxy, and on an induced node
subsample for the conditional fits, and its stated precision rests
on the atom-clustered standard errors.  Because the fitted
class is rejected, all model-implied quantities reported for the data
($r_{\mathrm{geo}}$, channels, $\Delta_{\mathrm{sort}}$) are
counterfactual diagnostics of the class, not descriptions of the
cities.  Graph-only precision is limited at the available sample
sizes, as the NYC diagnostic shows; city bootstrap bands use
\bootDraws\ draws.

Several extensions are natural, and this paper takes the first steps
along the ordered path the rejection sets.  The augmented class of
Section~\ref{sec:application} carries degree effects, two connection
scales, and closure jointly under the same
goodness-of-fit-before-attribution discipline; its residual pattern,
and the complementary failures of the competitor families, specify
the next enrichment as distance-limited closure or an explicit
community layer, the collaboration refusal sharpens that to a team
layer, papers as hyperedges with sizes drawn from the observed
distribution, and certified identification for the augmented
class, whose moment map is no longer polynomial, is open.  On the
theory side, the boundary laws are now uniform over the certified
null face (Theorem~S4.9) and the signed amplitude is first-order
identified under irregular designs (Proposition~S5.4); an estimator
that exploits the signed expansion, and certified identification for
inhomogeneous position designs, are the remaining theoretical gaps.
A rank-based analogue of the decomposition would remove the Pearson
restriction.  The multi-harmonic theory of Section~S5 of the
Supplement suggests estimating the harmonic energy profile itself
where sample sizes permit.  Section~\ref{sec:collab} supplies the
positive case for the direct view on a collaboration network built
by the replication package's fetch-and-format script; what remains
open there is the graph side, a class that carries the team
mechanism well enough to pass the gate, at which point the two-view
comparison of $\widehat\psi_{\mathrm{dir}}$ with
$\widehat\psi_{\mathrm{gr}}$ becomes available away from the
boundary on real data.  Across both domains the method's verdicts
so far are a certified negative with a rejection, and an interior
detection with a refusal.

\appendix
\section*{Appendix: limit objects}
For reference, the limit intensities used in
\eqref{eq:r-decomp-intro}:
\begin{align}
 \ell&=\rho\int q(w_1,w_2;u)H_x(dw_1)H_x(dw_2)\,dx\,du,\notag\\
 \Lambda(w,x)&=\rho\int q(w,\omega;v)H_x(d\omega)\,dv,
 \label{eq:ell-main-app}\\
 \Gamma(w_1,w_2,x;u)&=\rho\int q(w_1,\omega;v)\,
 q(w_2,\omega;v-u)\,H_x(d\omega)\,dv.\notag
\end{align}

\clearpage
\setcounter{section}{0}
\renewcommand{\thesection}{S\arabic{section}}
\setcounter{table}{0}
\renewcommand{\thetable}{S\arabic{table}}
\setcounter{figure}{0}
\renewcommand{\thefigure}{S\arabic{figure}}
\setcounter{equation}{0}
\renewcommand{\theequation}{S\arabic{equation}}
\setcounter{footnote}{0}
\makeatletter
\renewcommand{\theHsection}{S\arabic{section}}
\renewcommand{\theHtable}{S\arabic{table}}
\renewcommand{\theHfigure}{S\arabic{figure}}
\renewcommand{\theHequation}{S\arabic{equation}}
\makeatother
\begin{center}
{\LARGE Supplementary Material\par}
\end{center}
\bigskip

\noindent This supplement contains the complete proofs.
Section~\ref{sec:si-s1} proves the finite-sample
observational-equivalence result.  Section~\ref{sec:si-s2} proves
the joint edge-Palm, local-count, and realized-assortativity limit
theorem.  Section~\ref{sec:si-s3} contains the certified
identification results, including the global-injectivity
certificate.  Section~\ref{sec:si-s4} proves the interior central
limit theorem, the feasible optimal GMM theory, and the boundary
theory, with all model-side regularity conditions verified for the
certified submodel.  Section~\ref{sec:si-s5} contains the
multi-harmonic identification result, the conditional-design
theory, the graph-versus-direct specification test, and the
certified two-dimensional hard-disc submodel.

\section{Finite-sample observational equivalence and scale aliases}
\label{sec:si-s1}

\paragraph{One-harmonic submodel and observable law.}
Fix an integer $d\geq1$. Let
\[
\mathbb T^d=\mathbb R^d/\mathbb Z^d
\]
be equipped with normalized Haar probability measure $\mu$. Put
\[
F=[-1/2,1/2)^d,
\]
let $r\colon\mathbb T^d\to F$ be the coordinatewise centered
representative, and define
\[
\delta_{\mathbb T^d}(x,y)=r(y-x).
\]
Whenever a torus element appears inside an inner product with an
integer vector, the inner product is evaluated using any Euclidean
lift; the resulting value modulo $\mathbb Z$ is independent of that
lift.

Let
\[
m\in\mathbb Z^d\setminus\{0\},
\qquad
\varphi\in\mathbb R/(2\pi\mathbb Z),
\]
and define
\[
a(v)=\sqrt3(2v-1),
\qquad
\phi_{m,\varphi}(x)
=
\sqrt2\cos(2\pi m^{\trans} x-\varphi),
\]
and
\[
h_{b,\varphi,m}(v\mid x)
=
1+b\,a(v)\phi_{m,\varphi}(x),
\qquad
|b|\leq\frac{1}{\sqrt6}.
\]

Let
\[
\{(X_i,V_i):1\leq i\leq n\}
\]
be i.i.d. with joint density
\[
(x,v)\longmapsto h_{b,\varphi,m}(v\mid x)
\]
relative to $\mu(dx)\,dv$. Let
\[
W_i=\omega(V_i),
\]
where $\omega\colon[0,1]\to[0,\infty)$ is a finite-valued Borel
function.

Fix $\lambda>0$ and $\varepsilon_n>0$, and let
\[
k_{\eta,\nu}\colon\mathbb R^d\to[0,\infty)
\]
be a finite-valued, nonnegative, even Borel function. Conditional on
the latent variables, let
\[
\{A_{ij}:1\leq i<j\leq n\}
\]
be independent Bernoulli variables with
\begin{equation}
\label{eq:exact-exponential-link}
p_{ij}^{(n)}
=
1-\exp\left\{
-\frac{\lambda}{n\varepsilon_n^d}
W_iW_j
k_{\eta,\nu}\left(
\frac{\delta_{\mathbb T^d}(X_i,X_j)}{\varepsilon_n}
\right)
\right\}.
\end{equation}
Set
\[
A_{ji}=A_{ij},
\qquad
A_{ii}=0.
\]

Write
\[
\mathbb P_{\lambda,\varepsilon_n,k,\omega;
b,\varphi,m}^{(n)}
\]
for the unconditional law of the labeled adjacency matrix $A$.
Throughout the proposition and proof,
$\varepsilon_n$, $\omega$, and $k=k_{\eta,\nu}$ are held fixed.
Accordingly, abbreviate
\[
\mathbb P_{\lambda,\eta,\nu;b,\varphi,m}^{(n)}
:=
\mathbb P_{\lambda,\varepsilon_n,k_{\eta,\nu},\omega;
b,\varphi,m}^{(n)}.
\]

Define
\[
\mathcal H_k
=
\left\{
Q\in O(d,\mathbb Z):
k_{\eta,\nu}(Qu)=k_{\eta,\nu}(u)
\text{ for Lebesgue-a.e. }u\in\mathbb R^d
\right\},
\]
where
\[
O(d,\mathbb Z)
=
\left\{
Q\in\mathbb Z^{d\times d}:Q^{\trans} Q=I_d
\right\}.
\]

\begin{siproposition}[
Finite-sample adjacency-law equivalence and orbit-invariant
latent $X$--$V$ dependence]
\label{prop:finite-sample-equivalence}
For every admissible $(b,\varphi,m)$,
$h_{b,\varphi,m}(\cdot\mid x)$ is a probability density for every
$x\in\mathbb T^d$, and
\[
X_i\sim\mu,
\qquad
V_i\sim\operatorname{Unif}(0,1).
\]

Let $Q\in\mathcal H_k$ and $\tau\in\mathbb T^d$. Choose any lift
$\widetilde\tau\in\mathbb R^d$ satisfying
$\widetilde\tau+\mathbb Z^d=\tau$, and define
\[
g_{Q,\tau}(x)=Qx+\tau,
\qquad
m_g=Qm,
\]
and
\[
\varphi_g
=
\varphi+2\pi(Qm)^{\trans}\widetilde\tau
\pmod{2\pi}.
\]
Then $\varphi_g$ is independent of the chosen lift, and for every
$n\geq2$,
\begin{align}
\label{eq:graph-law-equivalence-full}
&
\mathbb P_{\lambda,\varepsilon_n,k_{\eta,\nu},\omega;
b,\varphi,m}^{(n)}
\nonumber\\
&\qquad =
\mathbb P_{\lambda,\varepsilon_n,k_{\eta,\nu},\omega;
b,\varphi_g,m_g}^{(n)}
=
\mathbb P_{\lambda,\varepsilon_n,k_{\eta,\nu},\omega;
-b,\varphi_g+\pi,m_g}^{(n)}.
\end{align}
Equivalently, in the abbreviated notation,
\begin{equation}
\label{eq:graph-law-equivalence}
\mathbb P_{\lambda,\eta,\nu;b,\varphi,m}^{(n)}
=
\mathbb P_{\lambda,\eta,\nu;b,\varphi_g,m_g}^{(n)}
=
\mathbb P_{\lambda,\eta,\nu;-b,\varphi_g+\pi,m_g}^{(n)}.
\end{equation}

For every
$\varphi'\in\mathbb R/(2\pi\mathbb Z)$,
\[
\mathbb P_{\lambda,\eta,\nu;b,\varphi,m}^{(n)}
=
\mathbb P_{\lambda,\eta,\nu;b,\varphi',m}^{(n)}
=
\mathbb P_{\lambda,\eta,\nu;-b,\varphi',m}^{(n)}.
\]
Thus phase is not identified, and the sign of a nonzero amplitude
$b$ is not identified, from the marginal adjacency law.

The set $\mathcal H_k$ is a subgroup of $O(d,\mathbb Z)$. Define
\[
\mathcal O_k(m)=\{Qm:Q\in\mathcal H_k\}.
\]
The marginal adjacency law is constant over all frequency
representatives in $\mathcal O_k(m)$, with the corresponding phase
transformation. Since $k_{\eta,\nu}$ is even,
\[
-I_d\in\mathcal H_k,
\]
so $\mathcal O_k(m)$ is sign-symmetric. No assertion is made that
$\mathcal O_k(m)$ is the complete frequency-equivalence class.

If $b=0$, then
\[
h_{0,\varphi,m}\equiv1,
\]
so every admissible $m\in\mathbb Z^d\setminus\{0\}$ and every
$\varphi$ are irrelevant to the latent and adjacency laws.

Finally, define
\[
\psi(b,\varphi,m)
=
\int_0^1\int_{\mathbb T^d}
\{h_{b,\varphi,m}(v\mid x)-1\}^2
\,\mu(dx)\,dv.
\]
Then
\[
\psi(b,\varphi,m)=b^2.
\]
Thus $\psi$ is invariant under the transformations in
\eqref{eq:graph-law-equivalence}. It measures latent $X$--$V$
dependence and need not be identified from the adjacency law under
the assumptions of this proposition.
\end{siproposition}

\begin{proof}
For every $v\in[0,1]$ and $x\in\mathbb T^d$,
\[
|a(v)|\leq\sqrt3,
\qquad
|\phi_{m,\varphi}(x)|\leq\sqrt2.
\]
Hence
\[
|b\,a(v)\phi_{m,\varphi}(x)|
\leq |b|\sqrt6
\leq1,
\]
so
\[
h_{b,\varphi,m}(v\mid x)\geq0.
\]

Direct calculation gives
\[
\int_0^1a(v)\,dv=0,
\qquad
\int_0^1a(v)^2\,dv=1.
\]
Consequently, for every $x$,
\[
\int_0^1h_{b,\varphi,m}(v\mid x)\,dv
=
1+b\phi_{m,\varphi}(x)\int_0^1a(v)\,dv
=1.
\]
Thus $h_{b,\varphi,m}(\cdot\mid x)$ is a probability density, and
the $X$-marginal is $\mu$.

Because $m\neq0$, the torus character
$x\mapsto e^{2\pi\mathrm{i}m^{\trans} x}$ is nonconstant. Therefore
\[
\int_{\mathbb T^d}
e^{2\pi\mathrm{i}m^{\trans} x}\,\mu(dx)=0.
\]
Since
\[
\phi_{m,\varphi}(x)
=
\sqrt2\operatorname{Re}\left\{
e^{-\mathrm{i}\varphi}
e^{2\pi\mathrm{i}m^{\trans} x}
\right\},
\]
we obtain
\[
\int_{\mathbb T^d}
\phi_{m,\varphi}(x)\,\mu(dx)=0.
\]
Furthermore,
\[
\phi_{m,\varphi}(x)^2
=
1+\cos(4\pi m^{\trans} x-2\varphi).
\]
Because $2m\neq0$, the same character-orthogonality argument yields
\[
\int_{\mathbb T^d}
\phi_{m,\varphi}(x)^2\,\mu(dx)=1.
\]
It follows that, for every $v$,
\[
\int_{\mathbb T^d}
h_{b,\varphi,m}(v\mid x)\,\mu(dx)=1,
\]
and hence
\[
V_i\sim\operatorname{Unif}(0,1).
\]

We next verify that $\mathcal H_k$ is a subgroup. The identity belongs
to $\mathcal H_k$. For $Q,R\in\mathcal H_k$, define
\[
N_Q=\{u:k_{\eta,\nu}(Qu)\neq k_{\eta,\nu}(u)\},
\qquad
N_R=\{u:k_{\eta,\nu}(Ru)\neq k_{\eta,\nu}(u)\}.
\]
These sets are Lebesgue-null. Outside
\[
N_R\cup R^{-1}N_Q
\]
we have
\[
k_{\eta,\nu}(QRu)
=
k_{\eta,\nu}(Ru)
=
k_{\eta,\nu}(u).
\]
Because orthogonal transformations preserve Lebesgue-null sets,
$QR\in\mathcal H_k$. Similarly, outside $QN_Q$, putting
$v=Q^{-1}u$ gives
\[
k_{\eta,\nu}(u)
=
k_{\eta,\nu}(Qv)
=
k_{\eta,\nu}(v)
=
k_{\eta,\nu}(Q^{-1}u).
\]
Thus $Q^{-1}\in\mathcal H_k$, proving the subgroup claim.

We now prove the graph-law equivalence. If
$\widetilde\tau'$ is another lift of $\tau$, then
\[
\widetilde\tau'-\widetilde\tau\in\mathbb Z^d.
\]
Since $Qm\in\mathbb Z^d$,
\[
2\pi(Qm)^{\trans}
(\widetilde\tau'-\widetilde\tau)
\in2\pi\mathbb Z.
\]
Thus $\varphi_g$ is well defined modulo $2\pi$.

Choose a lift $\widetilde x\in\mathbb R^d$ of
$x\in\mathbb T^d$. Then
\[
Q^{\trans}(\widetilde x-\widetilde\tau)
\]
is a lift of $g_{Q,\tau}^{-1}(x)$. Therefore
\begin{align*}
\phi_{m,\varphi}
\bigl(g_{Q,\tau}^{-1}(x)\bigr)
&=
\sqrt2\cos\left(
2\pi m^{\trans} Q^{\trans}
(\widetilde x-\widetilde\tau)-\varphi
\right)\\
&=
\sqrt2\cos\left(
2\pi(Qm)^{\trans}\widetilde x
-\{\varphi+2\pi(Qm)^{\trans}\widetilde\tau\}
\right)\\
&=
\phi_{m_g,\varphi_g}(x).
\end{align*}
Hence
\begin{equation}
\label{eq:latent-density-pushforward}
h_{b,\varphi,m}
\bigl(v\mid g_{Q,\tau}^{-1}(x)\bigr)
=
h_{b,\varphi_g,m_g}(v\mid x).
\end{equation}

Under the model indexed by $(b,\varphi,m)$, define
\[
X_i^\star=g_{Q,\tau}(X_i),
\qquad
V_i^\star=V_i,
\qquad
W_i^\star=\omega(V_i^\star)=W_i.
\]
Affine torus automorphisms preserve normalized Haar measure.
Therefore \eqref{eq:latent-density-pushforward} shows that
\[
\{(X_i^\star,V_i^\star):1\leq i\leq n\}
\]
are i.i.d. with the latent law indexed by
$(b,\varphi_g,m_g)$.

Every matrix in $O(d,\mathbb Z)$ is a signed permutation matrix:
each integer row has Euclidean norm one and therefore contains
exactly one nonzero entry, equal to $1$ or $-1$, while orthogonality
places these entries in distinct columns.

Define
\[
B_Q
=
\{z\in\mathbb T^d:r(Qz)\neq Qr(z)\}.
\]
Since $Q$ is a signed permutation,
\[
B_Q
\subseteq
\left\{
z:\text{some coordinate of }r(z)\text{ equals }-1/2
\right\},
\]
and therefore
\[
\mu(B_Q)=0.
\]
Consequently,
\[
\delta_{\mathbb T^d}
\bigl(g_{Q,\tau}(x),g_{Q,\tau}(y)\bigr)
=
Q\delta_{\mathbb T^d}(x,y)
\]
for $(\mu\otimes\mu)$-almost every $(x,y)$.

For the fixed $Q\in\mathcal H_k$, put
\[
N_Q
=
\{u\in\mathbb R^d:
k_{\eta,\nu}(Qu)\neq k_{\eta,\nu}(u)\}.
\]
Then $\operatorname{Leb}(N_Q)=0$. Since $X_i$ and $X_j$ are
independent and Haar-uniform, $X_j-X_i$ is Haar-uniform on
$\mathbb T^d$. Thus
\[
D_{ij}
=
\delta_{\mathbb T^d}(X_i,X_j)
\]
is uniform on $F$. With
\[
F/\varepsilon_n
=
\{u\in\mathbb R^d:\varepsilon_nu\in F\},
\]
the variable
\[
U_{ij}^{(n)}
=
\frac{D_{ij}}{\varepsilon_n}
\]
has Lebesgue density
\[
u\longmapsto
\varepsilon_n^d
\mathbf1_{F/\varepsilon_n}(u).
\]
It follows that
\[
\Pr\{U_{ij}^{(n)}\in N_Q\}=0.
\]

Combining the kernel exceptional event with the event
$\{X_j-X_i\in B_Q\}$ and taking a finite union over the
$\binom n2$ pairs shows that, with probability one,
\begin{align*}
k_{\eta,\nu}\left(
\frac{\delta_{\mathbb T^d}(X_i^\star,X_j^\star)}
     {\varepsilon_n}
\right)
&=
k_{\eta,\nu}\left(
Q\frac{\delta_{\mathbb T^d}(X_i,X_j)}
        {\varepsilon_n}
\right)\\
&=
k_{\eta,\nu}\left(
\frac{\delta_{\mathbb T^d}(X_i,X_j)}
     {\varepsilon_n}
\right)
\end{align*}
simultaneously for all $i<j$.

Let $p_{ij}^{\star(n)}$ denote the probability in
\eqref{eq:exact-exponential-link} computed from
$(X_i^\star,V_i^\star)$ under the transformed latent law. Since
$W_i^\star=W_i$, the preceding identity implies
\[
p_{ij}^{\star(n)}=p_{ij}^{(n)}
\]
simultaneously for all $i<j$, almost surely.

Let
\[
\{Z_{ij}:1\leq i<j\leq n\}
\]
be independent $\operatorname{Unif}(0,1)$ variables, independent of
the latent variables, and define
\[
A_{ij}
=
\mathbf1\{Z_{ij}\leq p_{ij}^{(n)}\},
\qquad
A_{ij}^\star
=
\mathbf1\{Z_{ij}\leq p_{ij}^{\star(n)}\}.
\]
Then
\[
A^\star=A
\qquad\text{almost surely}.
\]
The transformed latent array has the law indexed by
$(b,\varphi_g,m_g)$; hence this coupling proves the first equality in
\eqref{eq:graph-law-equivalence-full}.

Next,
\[
\phi_{m,\varphi+\pi}(x)
=
-\phi_{m,\varphi}(x),
\]
and therefore
\begin{align*}
h_{-b,\varphi+\pi,m}(v\mid x)
&=
1-b\,a(v)\phi_{m,\varphi+\pi}(x)\\
&=
1+b\,a(v)\phi_{m,\varphi}(x)\\
&=
h_{b,\varphi,m}(v\mid x).
\end{align*}
Applying this identity to
$(\varphi,m)=(\varphi_g,m_g)$ proves the second equality in
\eqref{eq:graph-law-equivalence-full}.

To prove arbitrary phase invariance, choose an index $r$ such that
$m_r\neq0$. Given
$\varphi'\in\mathbb R/(2\pi\mathbb Z)$, choose a real lift
$\widetilde\Delta$ of $\varphi'-\varphi$ and define a lift of
$\tau\in\mathbb T^d$ by
\[
\widetilde\tau_r
=
\frac{\widetilde\Delta}{2\pi m_r},
\qquad
\widetilde\tau_\ell=0
\quad(\ell\neq r).
\]
Then
\[
2\pi m^{\trans}\widetilde\tau
=
\varphi'-\varphi
\pmod{2\pi}.
\]
Because $I_d\in\mathcal H_k$,
\[
\mathbb P_{\lambda,\eta,\nu;b,\varphi,m}^{(n)}
=
\mathbb P_{\lambda,\eta,\nu;b,\varphi',m}^{(n)}.
\]
Furthermore,
\[
h_{b,\varphi'+\pi,m}
=
h_{-b,\varphi',m}.
\]
Using phase invariance once more gives
\[
\mathbb P_{\lambda,\eta,\nu;b,\varphi',m}^{(n)}
=
\mathbb P_{\lambda,\eta,\nu;b,\varphi'+\pi,m}^{(n)}
=
\mathbb P_{\lambda,\eta,\nu;-b,\varphi',m}^{(n)}.
\]

The exact identity
\[
\phi_{-m,-\varphi}(x)=\phi_{m,\varphi}(x)
\]
holds independently of the kernel. Since $k_{\eta,\nu}$ is even,
\[
-I_d\in\mathcal H_k.
\]
Thus $\mathcal O_k(m)$ is sign-symmetric. Combining the
$Q$-equivalences with phase invariance proves constancy of the
adjacency law over $\mathcal O_k(m)$. This argument supplies no
converse. If $b=0$, then
\[
h_{0,\varphi,m}\equiv1,
\]
so every admissible frequency and every phase disappear from the
latent and adjacency laws.

Finally, Tonelli's theorem gives
\begin{align*}
\psi(b,\varphi,m)
&=
b^2
\left\{\int_0^1a(v)^2\,dv\right\}
\left\{
\int_{\mathbb T^d}
\phi_{m,\varphi}(x)^2\,\mu(dx)
\right\}\\
&=b^2.
\end{align*}
With the convention
\[
\chi^2(P\|Q)
=
\int\left(\frac{dP}{dQ}-1\right)^2\,dQ,
\]
the uniform marginal calculations imply
\[
\psi(b,\varphi,m)
=
\chi^2\left(
\mathcal L(X_i,V_i)
\,\middle\|\,
\mu\otimes\operatorname{Unif}(0,1)
\right).
\]
This need not equal
\[
\chi^2\left(
\mathcal L(X_i,W_i)
\,\middle\|\,
\mu\otimes\mathcal L(W_i)
\right)
\]
when $\omega$ is non-injective, and it does not imply graph-law
identification.
\end{proof}

\paragraph{Scale aliases in an enlarged unnormalized model.}
This paragraph concerns the fully indexed family
\[
\mathbb P_{\lambda,\varepsilon_n,k,\omega;
b,\varphi,m}^{(n)}.
\]
Suppose that the admissible class is closed under the following
transformations for constants $\alpha,\beta,s>0$:
\[
\omega^\star(v)=\alpha\omega(v),
\qquad
\varepsilon_n^\star=s\varepsilon_n,
\]
\[
k^\star(u)=\beta k(su),
\qquad
\lambda^\star
=
\frac{s^d}{\alpha^2\beta}\lambda.
\]
Under the same paired latent variables $(X_i,V_i)$,
\[
W_i^\star=\omega^\star(V_i)=\alpha W_i.
\]
For every latent realization,
\begin{align*}
&
\frac{\lambda^\star}
     {n(\varepsilon_n^\star)^d}
W_i^\star W_j^\star
k^\star\left(
\frac{\delta_{\mathbb T^d}(X_i,X_j)}
     {\varepsilon_n^\star}
\right)\\
&\quad=
\frac{(s^d\lambda)/(\alpha^2\beta)}
     {n(s\varepsilon_n)^d}
(\alpha W_i)(\alpha W_j)
\beta k\left(
s\frac{\delta_{\mathbb T^d}(X_i,X_j)}
        {s\varepsilon_n}
\right)\\
&\quad=
\frac{\lambda}{n\varepsilon_n^d}
W_iW_j
k\left(
\frac{\delta_{\mathbb T^d}(X_i,X_j)}
     {\varepsilon_n}
\right).
\end{align*}
Therefore
\[
\mathbb P_{\lambda,\varepsilon_n,k,\omega;
b,\varphi,m}^{(n)}
=
\mathbb P_{\lambda^\star,\varepsilon_n^\star,k^\star,\omega^\star;
b,\varphi,m}^{(n)}.
\]

If the model is parameterized through
\[
H_x=\mathcal L(W\mid X=x),
\]
then the complete transformed conditional law is
\[
H_x^\star(B)=H_x(B/\alpha),
\qquad
B/\alpha=\{w\geq0:\alpha w\in B\}.
\]

The three special cases are
\[
(\omega,\lambda)
\longmapsto
(c\omega,\lambda/c^2),
\]
\[
(k,\lambda)
\longmapsto
(ck,\lambda/c),
\]
and
\[
(\varepsilon_n,k,\lambda)
\longmapsto
\bigl(s\varepsilon_n,k^{(s)},s^d\lambda\bigr),
\qquad
k^{(s)}(u)=k(su).
\]

Each equality is a within-model alias only when its corresponding
transformed objects remain admissible. The popularity-scale alias is
excluded for the joint law of $(A,W)$ when $W$ is observed, and is
excluded from the latent model when the scale of $\omega$ or $H_x$
is fixed. The kernel-amplitude alias is excluded when kernel
amplitude is normalized or the kernel family is not closed under
multiplication. The bandwidth alias is excluded when
$\varepsilon_n$ is fixed, when a prescribed value
\[
\rho=\lim_{n\to\infty}n\varepsilon_n^d\in(0,\infty)
\]
is held fixed, or when an independent kernel-range restriction not
preserved by dilation is imposed. Restrictions on one component do
not necessarily remove the other aliases.

If, additionally,
\[
k\in L^1(\mathbb R^d),
\qquad
\int_{\mathbb R^d}k(u)\,du=1,
\]
then integral normalization alone does not remove the bandwidth
alias. Indeed, setting
\[
k^\star(u)=s^d k(su),
\qquad
\varepsilon_n^\star=s\varepsilon_n,
\qquad
\lambda^\star=\lambda
\]
gives
\[
\int_{\mathbb R^d}k^\star(u)\,du
=
\int_{\mathbb R^d}k(z)\,dz
=1
\]
and preserves the edge exponent. A restriction not invariant under
kernel dilation, such as fixed bandwidth or fixed kernel support
radius, is needed to exclude this particular orbit. Excluding these
scaling directions does not by itself establish complete
identification.

\paragraph{Scope and nonidentification boundary.}
Proposition~\ref{prop:finite-sample-equivalence} concerns only the
unconditional labeled adjacency law after integration over latent
locations and marks. It does not assert equivalence conditional on
fixed or observed locations or completeness of the displayed
equivalence class.

Positive graph-law identification of $\psi$ can fail under the broad
assumptions above. For example, if $\omega$ is constant, then the
adjacency law depends only on the marginal Haar locations and not on
$h_{b,\varphi,m}$. The graph law is consequently identical for
different values of $b^2$. Any positive-identification theorem for
$\psi$ therefore requires additional assumptions and a separate
injectivity proof.

\clearpage
\section{Joint edge-Palm, local-count, and realized-assortativity
limits}
\label{sec:si-s2}

\begin{sitheorem}[Joint edge-Palm, local-count, and realized-assortativity limits]
\label{thm:joint-realized-local-limits}
Fix \(d\geq1\), and let
\[
\mathbb T^d=\mathbb R^d/\mathbb Z^d
\]
carry normalized Haar measure \(dx\). For \(x,y\in\mathbb T^d\), let
\[
\delta_{\mathbb T^d}(x,y)\in[-1/2,1/2)^d
\]
denote the centered representative of \(y-x\), with any measurable
convention on the Haar-null half-period boundary.
Throughout, \(d_{\mathbb T^d}\) is the Euclidean quotient metric.

Let \(\varepsilon_n>0\) and
\[
\rho_n=n\varepsilon_n^d\longrightarrow\rho\in(0,\infty).
\]
Let \((X_i,W_i)_{i\geq1}\) be i.i.d. with law
\[
dx\,H_x(dw)
\]
on \(\mathbb T^d\times[0,\infty)\), where \(W_i<\infty\) almost
surely and \(x\mapsto H_x\) is a probability kernel satisfying
\begin{equation}
\label{eq:uniform-TV-H-final}
\omega_H(a)
:=
\sup_{\substack{x,y\in\mathbb T^d\\
d_{\mathbb T^d}(x,y)\leq a}}
\|H_x-H_y\|_{\mathrm{TV}}
\longrightarrow0
\qquad(a\downarrow0).
\end{equation}

Let \(k=k_{\eta,\nu}:\mathbb R^d\to[0,\infty)\) be Borel measurable,
bounded, and even. Assume that, pointwise,
\[
k(u)=0\qquad\text{for }u\notin K,
\]
where \(K\subseteq B(0,R)\) is compact and \(R<\infty\). Fix
\(\lambda>0\), and define
\begin{align}
U_{ij,n}
&=
\frac{\delta_{\mathbb T^d}(X_i,X_j)}{\varepsilon_n},
\label{eq:scaled-displacement-final}\\
q_{\vartheta,n}(w,w';u)
&=
1-\exp\left\{
-\frac{\lambda}{\rho_n}ww'k_{\eta,\nu}(u)
\right\}.
\label{eq:finite-link-final}
\end{align}
Conditionally on \(\{(X_i,W_i)\}_{i=1}^n\), suppose that
\(\{A_{ij}:i<j\}\) are independent and
\[
A_{ij}\sim
\operatorname{Bernoulli}
\left\{
q_{\vartheta,n}(W_i,W_j;U_{ij,n})
\right\}.
\]
Set \(A_{ji}=A_{ij}\) and \(A_{ii}=0\).

Define the limiting link
\[
q_\vartheta(w,w';u)
=
1-\exp\left\{
-\frac{\lambda}{\rho}ww'k_{\eta,\nu}(u)
\right\}.
\]
All \(x\)-integrals below are with respect to normalized Haar
measure, and all \(u,v\)-integrals are with respect to Lebesgue
measure. Put
\begin{align}
\ell(\vartheta)
&=
\rho
\int
q_\vartheta(w_1,w_2;u)
H_x(dw_1)H_x(dw_2)\,dx\,du,
\label{eq:ell-final}\\
\Lambda_\vartheta(w,x)
&=
\rho
\int
q_\vartheta(w,\omega;v)H_x(d\omega)\,dv,
\label{eq:Lambda-final}\\
\Gamma_\vartheta(w_1,w_2,x;u)
&=
\rho
\int
q_\vartheta(w_1,\omega;v)
q_\vartheta(w_2,\omega;v-u)
H_x(d\omega)\,dv.
\label{eq:Gamma-final}
\end{align}
Assume
\[
\ell(\vartheta)>0,
\]
and define
\begin{equation}
\label{eq:Pi-final}
\Pi_\vartheta(dw_1,dw_2,dx,du)
=
\frac{\rho}{\ell(\vartheta)}
q_\vartheta(w_1,w_2;u)
H_x(dw_1)H_x(dw_2)\,dx\,du.
\end{equation}

Writing
\[
q_{ab}(r)=q_\vartheta(w_a,w_b;r),
\]
define
\begin{align}
t(\vartheta)
&=
\frac{\rho^2}{6}
\int
q_{12}(u)q_{13}(v)q_{23}(v-u)
\prod_{j=1}^3H_x(dw_j)\,dx\,du\,dv,
\label{eq:t-final}\\
s_2(\vartheta)
&=
\frac{\rho^2}{2}
\int
q_{12}(u)q_{13}(v)
\prod_{j=1}^3H_x(dw_j)\,dx\,du\,dv.
\label{eq:s2-final}
\end{align}
The identities proved below imply
\[
s_2(\vartheta)\geq\frac{\ell(\vartheta)^2}{2}>0.
\]
Consequently,
\begin{equation}
\label{eq:C-final}
C(\vartheta)=\frac{3t(\vartheta)}{s_2(\vartheta)}
\end{equation}
is well defined.

For
\[
z=(w_1,w_2,x,u),\qquad
L_j(z)=\Lambda_\vartheta(w_j,x),\qquad
G(z)=\Gamma_\vartheta(w_1,w_2,x;u),
\]
let, conditionally on \(z\),
\[
P_{11}\sim\operatorname{Pois}\{G(z)\},\qquad
P_{10}\sim\operatorname{Pois}\{L_1(z)-G(z)\},
\]
\[
P_{01}\sim\operatorname{Pois}\{L_2(z)-G(z)\}
\]
be mutually independent, and set
\begin{equation}
\label{eq:limiting-degrees-final}
D_1^\star=1+P_{10}+P_{11},
\qquad
D_2^\star=1+P_{01}+P_{11}.
\end{equation}
Let \(\mathsf Q_\vartheta\) be the joint law of
\[
(W_1,W_2,X,U,D_1^\star,D_2^\star)
\]
obtained by first drawing
\((W_1,W_2,X,U)\sim\Pi_\vartheta\) and then applying
\eqref{eq:limiting-degrees-final}.

Let
\[
\mathsf S
=
[0,\infty)^2\times\mathbb T^d\times\mathbb R^d
\times\mathbb N_0^2
\]
have the product topology, where \(\mathbb N_0\) has the discrete
topology. Define
\begin{equation}
\label{eq:C2-final}
\mathcal C_2
=
\left\{
f\in C(\mathsf S):
\|f\|_{(2)}
:=
\sup_{\mathsf S}
\frac{|f(w_1,w_2,x,u,d_1,d_2)|}
     {1+d_1^2+d_2^2}
<\infty
\right\}.
\end{equation}
Every assertion below involving \(\mathcal C_2\) is pointwise for
each fixed, \(n\)-independent \(f\in\mathcal C_2\). No uniform
statement over \(\mathcal C_2\) is asserted.  In particular, the
envelope in \eqref{eq:C2-final} is uniform in the unbounded mark and
displacement coordinates; no score growth in \(w_1,w_2\), or \(u\) is
covered.

For the observed graph, define
\[
D_i=\sum_{j\ne i}A_{ij},
\qquad
E_n=\sum_{i<j}A_{ij},
\]
\[
T_n=\sum_{i<j<k}A_{ij}A_{ik}A_{jk},
\qquad
S_{2,n}=\sum_{i=1}^n\binom{D_i}{2},
\]
and the unnormalized directed-edge measure
\begin{equation}
\label{eq:M-final}
\mathsf M_n
=
\frac1n
\sum_{i\ne j}
A_{ij}
\delta_{(W_i,W_j,X_i,U_{ij,n},D_i,D_j)}.
\end{equation}

Then the following assertions hold.

\begin{enumerate}
\item
If \(p_n=P_\vartheta(A_{12}=1)\), then
\[
np_n\longrightarrow\ell(\vartheta).
\]
In particular, \(p_n>0\) for all sufficiently large \(n\), and,
for those \(n\),
\begin{equation}
\label{eq:joint-Palm-TV-final}
\left\|
\mathcal L_\vartheta
\left(
W_1,W_2,X_1,U_{12,n},D_1,D_2
\mid A_{12}=1
\right)
-
\mathsf Q_\vartheta
\right\|_{\mathrm{TV}}
\longrightarrow0,
\end{equation}
where both laws are regarded as probability measures on the common
space \(\mathsf S\).

\item
For every fixed \(f\in\mathcal C_2\),
\begin{equation}
\label{eq:M-f-final}
E_\vartheta
\left|
\mathsf M_n(f)
-
\ell(\vartheta)\mathsf Q_\vartheta(f)
\right|^2
\longrightarrow0.
\end{equation}
Moreover,
\[
\mathsf M_n
\Rightarrow
\ell(\vartheta)\mathsf Q_\vartheta
\]
weakly in probability as finite Borel measures on \(\mathsf S\).
The scalar convergence in \eqref{eq:M-f-final} holds, in particular,
for
\[
1,\quad d_1,\quad d_2,\quad d_1^2,\quad d_2^2,\quad d_1d_2.
\]

\item
Jointly in \(L^2\),
\begin{equation}
\label{eq:count-LLN-final}
\left(
\frac{2E_n}{n},
\frac{T_n}{n},
\frac{S_{2,n}}{n}
\right)
\longrightarrow
\left(
\ell(\vartheta),
t(\vartheta),
s_2(\vartheta)
\right).
\end{equation}
This convergence is joint with \eqref{eq:M-f-final} for every fixed
finite collection of functions from \(\mathcal C_2\).

\item
Define
\[
\widehat C_n
=
\begin{cases}
3T_n/S_{2,n},&S_{2,n}>0,\\
0,&S_{2,n}=0.
\end{cases}
\]
On \(E_n>0\), define
\[
\widehat{\mathsf Q}_n
=
\frac{\mathsf M_n}{2E_n/n},
\qquad
\overline D_{e,n}
=
\int d_1\,d\widehat{\mathsf Q}_n,
\]
\[
\widehat V_{e,n}
=
\int(d_1-\overline D_{e,n})^2
\,d\widehat{\mathsf Q}_n,
\]
and
\[
\widehat\kappa_{e,n}
=
\int
(d_1-\overline D_{e,n})
(d_2-\overline D_{e,n})
\,d\widehat{\mathsf Q}_n.
\]
On \(E_n=0\), set
\[
\widehat{\mathsf Q}_n=0,
\qquad
\overline D_{e,n}=\widehat V_{e,n}=\widehat\kappa_{e,n}=0.
\]
Set
\[
\widehat r_n
=
\begin{cases}
\widehat\kappa_{e,n}/\widehat V_{e,n},
&
E_n>0\ \text{and}\ \widehat V_{e,n}>0,\\
0,&\text{otherwise}.
\end{cases}
\]
Then
\begin{equation}
\label{eq:r-final}
r(\vartheta)
=
\frac{
\operatorname{Cov}_{\Pi_\vartheta}(L_1,L_2)
+
E_{\Pi_\vartheta}G
}{
\operatorname{Var}_{\Pi_\vartheta}(L_1)
+
E_{\Pi_\vartheta}L_1
}
\end{equation}
is well defined, and
\begin{equation}
\label{eq:summary-vector-final}
\left(
\frac{2E_n}{n},
\frac{T_n}{n},
\frac{S_{2,n}}{n},
\widehat C_n,
\widehat r_n
\right)
\overset{\mathbb P}\longrightarrow
\left(
\ell(\vartheta),
t(\vartheta),
s_2(\vartheta),
C(\vartheta),
r(\vartheta)
\right).
\end{equation}
\end{enumerate}
\end{sitheorem}

\begin{proof}
Suppress \(\vartheta\) from the notation.

\medskip
\noindent
\paragraph{Eventual torus geometry and finiteness.}
Since \(n\varepsilon_n^d\to\rho\in(0,\infty)\), we have
\(\varepsilon_n\to0\). Hence there is \(n_0\) such that, for
\(n\geq n_0\),
\[
\frac{\rho}{2}\leq\rho_n\leq2\rho,
\qquad
3R\varepsilon_n<\frac12.
\]
The latter condition ensures that all centered-lift identities used
below are unambiguous on \(K\), \(K-K\), and the three-hop
neighborhoods appearing in the influence argument.
More explicitly, if \(z_0,\ldots,z_m\in\mathbb T^d\), \(m\leq3\),
and each successive centered displacement has Euclidean norm at most
\(R\varepsilon_n\), then the recursively lifted path in \(\mathbb R^d\)
has no coordinate wraparound.  Hence every centered displacement between
two vertices on that path equals the corresponding difference of the lifts.

Because \(0\leq q_n,q\leq1\) and both links vanish outside \(K\),
all limiting integrals are finite. In particular,
\[
0\leq\Lambda(w,x)\leq\rho|K|,
\qquad
0\leq\Gamma(w_1,w_2,x;u)\leq\rho|K|.
\]

\medskip
\noindent
\paragraph{Uniform convergence of the links.}
Put
\[
c_n=\frac{\lambda}{\rho_n},
\qquad
c=\frac{\lambda}{\rho}.
\]
For \(s=ww'k(u)\geq0\),
\[
|e^{-c_ns}-e^{-cs}|
\leq
|c_n-c|s e^{-(c_n\wedge c)s}
\leq
\frac{|c_n-c|}{e(c_n\wedge c)}.
\]
Consequently,
\begin{equation}
\label{eq:q-uniform-final}
\Delta_{q,n}
:=
\sup_{w,w'\geq0,\ u\in\mathbb R^d}
|q_n(w,w';u)-q(w,w';u)|
\longrightarrow0.
\end{equation}
Thus no moment or boundedness assumption on \(W_i\) is required for
this step.

\medskip
\noindent
\paragraph{The edge probability and its normalizer.}
Let
\[
\mathcal F_n
=
\left[-\frac1{2\varepsilon_n},
      \frac1{2\varepsilon_n}\right)^d.
\]
For fixed \(x\), the exact torus change of variables
\[
X_2=x+\varepsilon_nu
\]
gives
\[
p_n=\varepsilon_n^d I_n,
\]
where
\[
I_n
=
\int_{\mathbb T^d}\int_{\mathcal F_n}
q_n(w_1,w_2;u)
H_x(dw_1)H_{x+\varepsilon_nu}(dw_2)
\,du\,dx.
\]
For \(n\geq n_0\), \(K\subset\mathcal F_n\), and the integrand
vanishes off \(K\). Thus it may be extended by zero to
\(\mathbb R^d\).

Let
\[
I
=
\int
q(w_1,w_2;u)
H_x(dw_1)H_x(dw_2)\,dx\,du
=
\frac{\ell}{\rho}.
\]
A product-measure total-variation bound gives
\[
|I_n-I|
\leq
C|K|
\left\{
\Delta_{q,n}
+
\omega_H(R\varepsilon_n)
\right\}
\longrightarrow0.
\]
Therefore
\[
np_n
=
n\varepsilon_n^d I_n
=
\rho_n I_n
\longrightarrow
\rho I
=
\ell.
\]
Because \(\ell>0\), \(p_n>0\) for all sufficiently large \(n\).

\medskip
\noindent
\paragraph{The finite endpoint-Palm measure.}
The exact conditional law of
\[
Z_{12,n}=(W_1,W_2,X_1,U_{12,n})
\]
given \(A_{12}=1\) is
\begin{equation}
\label{eq:Pi-n-final}
\Pi_n(dw_1,dw_2,dx,du)
=
\frac1{I_n}
q_n(w_1,w_2;u)
H_x(dw_1)H_{x+\varepsilon_nu}(dw_2)
\,dx\,du,
\end{equation}
where the measure is extended by zero off \(\mathcal F_n\).

Let \(\nu_n\) and \(\nu\) denote the corresponding unnormalized
finite measures. The same argument as for the edge probability, now tested against
arbitrary measurable functions bounded by one, gives
\[
\|\nu_n-\nu\|_{\mathrm{TV}}
\leq
C|K|
\left\{
\Delta_{q,n}
+
\omega_H(R\varepsilon_n)
\right\}
\longrightarrow0.
\]
Since \(I_n\to I=\ell/\rho>0\),
\begin{equation}
\label{eq:Pi-TV-final}
\|\Pi_n-\Pi\|_{\mathrm{TV}}\longrightarrow0.
\end{equation}

\medskip
\noindent
\paragraph{The exact third-vertex category law.}
Fix \(z=(w_1,w_2,x,u)\) in a set of full \(\Pi_n\)-measure, and put
\[
x_2=x+\varepsilon_nu.
\]
For a third vertex, define the \(11,10,01\) category probabilities by
\begin{align}
a_{11,n}(z)
&=
\varepsilon_n^d
\int
q_n(w_1,\omega;v)
q_n(w_2,\omega;v-u)
H_{x+\varepsilon_nv}(d\omega)\,dv,
\label{eq:a11-final}\\
a_{10,n}(z)
&=
\varepsilon_n^d
\int
q_n(w_1,\omega;v)
\{1-q_n(w_2,\omega;v-u)\}
H_{x+\varepsilon_nv}(d\omega)\,dv,
\label{eq:a10-final}\\
a_{01,n}(z)
&=
\varepsilon_n^d
\int
\{1-q_n(w_1,\omega;v)\}
q_n(w_2,\omega;v-u)
H_{x+\varepsilon_nv}(d\omega)\,dv.
\label{eq:a01-final}
\end{align}
These formulas are exact for \(n\geq n_0\). Indeed, \(u\in K\)
\(\Pi_n\)-almost surely, and the relevant \(v\)-domain is contained
in
\[
K\cup(u+K)\subseteq B(0,2R).
\]
On this domain, the centered lift from \(x_2\) to
\(x+\varepsilon_nv\) is \(\varepsilon_n(v-u)\).

Conditionally on \(z\) and \(A_{12}=1\), the category vectors
associated with vertices \(3,\ldots,n\) are i.i.d. Hence
\[
(N_{10,n},N_{01,n},N_{11,n})
\sim
\operatorname{Mult}_{n-2}
(a_{10,n},a_{01,n},a_{11,n}),
\]
with the remaining probability assigned to \(00\), and
\[
D_1=1+N_{10,n}+N_{11,n},
\qquad
D_2=1+N_{01,n}+N_{11,n}.
\]

\medskip
\noindent
\paragraph{Convergence of the category means.}
Put
\[
\beta_{ab,n}(z)=(n-2)a_{ab,n}(z)
\]
and
\[
b_{11}(z)=G(z),\qquad
b_{10}(z)=L_1(z)-G(z),\qquad
b_{01}(z)=L_2(z)-G(z).
\]
Since \(0\leq q\leq1\),
\[
0\leq G(z)\leq L_1(z).
\]
After the change of variables \(r=v-u\),
\[
0\leq G(z)\leq L_2(z).
\]
Thus \(b_{10},b_{01},b_{11}\) are nonnegative.

Using the fixed finite \(v\)-domain, the product-measure
total-variation inequality, \eqref{eq:q-uniform-final}, and
\((n-2)\varepsilon_n^d\to\rho\), we obtain
\begin{equation}
\label{eq:category-means-final}
\int
\sum_{ab\in\{10,01,11\}}
|\beta_{ab,n}(z)-b_{ab}(z)|
\,\Pi_n(dz)
\leq
C
\left\{
|(n-2)\varepsilon_n^d-\rho|
+
\Delta_{q,n}
+
\omega_H(2R\varepsilon_n)
\right\}
\longrightarrow0.
\end{equation}

\medskip
\noindent
\paragraph{Multivariate Poisson approximation.}
Let
\[
a_{+,n}
=
a_{10,n}+a_{01,n}+a_{11,n}.
\]
By the union bound and the finite range of the link,
\[
a_{+,n}(z)
\leq
2|K|\varepsilon_n^d.
\]
The multivariate Le Cam inequality therefore yields
\begin{align}
&
\left\|
\mathcal L
\left(
N_{10,n},N_{01,n},N_{11,n}
\mid z,A_{12}=1
\right)
-
\bigotimes_{ab}
\operatorname{Pois}\{\beta_{ab,n}(z)\}
\right\|_{\mathrm{TV}}
\nonumber\\
&\qquad
\leq
2(n-2)a_{+,n}(z)^2
\leq
C n\varepsilon_n^{2d}
=
O(n^{-1}).
\label{eq:Le-Cam-final}
\end{align}
Moreover, a Poisson thinning coupling gives
\[
\left\|
\bigotimes_{ab}\operatorname{Pois}(\beta_{ab,n})
-
\bigotimes_{ab}\operatorname{Pois}(b_{ab})
\right\|_{\mathrm{TV}}
\leq
\sum_{ab}|\beta_{ab,n}-b_{ab}|.
\]

Let \(\mathsf K_n(z,\cdot)\) and \(\mathsf K(z,\cdot)\) be,
respectively, the finite and limiting conditional degree kernels.
Equations \eqref{eq:category-means-final} and
\eqref{eq:Le-Cam-final} imply
\[
\int
\|\mathsf K_n(z,\cdot)-\mathsf K(z,\cdot)\|_{\mathrm{TV}}
\,\Pi_n(dz)
\longrightarrow0.
\]
The disintegration inequality gives
\begin{align*}
\|\Pi_n\mathsf K_n-\Pi\mathsf K\|_{\mathrm{TV}}
&\leq
\int
\|\mathsf K_n(z,\cdot)-\mathsf K(z,\cdot)\|_{\mathrm{TV}}
\,\Pi_n(dz)\\
&\quad+
\|\Pi_n\mathsf K-\Pi\mathsf K\|_{\mathrm{TV}}\\
&\leq
\int
\|\mathsf K_n(z,\cdot)-\mathsf K(z,\cdot)\|_{\mathrm{TV}}
\,\Pi_n(dz)
+
\|\Pi_n-\Pi\|_{\mathrm{TV}}\\
&\longrightarrow0.
\end{align*}
This proves \eqref{eq:joint-Palm-TV-final}.

\medskip
\noindent
\paragraph{Exponential endpoint-degree bounds and uniform
integrability.}
On \(A_{12}=1\), let \(N_{12,n}\) count vertices \(k\geq3\) in
\[
B_{\mathbb T^d}(X_1,R\varepsilon_n)
\cup
B_{\mathbb T^d}(X_2,R\varepsilon_n).
\]
Then
\[
D_1+D_2\leq2+2N_{12,n}.
\]
Conditionally on the endpoint state, \(N_{12,n}\) is binomial with
success probability at most \(C\varepsilon_n^d\). Hence, for every
\(a>0\),
\begin{align*}
E\left[
e^{a(D_1+D_2)}
\mid Z_{12,n},A_{12}=1
\right]
&\leq
e^{2a}
\left[
1+C\varepsilon_n^d(e^{2a}-1)
\right]^{n-2}\\
&\leq
e^{2a}
\exp\left\{
C\rho_n(e^{2a}-1)
\right\}.
\end{align*}
Consequently,
\begin{equation}
\label{eq:Palm-all-moments-final}
\sup_{n\geq n_0}
E\left[
(D_1+D_2)^m
\mid A_{12}=1
\right]
<\infty
\qquad\text{for every fixed }m<\infty.
\end{equation}
The limiting Poisson kernel has the same property because its means
are uniformly bounded.

For \(f\in\mathcal C_2\),
\[
|f|^2
\leq
C_f(1+D_1+D_2)^4.
\]
Thus \eqref{eq:Palm-all-moments-final} supplies a genuine
de la Vall\'ee--Poussin bound and hence uniform integrability.
Truncating on \(\{D_1+D_2\leq B\}\), applying the total-variation
limit of the Poisson approximation to the bounded truncation, and then sending
\(B\to\infty\) gives
\begin{equation}
\label{eq:Palm-score-final}
E\left[
f(W_1,W_2,X_1,U_{12,n},D_1,D_2)
\mid A_{12}=1
\right]
\longrightarrow
\mathsf Q(f).
\end{equation}

\medskip
\noindent
\paragraph{Expected edge count.}
By exchangeability,
\[
E\left(\frac{2E_n}{n}\right)
=
(n-1)p_n
\longrightarrow\ell.
\]

\medskip
\noindent
\paragraph{Expected triangle count.}
For \(n\geq n_0\), substitute
\[
X_2=x+\varepsilon_nu,
\qquad
X_3=x+\varepsilon_nv.
\]
On the support of the product of the first two link probabilities,
\(u,v\in K\), and the centered displacement from \(X_2\) to \(X_3\)
is \(\varepsilon_n(v-u)\). Conditional edge independence therefore
gives
\[
ET_n
=
\binom n3\varepsilon_n^{2d}b_n,
\]
where
\begin{align*}
b_n
={}&
\int
q_n(w_1,w_2;u)
q_n(w_1,w_3;v)
q_n(w_2,w_3;v-u)\\
&\quad\times
H_x(dw_1)
H_{x+\varepsilon_nu}(dw_2)
H_{x+\varepsilon_nv}(dw_3)
\,dx\,du\,dv.
\end{align*}
The integrand is supported in \(K\times K\). A telescoping
product-measure bound, \eqref{eq:q-uniform-final}, and
\eqref{eq:uniform-TV-H-final} give
\[
b_n\longrightarrow
\int
q_{12}(u)q_{13}(v)q_{23}(v-u)
\prod_{j=1}^3H_x(dw_j)\,dx\,du\,dv.
\]
Since
\[
\frac1n\binom n3\varepsilon_n^{2d}
\longrightarrow
\frac{\rho^2}{6},
\]
we conclude that
\[
\frac{ET_n}{n}\longrightarrow t.
\]

\medskip
\noindent
\paragraph{Expected connected two-star count.}
For each center \(i\),
\[
\binom{D_i}{2}
=
\sum_{\substack{j<k\\j,k\ne i}}A_{ij}A_{ik}.
\]
The two incident edge innovations are conditionally independent.
The same two-coordinate change of variables therefore gives
\[
ES_{2,n}
=
n\binom{n-1}{2}\varepsilon_n^{2d}c_n,
\]
where
\[
c_n
=
\int
q_n(w_1,w_2;u)q_n(w_1,w_3;v)
H_x(dw_1)
H_{x+\varepsilon_nu}(dw_2)
H_{x+\varepsilon_nv}(dw_3)
\,dx\,du\,dv.
\]
Again,
\[
c_n\longrightarrow
\int
q_{12}(u)q_{13}(v)
\prod_{j=1}^3H_x(dw_j)\,dx\,du\,dv.
\]
Since
\[
\binom{n-1}{2}\varepsilon_n^{2d}
\longrightarrow\frac{\rho^2}{2},
\]
we obtain
\[
\frac{ES_{2,n}}n\longrightarrow s_2.
\]

\medskip
\noindent
\paragraph{Independent-input representation and local
occupancy moments.}
Realize the graph from the independent vertex inputs
\[
Z_i=(X_i,W_i)
\]
and the independent pair innovations
\[
\xi_{ij}\sim\operatorname{Unif}(0,1),
\qquad
A_{ij}
=
\mathbf1\{
\xi_{ij}
\leq
q_n(W_i,W_j;U_{ij,n})
\}.
\]
Conditionally on any fixed collection of at most two centers, the
number of other vertices in the union of their
\(3R\varepsilon_n\)-balls is dominated by
\[
2+\operatorname{Bin}\{n,\min(1,C\varepsilon_n^d)\}.
\]
Since \(n\varepsilon_n^d=O(1)\), for every fixed \(m<\infty\),
\begin{equation}
\label{eq:local-moments-final}
\sup_n
E\left[
\{2+\operatorname{Bin}(n,\min(1,C\varepsilon_n^d))\}^m
\right]
<\infty.
\end{equation}

\medskip
\noindent
\paragraph{Influence of replacing a vertex input.}
Replace \(Z_i\) by an independent copy
\(Z_i'=(X_i',W_i')\), keeping all other inputs fixed, and denote the
modified graph by a superscript \((i)\). Define
\[
M_i
=
1+
\#\left\{
v\ne i:
X_v\in
B_{\mathbb T^d}(X_i,3R\varepsilon_n)
\cup
B_{\mathbb T^d}(X_i',3R\varepsilon_n)
\right\}.
\]
Only edges incident to \(i\) can change. Let
\[
\mathcal V_i
=
\{i\}
\cup
\{j\ne i:A_{ij}\ne A_{ij}^{(i)}\}.
\]
Then \(|\mathcal V_i|\leq M_i\). Every directed-edge summand whose
value can change is incident, in the old or new graph, to a vertex
of \(\mathcal V_i\). Furthermore, every endpoint degree appearing
in such a summand is at most \(M_i\): a changed neighbor is within
one range step of \(X_i\) or \(X_i'\), the other endpoint of an
affected edge is within two range steps, and every vertex
contributing to its degree is within three range steps.

For fixed \(f\in\mathcal C_2\), put
\[
F_n(f)
=
\sum_{a\ne b}
A_{ab}
f(W_a,W_b,X_a,U_{ab,n},D_a,D_b).
\]
At most
\[
2\sum_{v\in\mathcal V_i}
(D_v+D_v^{(i)})
\leq
4M_i^2
\]
directed-edge summands can change. The sum of the old and new
absolute values of each such summand is at most
\[
C\|f\|_{(2)}(1+M_i)^2.
\]
Consequently,
\begin{equation}
\label{eq:vertex-influence-final}
|F_n(f)-F_n^{(i)}(f)|
\leq
C\|f\|_{(2)}(1+M_i)^4.
\end{equation}
No Lipschitz property of \(f\) is used.

By \eqref{eq:local-moments-final},
\begin{equation}
\label{eq:vertex-influence-square-final}
E|F_n(f)-F_n^{(i)}(f)|^2
\leq C_f.
\end{equation}

\medskip
\noindent
\paragraph{Influence of replacing a pair innovation.}
Replace only \(\xi_{ab}\), and denote the modified statistic by a
superscript \((ab)\). Define
\[
M_{ab}
=
2+
\#\left\{
v\notin\{a,b\}:
X_v\in
B_{\mathbb T^d}(X_a,3R\varepsilon_n)
\cup
B_{\mathbb T^d}(X_b,3R\varepsilon_n)
\right\}.
\]
The edge \(A_{ab}\) can change only on
\[
\mathcal A_{ab,n}
=
\{
d_{\mathbb T^d}(X_a,X_b)\leq R\varepsilon_n
\}.
\]
On this event, only \(A_{ab},D_a,D_b\) can change. At most
\(CM_{ab}\) directed-edge summands are affected, and every degree
appearing in such a summand is at most \(M_{ab}\). Hence
\begin{equation}
\label{eq:pair-influence-final}
|F_n(f)-F_n^{(ab)}(f)|
\leq
C\|f\|_{(2)}
\mathbf1_{\mathcal A_{ab,n}}
(1+M_{ab})^3.
\end{equation}
Conditionally on \(X_a,X_b\), all moments of \(M_{ab}\) are
uniformly bounded, while
\[
P(\mathcal A_{ab,n})\leq C\varepsilon_n^d.
\]
Therefore
\begin{equation}
\label{eq:pair-influence-square-final}
E|F_n(f)-F_n^{(ab)}(f)|^2
\leq C_f\varepsilon_n^d.
\end{equation}

\medskip
\noindent
\paragraph{Explicit motif-influence bounds.}
The same geometry gives the following deterministic bounds for
vertex replacement:
\begin{align}
|2E_n-2E_n^{(i)}|
&\leq CM_i,
\label{eq:edge-vertex-influence-final}\\
|T_n-T_n^{(i)}|
&\leq C(1+M_i)^2,
\label{eq:triangle-vertex-influence-final}\\
|S_{2,n}-S_{2,n}^{(i)}|
&\leq C(1+M_i)^2.
\label{eq:star-vertex-influence-final}
\end{align}
Indeed, an affected triangle must contain \(i\), and a change in
\(S_{2,n}\) can occur only at \(i\) or at a vertex whose incident
edge to \(i\) changes.

For pair-innovation replacement,
\begin{align}
|2E_n-2E_n^{(ab)}|
&\leq
2\mathbf1_{\mathcal A_{ab,n}},
\label{eq:edge-pair-influence-final}\\
|T_n-T_n^{(ab)}|
&\leq
C\mathbf1_{\mathcal A_{ab,n}}(1+M_{ab}),
\label{eq:triangle-pair-influence-final}\\
|S_{2,n}-S_{2,n}^{(ab)}|
&\leq
C\mathbf1_{\mathcal A_{ab,n}}(1+M_{ab}).
\label{eq:star-pair-influence-final}
\end{align}
Equations \eqref{eq:local-moments-final}--\eqref{eq:star-pair-influence-final}
imply uniform \(O(1)\) vertex-influence second moments and
\(O(\varepsilon_n^d)\) pair-influence second moments for all three
count statistics.

\medskip
\noindent
\paragraph{Efron--Stein concentration and the count laws of large numbers.}
Apply Efron--Stein to the complete independent input family
\[
\{Z_i:1\leq i\leq n\}
\cup
\{\xi_{ab}:1\leq a<b\leq n\}.
\]
For each
\[
F\in\{2E_n,T_n,S_{2,n},F_n(f)\},
\]
The three influence bounds give
\begin{align*}
\operatorname{Var}(F)
&\leq
\frac12\sum_{i=1}^n
E|F-F^{(i)}|^2
+
\frac12\sum_{a<b}
E|F-F^{(ab)}|^2\\
&\leq
C_f\{n+n^2\varepsilon_n^d\}
=
O(n).
\end{align*}
It follows that
\[
E\left|
\frac{F-EF}{n}
\right|^2
=
O(n^{-1}).
\]
Combining this bound with the expected motif counts yields
\[
\left(
\frac{2E_n}{n},
\frac{T_n}{n},
\frac{S_{2,n}}{n}
\right)
\overset{L^2}\longrightarrow
(\ell,t,s_2).
\]
Coordinatewise \(L^2\) convergence also implies joint \(L^2\)
convergence for every fixed finite vector.

\medskip
\noindent
\paragraph{The mean of a directed-edge score.}
Evenness of \(k\) gives
\[
q_n(W_i,W_j;U_{ij,n})=q_n(W_j,W_i;U_{ji,n})
\]
outside the Haar-null centered-representative boundary.  Thus the
unordered-edge construction is label-exchangeable, and
\begin{align}
E\{\mathsf M_n(f)\}
&=
(n-1)
E\left[
A_{12}
f(W_1,W_2,X_1,U_{12,n},D_1,D_2)
\right]
\nonumber\\
&=
(n-1)p_n
E\left[
f(W_1,W_2,X_1,U_{12,n},D_1,D_2)
\mid A_{12}=1
\right].
\label{eq:edge-score-mean-final}
\end{align}
By the edge-probability estimate,
\[
(n-1)p_n\longrightarrow\ell,
\]
and by the exponential endpoint-degree bounds,
\[
E\{\mathsf M_n(f)\}
\longrightarrow
\ell\mathsf Q(f).
\]

\medskip
\noindent
\paragraph{Scalar \(L^2\) convergence and weak
random-measure convergence.}
Since \(F_n(f)=n\mathsf M_n(f)\), the Efron--Stein concentration bound gives
\[
\operatorname{Var}\{\mathsf M_n(f)\}
=
\frac{\operatorname{Var}\{F_n(f)\}}{n^2}
=
O(n^{-1}).
\]
Consequently,
\begin{align*}
E\left|
\mathsf M_n(f)-\ell\mathsf Q(f)
\right|^2
&=
\operatorname{Var}\{\mathsf M_n(f)\}\\
&\quad+
\left|
E\{\mathsf M_n(f)\}-\ell\mathsf Q(f)
\right|^2
\longrightarrow0.
\end{align*}
This proves \eqref{eq:M-f-final}. The same argument applied to a
fixed finite collection of scores and the count statistics proves
the asserted joint convergence.

It remains to justify the weak random-measure statement. Let
\[
\mathsf R_n
=
\mathcal L
\left(
W_1,W_2,X_1,U_{12,n},D_1,D_2
\mid A_{12}=1
\right).
\]
By the multivariate Poisson approximation,
\[
\|\mathsf R_n-\mathsf Q\|_{\mathrm{TV}}\longrightarrow0.
\]
Fix \(a,b>0\).  Choose \(\eta>0\) small enough that
\(2(\ell+1)\eta/a<b\), and choose a compact set
\(B_\eta\subset\mathsf S\) satisfying
\[
\mathsf Q(B_\eta^c)<\eta.
\]
Then, for all large \(n\),
\[
\mathsf R_n(B_\eta^c)
\leq
\mathsf Q(B_\eta^c)
+
\|\mathsf R_n-\mathsf Q\|_{\mathrm{TV}}
\leq2\eta.
\]
By exchangeability,
\[
E\{\mathsf M_n(B_\eta^c)\}
=
(n-1)p_n\,\mathsf R_n(B_\eta^c).
\]
Because \((n-1)p_n\to\ell\), Markov's inequality proves
\[
\limsup_{n\to\infty}
P\{\mathsf M_n(B_\eta^c)>a\}
\leq \frac{2(\ell+1)\eta}{a}<b.
\]
This is asymptotic tightness in probability of
\(\{\mathsf M_n\}\).

For every fixed bounded continuous \(g\),
\[
\mathsf M_n(g)
\overset{\mathbb P}\longrightarrow
\ell\mathsf Q(g).
\]
Applying this to a countable convergence-determining class, together
with the preceding asymptotic tightness and the subsequence
characterization of convergence in probability, yields
\[
\mathsf M_n
\Rightarrow
\ell\mathsf Q
\]
weakly in probability.

The following identities record explicitly the unbounded score
limits used below:
\begin{align}
\mathsf M_n(1)
&=
\frac{2E_n}{n},
\label{eq:M-one-final}\\
\mathsf M_n(d_1)
=
\mathsf M_n(d_2)
&=
\frac1n\sum_{i=1}^nD_i^2,
\label{eq:M-first-final}\\
\mathsf M_n(d_1^2)
=
\mathsf M_n(d_2^2)
&=
\frac1n\sum_{i=1}^nD_i^3,
\label{eq:M-second-final}\\
\mathsf M_n(d_1d_2)
&=
\frac1n\sum_{i\ne j}A_{ij}D_iD_j.
\label{eq:M-cross-final}
\end{align}
Because these score functions belong to \(\mathcal C_2\), all four
identities have their corresponding joint \(L^2\) limits.

\medskip
\noindent
\paragraph{Normalization and population endpoint moments.}
Let
\[
\mu(dw,dx)=H_x(dw)\,dx.
\]
Fubini's theorem gives
\begin{equation}
\label{eq:ell-Lambda-final}
\ell
=
\int\Lambda(w,x)\,\mu(dw,dx)
\end{equation}
and
\begin{equation}
\label{eq:s2-Lambda-final}
s_2
=
\frac12
\int\Lambda(w,x)^2\,\mu(dw,dx).
\end{equation}
Since \(\mu\) is a probability measure, Jensen's inequality gives
\[
s_2\geq\frac{\ell^2}{2}>0.
\]

The first-endpoint marginal of \(\Pi\) is
\[
\Pi_1(dw,dx)
=
\frac{\Lambda(w,x)}{\ell}
H_x(dw)\,dx.
\]
Therefore
\begin{equation}
\label{eq:EPi-L-final}
E_\Pi L_1
=
\frac1\ell
\int\Lambda(w,x)^2H_x(dw)\,dx
=
\frac{2s_2}{\ell}
>0.
\end{equation}
Similarly,
\begin{equation}
\label{eq:EPi-G-final}
E_\Pi G
=
\frac{6t}{\ell},
\end{equation}
and hence
\[
C
=
\frac{3t}{s_2}
=
\frac{E_\Pi G}{E_\Pi L_1}.
\]

Evenness of \(k\) implies that \(\Pi\), and then \(\mathsf Q\), is
invariant under
\[
(w_1,w_2,x,u,d_1,d_2)
\longmapsto
(w_2,w_1,x,-u,d_2,d_1).
\]
Indeed,
\[
q(w_2,w_1;-u)=q(w_1,w_2;u),
\]
and the substitution \(r=v+u\) gives
\[
\Gamma(w_2,w_1,x;-u)
=
\Gamma(w_1,w_2,x;u).
\]
Thus the two endpoint marginals coincide. Put
\[
m=E_\Pi L_1=E_\Pi L_2.
\]

Conditionally on the Palm state \(z\),
\[
E(D_j^\star\mid z)=1+L_j,
\qquad
\operatorname{Var}(D_j^\star\mid z)=L_j,
\]
and
\[
\operatorname{Cov}(D_1^\star,D_2^\star\mid z)=G.
\]
Equivalently,
\begin{align*}
E\{(D_1^\star)^2\mid z\}
&=
1+3L_1+L_1^2,\\
E(D_1^\star D_2^\star\mid z)
&=
(1+L_1)(1+L_2)+G.
\end{align*}
The laws of total variance and covariance therefore give
\begin{align}
V
:=
\operatorname{Var}_{\mathsf Q}(D_1^\star)
&=
\operatorname{Var}_{\Pi}(L_1)+E_\Pi L_1,
\label{eq:population-V-final}\\
K_e
:=
\operatorname{Cov}_{\mathsf Q}(D_1^\star,D_2^\star)
&=
\operatorname{Cov}_{\Pi}(L_1,L_2)+E_\Pi G.
\label{eq:population-K-final}
\end{align}
By \eqref{eq:EPi-L-final},
\[
V\geq E_\Pi L_1=\frac{2s_2}{\ell}>0.
\]
Thus \(r=K_e/V\) is well defined and agrees with
\eqref{eq:r-final}.

\medskip
\noindent
\paragraph{Empirical assortativity and random denominators.}
On \(E_n>0\),
\[
\widehat{\mathsf Q}_n
=
\frac{\mathsf M_n}{\mathsf M_n(1)}
\]
is exactly the empirical distribution of a uniformly chosen
realized directed edge. Because each undirected edge occurs in both
orientations,
\begin{align}
\int d_1\,d\widehat{\mathsf Q}_n
=
\int d_2\,d\widehat{\mathsf Q}_n
&=
\frac{\sum_iD_i^2}{2E_n},
\label{eq:empirical-edge-mean-final}\\
\int d_1^2\,d\widehat{\mathsf Q}_n
=
\int d_2^2\,d\widehat{\mathsf Q}_n
&=
\frac{\sum_iD_i^3}{2E_n},
\label{eq:empirical-edge-second-final}\\
\int d_1d_2\,d\widehat{\mathsf Q}_n
&=
\frac{\sum_{i\ne j}A_{ij}D_iD_j}{2E_n}.
\label{eq:empirical-edge-cross-final}
\end{align}
Hence the two empirical endpoint means and variances are exactly
equal, not merely asymptotically equal.

From the concentration bounds and the scalar \(L^2\) convergence,
\[
\mathsf M_n(1)=\frac{2E_n}{n}
\overset{\mathbb P}\longrightarrow\ell>0.
\]
Slutsky's theorem and the joint score limits therefore give
\begin{align*}
\overline D_{e,n}
&\overset{\mathbb P}\longrightarrow
E_{\mathsf Q}D_1^\star,\\
\widehat V_{e,n}
&\overset{\mathbb P}\longrightarrow
\operatorname{Var}_{\mathsf Q}(D_1^\star)=V>0,\\
\widehat\kappa_{e,n}
&\overset{\mathbb P}\longrightarrow
\operatorname{Cov}_{\mathsf Q}(D_1^\star,D_2^\star)=K_e.
\end{align*}
Consequently,
\[
\widehat r_n
\overset{\mathbb P}\longrightarrow
\frac{K_e}{V}
=
r.
\]
Replacing each endpoint degree by its excess degree \(D_i-1\)
subtracts the same constant from both endpoint variables and
therefore leaves \(\widehat r_n\) unchanged exactly.

Finally,
\[
\frac{2E_n}{n}\to\ell>0,
\qquad
\frac{S_{2,n}}n\to s_2>0,
\qquad
\widehat V_{e,n}\to V>0
\]
in probability. Therefore
\[
P(E_n=0)\to0,\qquad
P(S_{2,n}=0)\to0,\qquad
P(\widehat V_{e,n}=0)\to0.
\]
The continuous mapping theorem yields
\[
\widehat C_n
=
\frac{3T_n/n}{S_{2,n}/n}
\overset{\mathbb P}\longrightarrow
\frac{3t}{s_2}
=
C
\]
and the preceding convergence of \(\widehat r_n\). Since the count
LLNs and the required finite list of edge-score LLNs were proved
jointly, this establishes \eqref{eq:summary-vector-final}.
\end{proof}

\begin{remark}[Scope]
The theorem is pointwise in \(\vartheta\). It does not provide a
uniform law of large numbers over a parameter set, a root-\(n\)
centering expansion, a central limit theorem, covariance
nondegeneracy, consistent covariance estimation, or a GMM inference
result. Each of those properties requires a separate theorem.
\end{remark}

\clearpage
\section{Certified identification: exact reduction and the global
theorem}
\label{sec:si-s3}

\subsection{Exact reduction and certified evaluator bounds}
\label{subsec:certified-three-parameter-region}

This result concerns the limiting population map in the rescaled tangent
coordinate.  Let $\mathbb T=\mathbb R/\mathbb Z$, let
$X\sim\operatorname{Haar}(\mathbb T)$, and let
\[
 z=\frac{(X_j-X_i)_{\mathbb T}}{\varepsilon_n}\in\mathbb R
\]
denote the tangent displacement in the two-scale limit.  Thus, the range
parameter below acts on $z\in\mathbb R$; it is not an unscaled distance on
the unit torus.  Fix $d=1$, $\rho=1$, set the popularity variable $W=U$, and
put
\begin{align*}
 a(u)&=\sqrt{3}(2u-1),
 &\varphi(x)&=\sqrt{2}\cos(2\pi x),\\
 h_b(u\mid x)&=1+b\,a(u)\varphi(x),
 &b&=\sqrt\psi,
\end{align*}
for $u\in[0,1]$ and $x\in\mathbb T$.  Conditional on a location $x$, marks
are independent with density $h_b(\cdot\mid x)$.  For
\[
 \theta=(\lambda,\eta,\psi)
 \in\Theta:=(0,\infty)^2\times(0,1/6),
\]
define
\[
 k_\eta(z)=\mathbf 1\{|z|\leq\eta\},
 \qquad
 H_\lambda(u,v)=1-e^{-\lambda uv},
\]
so the tangent edge link is
\[
 q_\theta(u,v;z)
 =1-\exp\{-\lambda uv k_\eta(z)\}
 =H_\lambda(u,v)\mathbf 1\{|z|\leq\eta\}.
\]
The kernel height, popularity scale, harmonic, and $\rho$ are therefore fixed.

For $2\leq q\leq4$, define the location-averaged mark density
\[
 \omega_q(u_1,\ldots,u_q;\psi)
 =\int_{\mathbb T}\prod_{j=1}^q h_{\sqrt\psi}(u_j\mid x)\,dx.
\]
Writing $a_i=a(u_i)$, the harmonic identities
$\int_{\mathbb T}\varphi\,dx=\int_{\mathbb T}\varphi^3\,dx=0$,
$\int_{\mathbb T}\varphi^2\,dx=1$, and
$\int_{\mathbb T}\varphi^4\,dx=3/2$ give
\begin{equation}
 \label{eq:omega-expansion-certified}
 \omega_q
 =1+\psi\sum_{1\leq i<j\leq q}a_i a_j
 +\mathbf 1\{q=4\}\frac32\psi^2a_1a_2a_3a_4.
\end{equation}
With all integrals below taken over the displayed unit cube, set
\begin{align*}
 E&=\int_{[0,1]^2}H_{12}\omega_2\,du_{1:2},\\
 S&=\int_{[0,1]^3}H_{12}H_{13}\omega_3\,du_{1:3},\\
 T&=\int_{[0,1]^3}H_{12}H_{13}H_{23}\omega_3\,du_{1:3},\\
 P&=\int_{[0,1]^4}H_{12}H_{13}H_{24}\omega_4\,du_{1:4},\\
 Q&=\int_{[0,1]^4}H_{12}H_{13}H_{14}\omega_4\,du_{1:4},
\end{align*}
where $H_{ij}=H_\lambda(u_i,u_j)$.  The limiting hard-range target map is
\begin{equation}
 \label{eq:certified-target-map}
 \ell=2\eta E,
 \qquad
 C=\frac{3T}{4S},
 \qquad
 r=
 \frac{2\eta(PE-S^2)+(3/4)TE}
      {2\eta(QE-S^2)+SE}.
\end{equation}

\begin{siproposition}[Exact reduction, denominator positivity, and
evaluator truncation bounds]
\label{prop:exact-reduction}
Let $M_3:\Theta\to\mathbb R^3$ be the raw-coordinate map
\[
 M_3(\lambda,\eta,\psi)
 =\{\ell(\lambda,\eta,\psi),
      C(\lambda,\eta,\psi),
      r(\lambda,\eta,\psi)\}^{\trans}.
\]
Then: \textup{(i)} $h_{\sqrt\psi}(\cdot\mid x)$ is a probability density
with uniform marginals for every $\psi<1/6$; \textup{(ii)} the three
observables reduce exactly to \eqref{eq:certified-target-map}, the
denominator of $r$ is strictly positive on $\Theta$, and $M_3$ is real
analytic on $\Theta$; \textup{(iii)} the degree-$N$ truncated motif
polynomials \eqref{eq:exact-truncated-motif} approximate the exact
motifs with the explicit remainder bounds
\eqref{eq:product-value-tail}--\eqref{eq:product-psi-tail}, which hold
for every truncation order $N$ with the tail constants
$\epsilon_0,\epsilon_1$ evaluated at the relevant order and range.
These are the inputs consumed by the global certificate of the next
subsection.
\end{siproposition}

\begin{proof}
The proof proceeds through the exact probabilistic reduction and the
truncation-error analysis.

\emph{*Validity of the mark model.*}
For $\psi<1/6$,
\[
 |b\,a(u)\varphi(x)|\leq\sqrt{6\psi}<1.
\]
Hence $0<h_b(u\mid x)<2$.  Moreover,
\[
 \int_0^1h_b(u\mid x)\,du=1,
 \qquad
 \int_{\mathbb T}h_b(u\mid x)\,dx=1,
\]
because $\int_0^1a(u)\,du=0$ and
$\int_{\mathbb T}\varphi(x)\,dx=0$.  Thus $h_b(\cdot\mid x)$ is a density
and $U$ is marginally uniform.  Expanding the product defining $\omega_q$
and using the four harmonic moments stated above proves
\eqref{eq:omega-expansion-certified}.

\emph{*Exact hard-range and Palm reduction.*}
The edge displacement set has length $2\eta$; the two-edge wedge set has area
$4\eta^2$; and the triangle set
\[
 \{(z_1,z_2):|z_1|\leq\eta,\ |z_2|\leq\eta,
                    \ |z_2-z_1|\leq\eta\}
\]
has area $3\eta^2$.  Therefore
\[
 \ell=2\eta E,
 \qquad
 t=\frac16(3\eta^2)T,
 \qquad
 s_2=\frac12(4\eta^2)S,
\]
and the transitivity target $C=3t/s_2$ equals $3T/(4S)$.

For the assortativity reduction, define
\begin{align*}
 g_x(u)&=\int_0^1H_\lambda(u,w)h_b(w\mid x)\,dw,\\
 R_x(u,v)&=\int_0^1H_\lambda(u,w)H_\lambda(v,w)
                    h_b(w\mid x)\,dw,
\end{align*}
and let the edge--Palm distribution of $(x,z,u_1,u_2)$ have density
\[
 \frac{\mathbf 1\{|z|\leq\eta\}
 H_{12}h_b(u_1\mid x)h_b(u_2\mid x)}{2\eta E}
\]
with respect to $dx\,dz\,du_1\,du_2$.  Its $(x,u_1,u_2)$ marginal has
density $H_{12}h_b(u_1\mid x)h_b(u_2\mid x)/E$.  Direct integration gives
\begin{equation}
 \label{eq:edge-palm-motif-identities}
 \mathbb E_{\mathrm e}g_X(U_1)=\frac SE,
 \quad
 \mathbb E_{\mathrm e}\{g_X(U_1)g_X(U_2)\}=\frac PE,
 \quad
 \mathbb E_{\mathrm e}g_X(U_1)^2=\frac QE,
 \quad
 \mathbb E_{\mathrm e}R_X(U_1,U_2)=\frac TE.
\end{equation}
Under the edge--Palm law, the tangent displacement is uniform on the interval
$[-\eta,\eta]$.  Its mean neighborhood-overlap length is
\[
 \frac1{2\eta}\int_{-\eta}^{\eta}(2\eta-|z|)\,dz
 =\frac{3\eta}{2}.
\]
The limiting mixed-Poisson edge--Palm representation has endpoint intensities
$\Lambda_i=2\eta g_X(U_i)$ and, conditionally on $(X,U_1,U_2,z)$,
common-neighbor covariance
$(2\eta-|z|)R_X(U_1,U_2)$.  Its conditional mean over $z$ is
$(3\eta/2)R_X(U_1,U_2)$.  The laws of total covariance and total variance,
together with \eqref{eq:edge-palm-motif-identities}, therefore yield
\begin{align*}
 \operatorname{Cov}_{\mathrm e}(D_1^-,D_2^-)
 &=\frac{4\eta^2}{E^2}(PE-S^2)+\frac{3\eta}{2E}T,\\
 \operatorname{Var}_{\mathrm e}(D_1^-)
 &=\frac{4\eta^2}{E^2}(QE-S^2)+\frac{2\eta}{E}S,
\end{align*}
where $D_i^-$ is the limiting excess degree at endpoint $i$.  Their ratio is
the expression for $r$ in \eqref{eq:certified-target-map}.  Thus $r$ is the
Pearson correlation of endpoint degrees under edge sampling.  Subtracting
one from both endpoint degrees leaves this correlation unchanged.

The preceding display also proves denominator positivity without numerical
assistance.  Indeed, $E>0$ and $S>0$, and
\[
 QE-S^2=E^2\operatorname{Var}_{\mathrm e}\{g_X(U_1)\}\geq0.
\]
Consequently, $2\eta(QE-S^2)+SE>0$ throughout $\Theta$.
Each motif is real analytic in $\lambda$ and polynomial in $\psi$ after the
fixed-domain reduction, while its dependence on $\eta$ in
\eqref{eq:certified-target-map} is rational with a positive denominator.
Thus $M_3$ is real analytic on $\Theta$.  Notice that this argument rescales
the hard-range domains; pointwise differentiation of
$\mathbf 1\{|z|\leq\eta\}$ would not justify differentiability in $\eta$.

\emph{*Truncation bounds for the motif integrals.*}
Fix a truncation order $N$ and a range bound $L$ with $\lambda\le L$,
and write, for $x\geq0$,
\[
 1-e^{-x}=P_N(x)+R_N(x),
 \qquad
 P_N(x)=\sum_{n=1}^{N}\frac{(-1)^{n+1}x^n}{n!}.
\]
The displays below take the illustrative instance $N=22$,
$L=301/100$; the global verifier of the next subsection uses the same
bounds at $N=40$, $L=5$.  Since $u,v\in[0,1]$,
$0\leq x=\lambda uv\leq L$.  Taylor's theorem,
applied separately to the function and its $\lambda$ derivative, gives
\begin{align}
 |R_N(\lambda uv)|
 &\leq\epsilon_0:=\frac{L^{23}}{23!},
 \label{eq:tail-epsilon-zero}\\
 |\partial_\lambda R_N(\lambda uv)|
 &=uv\left|e^{-\lambda uv}
   -\sum_{k=0}^{21}\frac{(-\lambda uv)^k}{k!}\right|
 \leq\epsilon_1:=\frac{L^{22}}{22!}.
 \label{eq:tail-epsilon-one}
\end{align}

For $d\geq0$, put
\[
 p_d=\int_0^1u^d\,du=\frac1{d+1},
 \qquad
 \bar q_d=\frac1{\sqrt3}\int_0^1u^da(u)\,du
 =\frac{d}{(d+1)(d+2)}.
\]
Consider any one of the five motif graphs, with $v\leq4$ vertices and
$m\leq3$ edges.  After expanding the $m$ copies of $P_N$, attach an exponent
$n_e\in\{1,\ldots,22\}$ to each edge and put
$d_i=\sum_{e\ni i}n_e$.  Exact monomial integration gives
\begin{align*}
 \Omega_{v,\mathbf d}(\psi)
 &=\prod_{i=1}^vp_{d_i}
 +3\psi\sum_{i<j}\bar q_{d_i}\bar q_{d_j}
              \prod_{k\notin\{i,j\}}p_{d_k}\\
 &\quad+
 \mathbf1\{v=4\}\frac{27}{2}\psi^2
              \prod_{i=1}^v\bar q_{d_i}.
\end{align*}
Hence the truncated motif is the rational polynomial
\begin{equation}
 \label{eq:exact-truncated-motif}
 \sum_{1\leq n_e\leq22,\ e\in\mathcal E(G)}
 \frac{(-1)^{\sum_en_e+m}\lambda^{\sum_en_e}}
      {\prod_en_e!}\,
 \Omega_{v,\mathbf d}(\psi).
\end{equation}

We next record the complete remainder propagation.  Write $H_j=P_j+R_j$
for the factors in a motif.  Since $0\leq H_j\leq1$ and
$|P_j|\leq1+\epsilon_0$, telescoping gives
\begin{equation}
 \label{eq:product-value-tail}
 \left|\prod_{j=1}^mH_j-\prod_{j=1}^mP_j\right|
 \leq m\epsilon_0(1+\epsilon_0)^{m-1}.
\end{equation}
For the derivative, the exact identity
\begin{align*}
 \partial_\lambda\!\left(\prod_jH_j-\prod_jP_j\right)
 &=\sum_{i=1}^m(\partial_\lambda R_i)\prod_{j\ne i}P_j\\
 &\quad+
 \sum_{i=1}^m(\partial_\lambda H_i)
 \left(\prod_{j\ne i}H_j-\prod_{j\ne i}P_j\right)
\end{align*}
and $|\partial_\lambda H_i|\leq1$ imply
\begin{equation}
 \label{eq:product-lambda-tail}
 \left|\partial_\lambda\!\left(\prod_jH_j-\prod_jP_j\right)\right|
 \leq m\left\{
 \epsilon_1(1+\epsilon_0)^{m-1}
 +(m-1)\epsilon_0(1+\epsilon_0)^{(m-2)_+}
 \right\},
\end{equation}
where $x_+=\max\{x,0\}$.  Multiplication by $2^v$ bounds the corresponding
integrated errors because $0<h_b<2$ and hence $0\leq\omega_v\leq2^v$.
Finally,
\[
 \sup|\partial_\psi\omega_v|\leq B_v,
 \qquad (B_2,B_3,B_4)=(3,9,23).
\]
For $v=4$, the bound follows from
$6\cdot3+3\psi(\sqrt3)^4<18+9/2<23$; the cases $v=2,3$ are immediate.
Thus the integrated $\psi$-derivative error is at most
\begin{equation}
 \label{eq:product-psi-tail}
 B_vm\epsilon_0(1+\epsilon_0)^{m-1}.
\end{equation}

For general $N$ and $L$ the tail constants are
$\epsilon_0=L^{N+1}/(N+1)!$ and $\epsilon_1=L^{N}/N!$, and the
argument above is unchanged.  This proves (iii) and completes the
proof.
\end{proof}


\subsection{Certified global identification}
\label{subsec:certified-global-region}

This section proves \emph{global injectivity} and uniform Jacobian
nondegeneracy of the three-moment map on
\begin{equation}
 \label{eq:big-box}
 B \;=\; [1,5]\times[5,15]\times[0,\tfrac{3}{20}]
 \;\subset\;(0,\infty)^2\times[0,\tfrac16),
\end{equation}
with coordinates $\theta=(\lambda,\eta,\psi)$.  The box includes the
independence boundary $\psi=0$; its $\ell$-image spans mean degrees from
about $2.03$ to about $17.04$.

\subsubsection{Model reduction}

The submodel, the motif integrals
\[
 E=\textstyle\int H_{12}\omega_2,\quad
 S=\int H_{12}H_{13}\omega_3,\quad
 T=\int H_{12}H_{13}H_{23}\omega_3,\quad
 P=\int H_{12}H_{13}H_{24}\omega_4,\quad
 Q=\int H_{12}H_{13}H_{14}\omega_4,
\]
with $H_{ij}=1-e^{-\lambda u_iu_j}$ and
$\omega_q=1+\psi\sum_{i<j}a_ia_j+\mathbf 1\{q=4\}\tfrac32\psi^2a_1a_2a_3a_4$,
and the exact hard-range and edge--Palm reduction
\begin{equation}
 \label{eq:target-map-v2}
 \ell=2\eta E,
 \qquad
 C=\frac{3T}{4S},
 \qquad
 r=\frac{2\eta(PE-S^2)+\tfrac34 TE}{2\eta(QE-S^2)+SE},
\end{equation}
are exactly as in the preceding subsection; the derivation, the
mixed-Poisson edge--Palm representation behind $r$, the identity
$QE-S^2=E^2\operatorname{Var}_{\mathrm e}\{g_X(U_1)\}\ge0$, and real
analyticity of $M_3=(\ell,C,r)$ on $\Theta$ are proved in
Proposition~\ref{prop:exact-reduction} and are not repeated.

Two structural facts drive the global result.  First, $C$ in
\eqref{eq:target-map-v2} \emph{does not depend on} $\eta$: the hard-range
areas cancel in the ratio $3t/s_2$.  Second, for $E>0$ the substitution
$2\eta=\ell/E$ eliminates $\eta$ from $r$: multiplying numerator and
denominator of \eqref{eq:target-map-v2} by $E>0$,
\begin{equation}
 \label{eq:reduced-r}
 r=\widetilde r(\lambda,\psi;\ell)
 :=\frac{\ell\,(PE-S^2)+\tfrac34\,TE^{2}}{\ell\,(QE-S^2)+S E^{2}}
 \qquad\text{whenever }2\eta E=\ell .
\end{equation}

\subsubsection{A monotone--triangular injectivity criterion}

Write $B_2=[1,5]\times[0,\tfrac{3}{20}]$ for the $(\lambda,\psi)$
rectangle and
\[
 \operatorname{den}'(\lambda,\psi;\ell)=\ell\,(QE-S^2)+SE^2 ,
 \qquad
 \Delta(\lambda,\psi;\ell)
 =\partial_\psi\widetilde r\,\partial_\lambda C
 -\partial_\lambda\widetilde r\,\partial_\psi C .
\]

\begin{silemma}[Monotone--triangular global injectivity]
\label{lem:mono-triangular}
Let $[\ell_-,\ell_+]$ contain the $\ell$-image
$\{2\eta E(\lambda,\psi):(\lambda,\eta,\psi)\in B\}$.  Assume that on
$B_2\times[\ell_-,\ell_+]$,
\begin{enumerate}[label=\textup{(I\arabic*)}]
\item $S>0$ and $\operatorname{den}'>0$;
\item $\partial_\lambda C>0$ on $B_2$;
\item $\partial_\psi C>0$ on $B_2$;
\item $\Delta>0$.
\end{enumerate}
Then $M_3$ is injective on $B$.
\end{silemma}

\begin{proof}
Let $\theta^{(i)}=(\lambda_i,\eta_i,\psi_i)\in B$, $i=1,2$, satisfy
$M_3(\theta^{(1)})=M_3(\theta^{(2)})=(\ell,C_0,r_0)$.  Then $\ell$ lies
in the $\ell$-image of $B$, hence in $[\ell_-,\ell_+]$, and since $E>0$,
each point satisfies $2\eta_i=\ell/E(\lambda_i,\psi_i)$ and, by
\eqref{eq:reduced-r}, $r_0=\widetilde r(\lambda_i,\psi_i;\ell)$.  It
therefore suffices to prove $(\lambda_1,\psi_1)=(\lambda_2,\psi_2)$;
$\eta_i=\ell/\{2E(\lambda_i,\psi_i)\}$ is then determined.

By (I2), for each fixed $\psi$ the map $\lambda\mapsto C(\lambda,\psi)$
is strictly increasing on $[1,5]$, so the equation $C(\lambda,\psi)=C_0$
has at most one root $\lambda_{C_0}(\psi)\in[1,5]$, and a root exists
exactly when $C(1,\psi)\le C_0\le C(5,\psi)$.  Let
$D=\{\psi\in[0,\tfrac3{20}]:C(1,\psi)\le C_0\le C(5,\psi)\}$.  By (I3)
the endpoint maps $\psi\mapsto C(1,\psi)$ and $\psi\mapsto C(5,\psi)$
are strictly increasing, so
$\{\psi:C(1,\psi)\le C_0\}$ is a closed lower interval,
$\{\psi:C(5,\psi)\ge C_0\}$ a closed upper interval, and $D$ is a closed
interval.  On $D$ the root map $\psi\mapsto\lambda_{C_0}(\psi)$ is
continuous (uniqueness of roots of a continuous strictly monotone family)
and, on the interior of $D$ where $\lambda_{C_0}\in(1,5)$, continuously
differentiable with
$\lambda_{C_0}'=-\partial_\psi C/\partial_\lambda C$ by the implicit
function theorem and (I2).

Define $\Phi(\psi)=\widetilde r\{\lambda_{C_0}(\psi),\psi;\ell\}$ on
$D$; $\Phi$ is continuous on $D$, and at interior points,
\[
 \Phi'(\psi)
 =\partial_\psi\widetilde r
 +\partial_\lambda\widetilde r\;\lambda_{C_0}'(\psi)
 =\frac{\partial_\psi\widetilde r\,\partial_\lambda C
        -\partial_\lambda\widetilde r\,\partial_\psi C}
       {\partial_\lambda C}
 =\frac{\Delta}{\partial_\lambda C}>0
\]
by (I2) and (I4).  If $\lambda_{C_0}$ touches an endpoint
$\lambda\in\{1,5\}$ on a subinterval of $D$, then on that subinterval
$C(\lambda_{C_0}(\psi),\psi)=C_0$ forces, by (I3), that the subinterval
is a single point; hence $\Phi'$ exists and is positive off a finite set,
and $\Phi$ is strictly increasing on $D$.  Both $\psi_1,\psi_2\in D$ with
$\Phi(\psi_1)=r_0=\Phi(\psi_2)$, so $\psi_1=\psi_2$; then
$\lambda_1=\lambda_{C_0}(\psi_1)=\lambda_2$, and finally
$\eta_1=\eta_2$.
\end{proof}

\subsubsection{The certified theorem}

\begin{sitheorem}[Certified global identification on $B$]
\label{thm:certified-global}
Let $M_3=(\ell,C,r):\Theta\to\mathbb R^3$ be the map
\eqref{eq:target-map-v2} and let $B$ be the box \eqref{eq:big-box}.
Then:
\begin{enumerate}[label=\textup{(\roman*)}]
\item \textup{(Global injectivity)} $M_3$ is injective on $B$; in
particular the parameter $(\lambda,\eta,\psi)$, including the boundary
value $\psi=0$, is uniquely determined on $B$ by the population mean
degree, transitivity, and edge endpoint-degree correlation.
\item \textup{(Uniform nondegeneracy)} uniformly on $B$,
\[
 \det DM_3<-\tfrac1{64},
 \qquad
 \sigma_{\min}\{DM_3(\theta)\}\ \ge\ \tfrac1{100}.
\]
Hence $M_3$ is a real-analytic diffeomorphism of the interior of $B$
onto its image, and the inverse is Lipschitz on compacts of the image.
\item \textup{(Compensated solvability)} on the interior of $B$, along
any profile that holds $(\ell,C)$ fixed, the derivative of $r$ with
respect to $\psi$ equals $\Delta/\partial_\lambda C>0$; endpoint mixing
strictly increases with dependence strength after compensating intensity
and range.
\end{enumerate}
\end{sitheorem}

\begin{proof}
The archived exact interval verifier
(Section~\ref{subsec:verifier-v2}; replication package)
certifies the following inequalities, each uniformly on the displayed
region, by adaptive bisection with inclusion-isotone interval
arithmetic:
\begin{center}
\begin{tabular}{@{}lll@{}}
\toprule
Stage & Region & Certified gates\\
\midrule
0 & $B_2$ &
$0.20340\,{\le}\,E\,{\le}\,0.56805$, hence
$\ell\text{-image}(B)\subset[2.0340,\,17.0415]$\\
1 & $B_2\times[2.0340,17.0415]$ &
$S>0$, $\operatorname{den}'>\tfrac1{512}$,
$\partial_\lambda C>\tfrac1{50}$,
$\partial_\psi C>\tfrac1{12}$,
$\Delta>\tfrac1{800}$\\
2 & $B$ &
$\det DM_3<-\tfrac1{64}$,
$\sigma_{\min}\ge\tfrac1{100}$\\
\bottomrule
\end{tabular}
\end{center}
Stage~0 and the identity $\ell=2\eta E$ show that the $\ell$-image of
$B$ is contained in the Stage~1 $\ell$-range.  Stage~1 supplies
hypotheses (I1)--(I4) of Lemma~\ref{lem:mono-triangular}, which yields
(i).  For (ii), the verifier certifies, on every leaf box of an adaptive
partition of $B$, an enclosure of $\det DM_3$ with upper endpoint below
$-1/64$ and the bound
\[
 \sigma_{\min}(J)
 =\frac{|\det J|}{\sigma_1\sigma_2}
 \ \ge\
 \frac{|\det J|}{\|\Lambda^2 J\|_F},
 \qquad
 \|\Lambda^2J\|_F^2=\sum_{\text{all }2\times2\text{ minors }\mu}\mu^2
 =\sigma_1^2\sigma_2^2+\sigma_1^2\sigma_3^2+\sigma_2^2\sigma_3^2 ,
\]
whose certified per-box values exceed $1/100$; the identity holds
because the squared Frobenius norm of the second compound matrix is the
second elementary symmetric function of the squared singular values.
Real-analytic local invertibility at every interior point plus global
injectivity give the diffeomorphism statement; Lipschitz continuity of
the inverse on compacts follows from
$\sigma_{\min}\ge1/100$ and the mean value inequality applied to
$M_3^{-1}$ along paths in the open image.  Part (iii) is the identity
$\Phi'=\Delta/\partial_\lambda C$ from the proof of
Lemma~\ref{lem:mono-triangular}, certified positive by Stage~1 (its
floating-point value at the center design $(3,9,0.0835)$ is about
$1.21$; over $B$ it varies but remains strictly positive).
\end{proof}

\begin{remark}[Scope]
The theorem fixes the popularity marginal, the one-harmonic
family, the exact exponential link, the hard-range kernel form, kernel
height, $\rho=1$, and $d=1$; changing any of these requires a new
certificate.  Within that submodel, identification is global on a box
whose sides span a factor of $5$ in intensity, a factor of $3$ in
range, and the full dependence interval $[0,0.15]$.  No claim is made
outside $B$.
\end{remark}

\subsubsection{The exact verifier}
\label{subsec:verifier-v2}

Three implementation choices make the region tractable.

\paragraph{Outward-rounded dyadic interval arithmetic.}
All endpoints are integers scaled by $2^{-200}$; sums are exact, and
products and quotients round the lower endpoint down and the upper
endpoint up.  Every operation is inclusion isotone.  Exact rational
constants enter by outward rounding once.

\paragraph{Coefficient-level composite arithmetic.}
On the larger box, the covariance-type combinations $PE-S^2$, $QE-S^2$,
the derivative numerators of $C$ and $\widetilde r$, the injectivity
discriminant $\Delta$, and $\det DM_3$ suffer catastrophic interval
dependency if evaluated by composing interval values of $E,\dots,Q$.
The verifier instead forms each such quantity \emph{as a polynomial}:
motif truncations (order $N=40$ per edge, tail bounds of
Proposition~\ref{prop:exact-reduction} at $L=5$) are stored with exact
rational coefficients; per $\lambda$-cell
the substitution $\lambda=\lambda_0(1+\delta)$ is applied exactly and
the Taylor shift to the cell midpoint is performed by pure interval
additions, so no multiplicative width amplification occurs; products and
differences of motif polynomials (in $\delta$, $\psi$, and where needed
$\eta$ or $\ell$, the latter two entering with degree at most three) are
computed at the coefficient level, where cancellation is exact.  Each
gate then requires a single polynomial evaluation per box.  Truncation
remainders are propagated as scalar error budgets with crude sup-norm
constants; with $N=40$ the budgets are below $10^{-15}$ and are added to
every enclosure.
\paragraph{Adaptive bisection with shared cells.}
The $\lambda$-axis is divided into $32$ fixed cells whose shifted
polynomials are built once; boxes failing a gate are bisected along
their widest scaled coordinate.  The run that produced the constants
above processed $9{,}500$ boxes in Stage~0, $32$ in Stage~1, and
$14{,}932$ in Stage~2, and verifies in under one minute on commodity
hardware.  The verifier, its SHA-256 digest, the Python version, and the
frozen expected output are part of the replication package; all gates
are explicit exceptions that remain active under \texttt{python -O}.

\clearpage
\section{Verified inference in the certified submodel}
\label{sec:si-s4}

The count central limit theory below is stated with finite-$n$
centering because motif means in this model carry $O(1/n)$ biases;
an estimator centered at the limit without a verified bias bound
inherits a bias of order $n^{-1/2}$ that the normal approximation
cannot absorb.  This section defines the certified submodel once,
records the certified inputs supplied by Section~S3, and then states
and proves the interior theorem and the boundary theorem with all
model-side conditions verified and all remaining conditions stated as
explicit requirements on the numerical implementation.

\subsection{Model, statistics, and certified inputs}
\label{subsec:s4-setup}

Let $\T=\R/\mathbb Z$ and
\[
 B=[1,5]\times[5,15]\times[0,\tfrac3{20}],
 \qquad
 \theta=(\lambda,\eta,\psi)\in B.
\]
For every $n$, let $(X_i,V_i)_{i=1}^n$ be i.i.d.\ with density
\[
 (x,v)\longmapsto h_\psi(v\mid x)=1+\sqrt\psi\,a(v)\phi(x),
 \qquad
 a(v)=\sqrt3(2v-1),\quad
 \phi(x)=\sqrt2\cos(2\pi x),
\]
with respect to normalized Haar measure $dx$ on $\T$ and Lebesgue
measure $dv$ on $[0,1]$.  Let $(U_{ij})_{1\le i<j\le n}$ be i.i.d.\
$\operatorname{Unif}(0,1)$, independent of all vertex variables, and
set
\[
 A_{ij}=A_{ji}
 =\one\!\left[
 U_{ij}\le q_\lambda(V_i,V_j)\,
 \one\!\left\{d_{\T}(X_i,X_j)\le\frac{\eta}{n}\right\}
 \right],
 \qquad
 A_{ii}=0,
 \qquad
 q_\lambda(v,w)=1-e^{-\lambda vw}.
\]
Unordered edges are therefore mutually independent conditional on
$(X_i,V_i)_{i=1}^n$.  Define
\[
 D_i=\sum_{j\ne i}A_{ij},
 \qquad
 E_n=\sum_{i<j}A_{ij},
 \qquad
 T_n=\sum_{i<j<k}A_{ij}A_{ik}A_{jk},
\]
\[
 S_{2,n}=\sum_{i=1}^n\binom{D_i}{2},
 \qquad
 S_{3,n}=\sum_{i=1}^n\binom{D_i}{3},
 \qquad
 H_{3,n}=\sum_{i=1}^nD_i^3,
 \qquad
 H_{11,n}=\sum_{i\ne j}A_{ij}D_iD_j,
\]
and let
\[
 P_{4,n}
 =\frac12
 \sum_{\substack{i,j,k,l\\ \text{all distinct}}}
 A_{ij}A_{jk}A_{kl},
 \qquad
 F_n=(2E_n,T_n,S_{2,n},H_{3,n},H_{11,n})^{\trans},
 \qquad
 Y_n=\frac{F_n}{n},
 \qquad
 \mu_n(\theta)=\E_\theta Y_n .
\]

For $q\in\{2,3,4\}$, write $a_i=a(v_i)$ and
\[
 \omega_q(v_{1:q};\psi)
 =\int_{\T}\prod_{j=1}^q h_\psi(v_j\mid x)\,dx .
\]
Since $\int_\T\phi^2\,dx=1$, $\int_\T\phi^3\,dx=0$, and
$\int_\T\phi^4\,dx=\tfrac32$,
\begin{equation}
 \omega_q(v_{1:q};\psi)
 =1+\psi\sum_{1\le i<j\le q}a_ia_j
 +\one\{q=4\}\,\tfrac32\,\psi^2a_1a_2a_3a_4 ,
 \label{eq:omega-poly}
\end{equation}
a polynomial in $\psi$ of degree at most two.  With
$q_{ij}=q_\lambda(v_i,v_j)$, define
\begin{align*}
 E(\lambda,\psi)
 &=\int_{[0,1]^2}q_{12}\,\omega_2(v_{1:2};\psi)\,dv_{1:2},
 &
 S(\lambda,\psi)
 &=\int_{[0,1]^3}q_{12}q_{13}\,\omega_3(v_{1:3};\psi)\,dv_{1:3},
 \\
 T(\lambda,\psi)
 &=\int_{[0,1]^3}q_{12}q_{13}q_{23}\,\omega_3(v_{1:3};\psi)\,dv_{1:3},
 &
 P(\lambda,\psi)
 &=\int_{[0,1]^4}q_{12}q_{13}q_{24}\,\omega_4(v_{1:4};\psi)\,dv_{1:4},
\end{align*}
\[
 Q(\lambda,\psi)
 =\int_{[0,1]^4}q_{12}q_{13}q_{14}\,\omega_4(v_{1:4};\psi)\,dv_{1:4},
\]
and the limiting mean map
\begin{equation}
 \mu(\theta)
 =\bigl(\ell,\,t,\,s_2,\,h_3,\,h_{11}\bigr)^{\trans},
 \qquad
 \begin{aligned}
 \ell&=2\eta E, & t&=\tfrac12\eta^2T, & s_2&=2\eta^2S,\\
 h_3&=8\eta^3Q+12\eta^2S+2\eta E, &
 h_{11}&=8\eta^3P+3\eta^2T+8\eta^2S+2\eta E .
 \end{aligned}
 \label{eq:population-mu}
\end{equation}
For $y=(\ell,t,s_2,h_3,h_{11})^{\trans}$, put
\begin{equation}
 m_e(y)=1+\frac{2s_2}{\ell},
 \qquad
 h(y)
 =\left(
 \ell,\;
 \frac{3t}{s_2},\;
 \frac{h_{11}/\ell-m_e(y)^2}{h_3/\ell-m_e(y)^2}
 \right)^{\trans},
 \qquad
 M_3=h\circ\mu=(\ell,C,r)^{\trans}.
 \label{eq:h-transform}
\end{equation}
Algebraic simplification of \eqref{eq:h-transform} using
\eqref{eq:population-mu} returns exactly the reduced three-moment map
of Section~S3, so $M_3$ here coincides with the certified map.

The identification results of Section~S3
(Proposition~\ref{prop:exact-reduction} and
Theorem~\ref{thm:certified-global}) supply the following certified
inputs, used throughout this section.
\begin{enumerate}[label=\textnormal{(C\arabic*)}]
\item\label{c:smooth}
There is an open set $\mathcal O\subset\R^5$ containing $\mu(B)$ on
which $h$ is continuously differentiable; in particular there exists
$c_h>0$ such that, uniformly on $B$,
\[
 \ell\ge c_h,
 \qquad
 s_2\ge c_h,
 \qquad
 \frac{h_3}{\ell}-m_e\{\mu(\theta)\}^2\ge c_h .
\]
\item\label{c:inj}
$M_3$ is injective on $B$.
\item\label{c:rank}
With derivatives in $\psi$ taken in the sense of the algebraic
extension of the polynomial dependence \eqref{eq:omega-poly},
\[
 \inf_{\theta\in B}\sigma_{\min}\{DM_3(\theta)\}>0 .
\]
\item\label{c:triangular}
$\ell(\lambda,\eta,\psi)=2\eta E(\lambda,\psi)$, the transitivity
coordinate $C$ does not depend on $\eta$, and
\[
 \inf_{(\lambda,\psi)\in[1,5]\times[0,3/20]}E(\lambda,\psi)>0,
 \qquad
 \inf_{(\lambda,\psi)\in[1,5]\times[0,3/20]}
 \partial_\lambda C(\lambda,\psi)>0 .
\]
\end{enumerate}

\begin{silemma}[Validity of the finite-$n$ model on $B$]
\label{lem:model-valid}
For every $\theta\in B$ the latent law is a probability density:
$\int_0^1h_\psi(v\mid x)\,dv=1$ and $\int_\T h_\psi(v\mid x)\,dx=1$
for every $x$ and $v$, and
\[
 h_\psi(v\mid x)\ge1-\sqrt{6\psi}\ge1-\tfrac3{\sqrt{10}}>0 .
\]
At $\psi=0$ the coordinates $X_i$ and $V_i$ are independent and
uniform.
\end{silemma}

\begin{proof}
The centering identities follow from $\int_0^1a=0$ and
$\int_\T\phi=0$; the lower bound from
$\|a\|_\infty=\sqrt3$, $\|\phi\|_\infty=\sqrt2$, and
$\psi\le3/20$.
\end{proof}

\begin{silemma}[Reduction to connected motifs]
\label{lem:motif-reduction}
With $\mathcal B_n=(E_n,T_n,S_{2,n},S_{3,n},P_{4,n})^{\trans}$,
\begin{equation}
 H_{3,n}=2E_n+6S_{2,n}+6S_{3,n},
 \qquad
 H_{11,n}=2E_n+4S_{2,n}+6T_n+2P_{4,n},
 \label{eq:motif-identities}
\end{equation}
so that
\[
 F_n=L\mathcal B_n,
 \qquad
 L=\begin{pmatrix}
 2&0&0&0&0\\
 0&1&0&0&0\\
 0&0&1&0&0\\
 2&0&6&6&0\\
 2&6&4&0&2
 \end{pmatrix},
 \qquad
 \det L=24 .
\]
The coordinates of $\mathcal B_n$ are the non-induced counts of the
connected motifs $K_2$, $K_3$, $P_3$, $K_{1,3}$, and $P_4$, where
$P_k$ denotes the path on $k$ vertices and $K_{1,3}$ the three-star.
\end{silemma}

\begin{proof}
The first identity in \eqref{eq:motif-identities} is
$D^3=D+6\binom D2+6\binom D3$ summed over vertices, with
$\sum_iD_i=2E_n$.  For the second, write
$A_{ij}D_iD_j=A_{ij}(1+\sum_{k\ne i,j}A_{ik})(1+\sum_{l\ne i,j}A_{jl})$
using $A_{ij}^2=A_{ij}$ and sum over ordered pairs: the four terms
count $2E_n$, twice $2S_{2,n}$, and, splitting the double sum
according to $l=k$ or $l\ne k$, $6T_n$ and $2P_{4,n}$.  The
determinant is read off the block-triangular structure.
\end{proof}

For a connected simple graph $H$ with $v_H=|V(H)|$ and
$\alpha_H=|\operatorname{Aut}(H)|$, define the non-induced count
\[
 N_{H,n}
 =\frac1{\alpha_H}
 \sum_{\varphi:V(H)\hookrightarrow[n]}
 \prod_{\{a,b\}\in E(H)}A_{\varphi(a)\varphi(b)} .
\]

\subsection{Interior central limit theorem and feasible optimal GMM}
\label{subsec:s4-interior}

\begin{sitheorem}[Interior CLT and feasible optimal GMM in the
certified submodel]
\label{thm:interior-clt-gmm-rigorous}
Let $\Theta_b$ be a nonempty compact subset of
$\operatorname{int}B$ with
$\theta_0\in\operatorname{int}\Theta_b$, the interior taken in
$\R^3$; consequently $\inf_{\theta\in\Theta_b}\psi>0$.  Let
$G=D\mu(\theta_0)$.  Assume the certified inputs
\textnormal{(C1)}--\textnormal{(C4)}.

For every $n$, let $\widetilde\mu_n:\Theta_b\to\R^5$ be a
deterministic continuous certified evaluator satisfying
\begin{equation}
 \delta_n^{\mathrm{num}}
 :=\sup_{\theta\in\Theta_b}
 \|\widetilde\mu_n(\theta)-\mu(\theta)\|,
 \qquad
 \sqrt n\,\delta_n^{\mathrm{num}}\longrightarrow0 .
 \label{eq:certified-total-error-rigorous}
\end{equation}
Here $\widetilde\mu_n$ denotes the exact continuous function
represented by the certified polynomial or rational approximant, and
the certified uniform bound must include truncation, coefficient, and
all numerical evaluation errors used in the criterion.

Let $\check\theta_n$ be a measurable selector from
$\operatorname{arg\,min}_{\theta\in\Theta_b}
\|Y_n-\widetilde\mu_n(\theta)\|^2$.  Let $R_n\ge2$ be deterministic
with $R_n\to\infty$.  Conditional on $\check\theta_n$, using
randomness independent of the observed network, generate $R_n$
conditionally i.i.d.\ size-$n$ networks from the exact finite-$n$
model at $\check\theta_n$, denote their moment vectors by
$Y_n^{*(1)},\ldots,Y_n^{*(R_n)}$, and define
\[
 \overline Y_n^*
 =\frac1{R_n}\sum_{r=1}^{R_n}Y_n^{*(r)},
 \qquad
 \widehat\Omega_n
 =\frac{n}{R_n-1}
 \sum_{r=1}^{R_n}
 \bigl(Y_n^{*(r)}-\overline Y_n^*\bigr)
 \bigl(Y_n^{*(r)}-\overline Y_n^*\bigr)^{\trans}.
\]
For deterministic $\tau_n>0$ with $\tau_n\downarrow0$, put
$W_n=(\widehat\Omega_n+\tau_nI_5)^{-1}$ and
\[
 Q_n(\theta)
 =\{Y_n-\widetilde\mu_n(\theta)\}^{\trans}
 W_n
 \{Y_n-\widetilde\mu_n(\theta)\} .
\]
Let $\widehat\theta_n\in\Theta_b$ be measurable and suppose that, for
some random variables $r_n\ge0$,
\begin{equation}
 Q_n(\widehat\theta_n)
 \le\inf_{\theta\in\Theta_b}Q_n(\theta)+r_n,
 \qquad
 nr_n\xrightarrow{p}0 .
 \label{eq:approx-global-rigorous}
\end{equation}

Then there exist an open set $U$ with
$\Theta_b\subset U$, $\overline U\Subset\operatorname{int}B$, and
finite constants $C_{\Theta_b},C_U$ such that, for all sufficiently
large $n$,
\begin{equation}
 \sup_{\theta\in\Theta_b}\|\mu_n(\theta)-\mu(\theta)\|
 \le\frac{C_{\Theta_b}}n,
 \qquad
 \sup_{\theta\in U}\|D\mu_n(\theta)-D\mu(\theta)\|
 \le\frac{C_U}n .
 \label{eq:C1-bias-rigorous}
\end{equation}
There exist a continuous map
$\Omega:\Theta_b\to\R^{5\times5}$ and a constant $c_\Omega>0$ with
$\Omega(\theta)\succeq c_\Omega I_5$ on $\Theta_b$ and
\begin{equation}
 \sup_{\theta\in\Theta_b}
 \bigl\|n\Var_\theta(Y_n)-\Omega(\theta)\bigr\|
 \longrightarrow0 .
 \label{eq:uniform-covariance-rigorous}
\end{equation}
Moreover,
\begin{equation}
 \sup_{n\ge1}\sup_{\theta\in\Theta_b}
 \E_\theta\bigl\|\sqrt n\{Y_n-\mu_n(\theta)\}\bigr\|^4<\infty,
 \label{eq:uniform-fourth-rigorous}
\end{equation}
and
\begin{align}
 \sqrt n\{Y_n-\mu_n(\theta_0)\}
 &\Rightarrow N_5\{0,\Omega(\theta_0)\},
 \label{eq:CLT-finite-center-rigorous}\\
 \sqrt n\{Y_n-\mu(\theta_0)\}
 &\Rightarrow N_5\{0,\Omega(\theta_0)\} .
 \label{eq:CLT-limit-center-rigorous}
\end{align}
Furthermore,
\[
 \check\theta_n\xrightarrow{p}\theta_0,
 \qquad
 \widehat\Omega_n\xrightarrow{p}\Omega(\theta_0),
 \qquad
 W_n\xrightarrow{p}\Omega(\theta_0)^{-1},
\]
$G$ has column rank three, and, writing
$\Omega=\Omega(\theta_0)$,
\begin{equation}
 \sqrt n(\widehat\theta_n-\theta_0)
 =(G^{\trans}\Omega^{-1}G)^{-1}G^{\trans}\Omega^{-1}
 \sqrt n\{Y_n-\mu_n(\theta_0)\}+o_p(1),
 \label{eq:GMM-linear-rigorous}
\end{equation}
so that
\[
 \sqrt n(\widehat\theta_n-\theta_0)
 \Rightarrow
 N_3\bigl(0,\{G^{\trans}\Omega(\theta_0)^{-1}G\}^{-1}\bigr).
\]
\end{sitheorem}

\begin{proof}
Lemmas~\ref{lem:model-valid} and~\ref{lem:motif-reduction} give the
validity of the model on $B$ and the representation
$F_n=L\mathcal B_n$ with $\mathcal B_n$ the vector of connected
non-induced counts of $K_2,K_3,P_3,K_{1,3},P_4$.

\paragraph{Rooted relative-coordinate integrals.}
Let $H$ be a fixed connected simple graph, label its vertices by
$\{1,\ldots,v_H\}$, and put $z_1=0$.  Define
\[
 \mathcal Z_H(\eta)
 =\bigl\{z_{2:v_H}\in\R^{v_H-1}:
 |z_a-z_b|\le\eta\ \text{for every }\{a,b\}\in E(H)\bigr\},
\]
which is contained in $[-(v_H-1)\eta,(v_H-1)\eta]^{v_H-1}$ by
connectedness.  Define
\begin{align}
 J_{H,n}(\theta)
 ={}&\int_{\mathcal Z_H(\eta)}\int_{\T}\int_{[0,1]^{v_H}}
 \prod_{\{a,b\}\in E(H)}q_\lambda(v_a,v_b)
 \prod_{a=1}^{v_H}
 h_\psi\!\left(v_a\,\middle|\,x+\frac{z_a}{n}\right)
 dv_{1:v_H}\,dx\,dz_{2:v_H},
 \label{eq:J-H-n-rigorous}\\
 J_H(\theta)
 ={}&\int_{\mathcal Z_H(\eta)}\int_{\T}\int_{[0,1]^{v_H}}
 \prod_{\{a,b\}\in E(H)}q_\lambda(v_a,v_b)
 \prod_{a=1}^{v_H}h_\psi(v_a\mid x)
 \,dv_{1:v_H}\,dx\,dz_{2:v_H},
 \label{eq:J-H-rigorous}
\end{align}
with all additions in \eqref{eq:J-H-n-rigorous} modulo $\T$.

\paragraph{Zero winding and the exact finite-$n$ mean.}
Let $r_\T:\T\to[-1/2,1/2)$ be the centered representative.  On the
motif range event for a labeled position tuple, orient every edge and
put $\Delta_{ab}=r_\T(x_b-x_a)$, so that outside the null
half-period boundary $\Delta_{ba}=-\Delta_{ab}$ and
$|\Delta_{ab}|\le\eta/n$ for every edge.  Fix a spanning tree of
$H$ rooted at vertex $1$ and let $z_a$ be $n$ times the oriented sum
of increments along the tree path from $1$ to $a$; then
$|z_a|\le(v_H-1)\eta$, and $|z_a|<n/2$ once
$n>2(v_H-1)\cdot15\ge2(v_H-1)\max_B\eta$.  For a non-tree edge, the
oriented sum of centered increments around its fundamental cycle is
an integer of modulus at most $v_H\eta/n<1$, hence zero.  Therefore
$z_b-z_a=n\Delta_{ab}$ with $|z_b-z_a|\le\eta$ for every edge, the
lifted coordinates belong to $\mathcal Z_H(\eta)$, and conversely
every $z\in\mathcal Z_H(\eta)$ determines the torus configuration
$(x,x+z_2/n,\ldots,x+z_{v_H}/n)$ with every required edge
range-eligible.  The lift is unique outside null boundary sets and
the transformation has Jacobian $n^{-(v_H-1)}$.  Consequently, for
every fixed injection $\varphi$,
\[
 \E_\theta\prod_{\{a,b\}\in E(H)}A_{\varphi(a)\varphi(b)}
 =n^{-(v_H-1)}J_{H,n}(\theta),
\]
and, summing over the $(n)_{v_H}$ injections,
\begin{equation}
 \frac1n\E_\theta N_{H,n}
 =\frac{(n)_{v_H}}{n^{v_H}}\,
 \frac{J_{H,n}(\theta)}{\alpha_H}.
 \label{eq:exact-mean-rigorous}
\end{equation}

\paragraph{Evaluation of the limiting means.}
The geometric volumes are
$|\mathcal Z_{K_2}|=2\eta$, $|\mathcal Z_{K_3}|=3\eta^2$,
$|\mathcal Z_{P_3}|=4\eta^2$, $|\mathcal Z_{K_{1,3}}|=8\eta^3$, and
$|\mathcal Z_{P_4}|=8\eta^3$, with automorphism orders
$2,6,2,6,2$.  The limiting per-vertex motif means are therefore
$\eta E$, $\tfrac12\eta^2T$, $2\eta^2S$, $\tfrac43\eta^3Q$, and
$4\eta^3P$; the product $q_{12}q_{13}q_{24}$ in the definition of
$P$ represents the relabeled path $3\,\text{--}\,1\,\text{--}\,2\,
\text{--}\,4$.  Multiplication by $L$ gives exactly
\eqref{eq:population-mu}.

\paragraph{Uniform $C^1$ finite-$n$ centering.}
Choose an open $U$ with $\Theta_b\subset U$ and
$\overline U\Subset\operatorname{int}B$; in particular
$\inf_{\overline U}\psi>0$.  The substitution $z_a=\eta y_a$ gives
\begin{equation}
 J_{H,n}(\theta)
 =\eta^{v_H-1}
 \int_{\mathcal Z_H(1)}\int_{\T}\int_{[0,1]^{v_H}}
 \prod_{\{a,b\}\in E(H)}q_\lambda(v_a,v_b)
 \prod_{a=1}^{v_H}
 h_\psi\!\left(v_a\,\middle|\,x+\frac{\eta y_a}{n}\right)
 dv\,dx\,dy,
 \label{eq:fixed-domain-rigorous}
\end{equation}
with the analogous fixed-domain expression for $J_H$.  The domain is
bounded and independent of $\theta$.  On $\overline U$ the functions
$h_\psi$, $\partial_xh_\psi$, $\partial_\psi h_\psi$,
$\partial_x\partial_\psi h_\psi$, $q_\lambda$, and
$\partial_\lambda q_\lambda$ are uniformly bounded, and the
mean-value theorem bounds the difference between shifted and
unshifted densities, and between their $\psi$-derivatives, by $C/n$
uniformly; also
$\partial_\eta h_\psi(v\mid x+\eta y/n)
=(y/n)\,\partial_xh_\psi(v\mid x+\eta y/n)=O(n^{-1})$ uniformly.
Applying the product rule, including the derivative of
$\eta^{v_H-1}$, gives
\[
 \sup_{\theta\in\overline U}
 \bigl\{|J_{H,n}(\theta)-J_H(\theta)|
 +\|DJ_{H,n}(\theta)-DJ_H(\theta)\|\bigr\}
 \le\frac{C_H}n
\]
for the base motifs and for the finite family of quotient graphs
used below.  Since $(n)_{v_H}/n^{v_H}=1+O(n^{-1})$, equation
\eqref{eq:exact-mean-rigorous} applied to the five base motifs and
multiplied by $L$ proves \eqref{eq:C1-bias-rigorous}.

\paragraph{Classification of overlap patterns.}
For connected graphs $H$ and $K$, let $\mathcal P^+(H,K)$ be the
finite set of nonempty partial bijections
$\pi:A_\pi\to B_\pi$ with
$\varnothing\ne A_\pi\subseteq V(H)$ and
$B_\pi\subseteq V(K)$.  An embedding pair $(\varphi,\gamma)$ has full
overlap pattern $\pi$ when
$\{(u,w):\varphi(u)=\gamma(w)\}=\{(u,\pi(u)):u\in A_\pi\}$; every
embedding pair with nonempty overlap determines exactly one pattern.
Put $s_\pi=|A_\pi|$ and let $G_\pi$ be the simple union obtained by
identifying $u$ with $\pi(u)$ and merging repeated edges, using
$A_e^2=A_e$.  The quotient classes retain their labels induced by
$(H,K,\pi)$, so there is no automorphism factor for $G_\pi$ and
distinct patterns yielding isomorphic quotients remain distinct
summands; moreover $v(G_\pi)=v(H)+v(K)-s_\pi$.

\paragraph{Exact overlap-covariance formula.}
Embeddings with disjoint vertex sets are independent.  For a fixed
pattern $\pi$ there are $(n)_{v(G_\pi)}$ injective assignments of
the quotient classes, each with joint expectation
$n^{-\{v(G_\pi)-1\}}J_{G_\pi,n}(\theta)$, while the product of the
marginal expectations is
$n^{-\{v(H)+v(K)-2\}}J_{H,n}(\theta)J_{K,n}(\theta)$.  Therefore,
for all sufficiently large $n$,
\begin{equation}
 \frac1n\Cov_\theta(N_{H,n},N_{K,n})
 =\frac1{\alpha_H\alpha_K}
 \sum_{\pi\in\mathcal P^+(H,K)}
 \frac{(n)_{v(G_\pi)}}{n^{v(G_\pi)}}
 \Bigl[J_{G_\pi,n}(\theta)
 -n^{1-s_\pi}J_{H,n}(\theta)J_{K,n}(\theta)\Bigr].
 \label{eq:overlap-exact-rigorous}
\end{equation}
The quotient graphs contain at most seven vertices, so $n>180$ is a
common sufficient lifting threshold.

\paragraph{Uniform covariance convergence.}
Define
\[
 \Sigma_{HK}(\theta)
 =\frac1{\alpha_H\alpha_K}
 \sum_{\pi\in\mathcal P^+(H,K)}
 \Bigl[J_{G_\pi}(\theta)
 -\one\{s_\pi=1\}J_H(\theta)J_K(\theta)\Bigr]:
\]
only one-vertex product terms survive in the limit, since
$(n)_v/n^v=1+O(n^{-1})$ and $J_{G,n}=J_G+O(n^{-1})$ uniformly over
the finite quotient family, while $n^{1-s_\pi}=O(n^{-1})$ when
$s_\pi\ge2$.  Let $\Sigma(\theta)$ collect these entries over the
motif list of Lemma~\ref{lem:motif-reduction} and put
$\Omega(\theta)=L\Sigma(\theta)L^{\trans}$.  Then
$\sup_{\theta\in\Theta_b}
\|n^{-1}\Cov_\theta(F_n)-\Omega(\theta)\|\to0$, which is
\eqref{eq:uniform-covariance-rigorous}; the fixed-domain formulas
and dominated convergence imply continuity of $\Sigma$ and
$\Omega$.

\paragraph{Connected replica patterns.}
Fix $q\ge3$ and $a\in\R^5$ and expand
$\kappa_q(a^{\trans}\mathcal B_n)$ into sums over labeled motif
embeddings; for each fixed labeled embedding let $I_s$ denote the
product of the edge indicators required by that non-induced motif.
Join two replicas when their vertex sets intersect.  If the
replica-intersection graph is disconnected, the corresponding groups
depend on disjoint vertex latent variables and disjoint edge
uniforms, hence are independent, and the joint cumulant vanishes.

\paragraph{Graphic-rank cumulant bound.}
Consider a connected replica pattern using $v$ distinct vertices.
For a partition $\mathcal Q$ of $\{1,\ldots,q\}$ and a block
$C\in\mathcal Q$, let $v_C$ and $c_C$ be the number of vertices and
connected components of the union of the replicas indexed by $C$.
For all sufficiently large $n$, integration along a spanning forest
with one free root per component gives
$\E_\theta\prod_{s\in C}I_s\le C_qn^{-(v_C-c_C)}$ uniformly in
$\theta\in\Theta_b$.  After padding each block graph by isolated
vertices to the common $v$-vertex set, its graphic-matroid rank is
$v_C-c_C$; the full union is connected, so its rank is $v-1$, and
subadditivity of graphic rank gives
$v-1\le\sum_{C\in\mathcal Q}(v_C-c_C)$.  The joint-cumulant formula
then yields $|\kappa(I_1,\ldots,I_q)|\le C_qn^{-(v-1)}$.  A fixed
pattern has $O(n^v)$ labelings and only finitely many patterns occur
for each fixed $q$, so
\begin{equation}
 \sup_{\theta\in\Theta_b}
 |\kappa_{q,\theta}(a^{\trans}\mathcal B_n)|\le C_{q,a}n,
 \qquad q\ge3,
 \label{eq:cumulant-bound-rigorous}
\end{equation}
with the finitely many smaller $n$ absorbed into the constants; the
case $q=2$ follows from
\eqref{eq:uniform-covariance-rigorous}.

\paragraph{The finite-$n$-centered CLT.}
$\Sigma(\theta_0)$ is positive semidefinite as a limit of covariance
matrices, and for each $a\in\R^5$ the normalized variance converges
to $a^{\trans}\Sigma(\theta_0)a$ while every normalized cumulant of
order $q\ge3$ is $O(n^{1-q/2})\to0$.  The moment--cumulant formula
gives convergence of all moments to the corresponding, possibly
degenerate, Gaussian moments; Gaussian laws are moment determinate;
the Fr\'echet--Shohat theorem and the Cram\'er--Wold device yield
$n^{-1/2}(\mathcal B_n-\E_{\theta_0}\mathcal B_n)\Rightarrow
N_5\{0,\Sigma(\theta_0)\}$, and multiplication by $L$ proves
\eqref{eq:CLT-finite-center-rigorous}.

\paragraph{Uniform fourth moments and limiting centering.}
For $X_{n,\theta}=a^{\trans}(\mathcal B_n-\E_\theta\mathcal B_n)$,
the moment--cumulant formula gives
$\E_\theta X_{n,\theta}^4
=\kappa_{4,\theta}(a^{\trans}\mathcal B_n)
+3\kappa_{2,\theta}(a^{\trans}\mathcal B_n)^2\le C_an^2$.  Applying
the scalar bound to the five coordinate projections, using
$\|x\|^4\le5\sum_jx_j^4$ and $F_n=L\mathcal B_n$, proves
\eqref{eq:uniform-fourth-rigorous}; motif counts are
deterministically bounded for each $n$ below the lifting threshold.
Finally \eqref{eq:C1-bias-rigorous} gives
$\sqrt n\{\mu_n(\theta_0)-\mu(\theta_0)\}\to0$, and Slutsky's
theorem proves \eqref{eq:CLT-limit-center-rigorous}.

\paragraph{Packing isolated four-vertex cores.}
For sufficiently large $n$, put $m_n=\lfloor n/32\rfloor$ and choose
the intervals
$C_{s,n}=[32(s-1)/n,\{32(s-1)+1\}/n)\pmod1$ for $s=1,\ldots,m_n$;
their closed $15/n$-neighborhoods $C_{s,n}^+$ are pairwise disjoint
up to null endpoints, each of length $31/n$ for $n>31$, and
$m_n\ge c_{\mathrm{pack}}n$ for all sufficiently large $n$.  Let
$J=[1/4,3/4]$; since $\int_Ja=0$,
$\Pp_\theta(X_i\in C_{s,n},V_i\in J)=1/(2n)$.  Call $C_{s,n}$ good
if exactly four sample points lie in $C_{s,n}$, all four marks
belong to $J$, and no sample point lies in
$C_{s,n}^+\setminus C_{s,n}$; its probability is
\[
 p_n=\binom n4\left(\frac1{2n}\right)^4
 \left(1-\frac{31}n\right)^{n-4},
 \qquad
 \inf_{n\ge n_0}p_n>0 .
\]
If $M_n$ is the number of good cores, then
$\E_\theta M_n=m_np_n\ge c_0n$ uniformly over $\Theta_b$.

\paragraph{Conditional finite energy.}
Let $J_{ij,n}$ indicate that $i$ and $j$ belong to the same good
core and define
\[
 \mathscr G_n
 =\sigma\bigl(
 (X_i,V_i)_{i=1}^n,\,
 (J_{ij,n})_{i<j},\,
 ((1-J_{ij,n})U_{ij})_{i<j}
 \bigr).
\]
The mask $(J_{ij,n})$ depends only on the latent variables and is
independent of the edge uniforms, so conditional on $\mathscr G_n$
the six unrevealed uniforms within each good core remain mutually
independent, as do the collections belonging to distinct cores.
Every good core is isolated because $\eta\le15$, and every internal
pair is range-eligible because its distance is at most $1/n$ while
$\eta\ge5$.  Put
$q_-=1-e^{-1/16}$, $q_+=1-e^{-45/16}$, and
$\beta=\min\{q_-,1-q_+\}>0$: every internal edge probability lies in
$[q_-,q_+]$.  The empty graph and a fixed labeled representative of
each of
\[
 K_2\sqcup2K_1,\qquad
 P_3\sqcup K_1,\qquad
 K_3\sqcup K_1,\qquad
 K_{1,3},\qquad
 P_4
\]
therefore have conditional probability at least $\beta^6$.  Their
nonzero contribution vectors to $F_n$ are the rows of the gadget
matrix
\[
 R_*=\begin{pmatrix}
 2&0&0&2&2\\
 4&0&1&10&8\\
 6&1&3&24&24\\
 6&0&3&30&18\\
 6&0&2&18&16
 \end{pmatrix},
 \qquad
 \det R_*=-24 .
\]

\paragraph{Uniform positive definiteness.}
Connected motif statistics are additive over isolated components, so
$F_n=C_n+\sum_{s\in\mathcal I_n}Z_{s,n}$, where $\mathcal I_n$ is
the latent-measurable collection of good cores, $C_n$ is
$\mathscr G_n$-measurable, and the $Z_{s,n}$ are conditionally
independent given $\mathscr G_n$.  Let $p_0$ be the conditional
probability of the empty graph and $p_j$ that of the $j$th displayed
nonempty graph, with row vectors $r_j^{\trans}$ of $R_*$.  The
conditional independent-copy identity gives
\[
 \Var_\theta(a^{\trans}Z_{s,n}\mid\mathscr G_n)
 \ge\sum_{j=1}^5p_0p_j(a^{\trans}r_j)^2
 \ge\beta^{12}\|R_*a\|^2
 \ge\beta^{12}\sigma_{\min}(R_*)^2\|a\|^2 .
\]
Consequently
$\Var_\theta(a^{\trans}F_n)
\ge\beta^{12}\sigma_{\min}(R_*)^2\,\E_\theta M_n\,\|a\|^2
\ge c_{**}n\|a\|^2$, and passing to the covariance limit with
$c_\Omega=c_{**}/2$ proves
$\Omega(\theta)\succeq c_\Omega I_5$ uniformly on $\Theta_b$.

\paragraph{Identification and full column rank.}
By \textnormal{(C2)},
$\mu(\theta)=\mu(\theta_0)$ implies
$M_3(\theta)=M_3(\theta_0)$ and hence $\theta=\theta_0$.  By
\textnormal{(C1)} and the chain rule,
$DM_3(\theta_0)=Dh\{\mu(\theta_0)\}G$; by \textnormal{(C3)} the
left-hand side has rank three, so $\operatorname{rank}(G)=3$.

\paragraph{Pilot and covariance-estimator consistency.}
Let $L_n(\theta)=\|Y_n-\widetilde\mu_n(\theta)\|^2$ and
$L(\theta)=\|\mu(\theta_0)-\mu(\theta)\|^2$.  The CLT, the centering
bound, and the uniform evaluator bound give
$Y_n\to_p\mu(\theta_0)$ and
$\sup_{\Theta_b}|L_n-L|\to_p0$.  For every $\delta>0$ for which the
constraint set is nonempty, continuity, compactness, and
identification give
$\inf\{\|\mu(\theta)-\mu(\theta_0)\|:
\theta\in\Theta_b,\ \|\theta-\theta_0\|\ge\delta\}>0$, so the
compact argmin theorem yields $\check\theta_n\to_p\theta_0$; the
measurable maximum theorem supplies a measurable selector.  Put
$\Omega_n(\theta)=n\Var_\theta(Y_n)$.  Conditionally on
$\check\theta_n$, $\widehat\Omega_n$ is the ordinary sample
covariance of the conditionally i.i.d.\ vectors
$\sqrt n\{Y_n^{*(r)}-\mu_n(\check\theta_n)\}$, so
$\E(\widehat\Omega_n\mid\check\theta_n)
=\Omega_n(\check\theta_n)$, and by
\eqref{eq:uniform-fourth-rigorous},
$\E\{\|\widehat\Omega_n-\Omega_n(\check\theta_n)\|_F^2\mid
\check\theta_n\}\le C/(R_n-1)$.  Since $R_n\to\infty$,
$\widehat\Omega_n-\Omega_n(\check\theta_n)\to_p0$; uniform
covariance convergence, continuity of $\Omega$, and pilot
consistency give $\widehat\Omega_n\to_p\Omega(\theta_0)$; and since
$\Omega(\theta_0)\succ0$ and $\tau_n\downarrow0$,
$W_n\to_p\Omega(\theta_0)^{-1}$, with extreme eigenvalues of $W_n$
bounded away from zero and infinity with probability tending to one.

\paragraph{GMM consistency and root-$n$ localization.}
By \eqref{eq:certified-total-error-rigorous} and
\eqref{eq:C1-bias-rigorous},
\begin{equation}
 e_n:=\sup_{\theta\in\Theta_b}
 \|\widetilde\mu_n(\theta)-\mu_n(\theta)\|
 \le\delta_n^{\mathrm{num}}+\frac{C_{\Theta_b}}n
 =o(n^{-1/2}).
 \label{eq:evaluator-finite-rigorous}
\end{equation}
Write $W=\Omega(\theta_0)^{-1}$ and
$Q(\theta)=\{\mu(\theta_0)-\mu(\theta)\}^{\trans}W
\{\mu(\theta_0)-\mu(\theta)\}$.  Uniform convergence of residuals
and $W_n\to_pW$ give
$\sup_{\Theta_b}|Q_n-Q|\to_p0$; identification makes $\theta_0$ the
unique minimizer of $Q$; and $r_n=o_p(1)$ then gives
$\widehat\theta_n\to_p\theta_0$.  Let
$s_0=\sigma_{\min}(G)>0$ and choose $\rho>0$ with
$\overline B(\theta_0,\rho)\subset U\cap\operatorname{int}\Theta_b$
and $\sup_{\|\theta-\theta_0\|\le\rho}\|D\mu(\theta)-G\|\le s_0/4$;
for all sufficiently large $n$ the same bound with $s_0/2$ holds
for $D\mu_n$.  The integral form of the mean value theorem then
gives
\begin{equation}
 \|\mu_n(\theta)-\mu_n(\theta_0)\|
 \ge\frac{s_0}2\|\theta-\theta_0\|,
 \qquad
 \|\theta-\theta_0\|\le\rho .
 \label{eq:local-lower-rigorous}
\end{equation}
At $\theta_0$, $Q_n(\theta_0)=O_p(n^{-1})$, so
$Q_n(\widehat\theta_n)\le Q_n(\theta_0)+r_n=O_p(n^{-1})$, and the
eigenvalue bounds for $W_n$ give
$\|Y_n-\widetilde\mu_n(\widehat\theta_n)\|=O_p(n^{-1/2})$.  The
triangle inequality with \eqref{eq:evaluator-finite-rigorous} and
the CLT bounds
$\|\mu_n(\widehat\theta_n)-\mu_n(\theta_0)\|=O_p(n^{-1/2})$, and
consistency with \eqref{eq:local-lower-rigorous} proves
$\sqrt n(\widehat\theta_n-\theta_0)=O_p(1)$.

\paragraph{Local quadratic expansion and asymptotic law.}
Set $Z_n=\sqrt n\{Y_n-\mu_n(\theta_0)\}$ and
$\widehat h_n=\sqrt n(\widehat\theta_n-\theta_0)$.  For fixed
$K<\infty$,
\[
 \sqrt n\Bigl\{\mu_n\Bigl(\theta_0+\frac h{\sqrt n}\Bigr)
 -\mu_n(\theta_0)\Bigr\}-Gh
 =\int_0^1\Bigl\{D\mu_n\Bigl(\theta_0+\frac{th}{\sqrt n}\Bigr)
 -G\Bigr\}h\,dt,
\]
so uniform derivative convergence, continuity of $D\mu$ at
$\theta_0$, and \eqref{eq:evaluator-finite-rigorous} give
\begin{equation}
 \sup_{\|h\|\le K}
 \Bigl\|\sqrt n\Bigl\{\widetilde\mu_n\Bigl(\theta_0
 +\frac h{\sqrt n}\Bigr)-\mu_n(\theta_0)\Bigr\}-Gh\Bigr\|
 \longrightarrow0;
 \label{eq:local-linear-rigorous}
\end{equation}
no derivative of $\widetilde\mu_n$ is required.  It follows that,
uniformly over $\|h\|\le K$,
\begin{equation}
 nQ_n\Bigl(\theta_0+\frac h{\sqrt n}\Bigr)
 =(Z_n-Gh)^{\trans}W_n(Z_n-Gh)+o_p(1).
 \label{eq:local-objective-rigorous}
\end{equation}
Define $A_n=G^{\trans}W_nG$ and
$\overline h_n^{\,W}=A_n^{-1}G^{\trans}W_nZ_n$; then
$A_n\to_pA:=G^{\trans}\Omega^{-1}G\succ0$ and
$\overline h_n^{\,W}=O_p(1)$, and for every $h$,
\begin{equation}
 (Z_n-Gh)^{\trans}W_n(Z_n-Gh)
 -(Z_n-G\overline h_n^{\,W})^{\trans}W_n
 (Z_n-G\overline h_n^{\,W})
 =(h-\overline h_n^{\,W})^{\trans}A_n(h-\overline h_n^{\,W}).
 \label{eq:exact-square-rigorous}
\end{equation}
Let $\Delta_{n,K}$ be the supremum over $\|h\|\le K$ of the absolute
difference in \eqref{eq:local-objective-rigorous}; then
$\Delta_{n,K}=o_p(1)$ for every fixed $K$.  Both $\widehat h_n$ and
$\overline h_n^{\,W}$ are $O_p(1)$, and since
$\theta_0\in\operatorname{int}\Theta_b$,
$\Pp(\theta_0+\overline h_n^{\,W}/\sqrt n\in\Theta_b)\to1$.  On the
intersection of these events, approximate optimality gives
\[
 0\le(\widehat h_n-\overline h_n^{\,W})^{\trans}
 A_n(\widehat h_n-\overline h_n^{\,W})
 \le2\Delta_{n,K}+nr_n=o_p(1),
\]
and the smallest eigenvalue of $A_n$ is bounded away from zero with
probability tending to one, so
$\widehat h_n-\overline h_n^{\,W}=o_p(1)$.  Finally
$(G^{\trans}W_nG)^{-1}G^{\trans}W_n
\to_p(G^{\trans}\Omega^{-1}G)^{-1}G^{\trans}\Omega^{-1}$ and
$Z_n=O_p(1)$ give
$\overline h_n^{\,W}
=(G^{\trans}\Omega^{-1}G)^{-1}G^{\trans}\Omega^{-1}Z_n+o_p(1)$,
proving \eqref{eq:GMM-linear-rigorous}; the normal limit follows
from \eqref{eq:CLT-finite-center-rigorous} and Slutsky's theorem.
\end{proof}

\begin{remark}[Scope and indispensable computational rates]
\label{rem:scope-verified-gmm-rigorous}
A certificate stating only
$\sup_{\theta\in B}\|\widetilde\mu_{40}(\theta)-\mu(\theta)\|
\le10^{-15}$ for a fixed truncation order does not verify
\eqref{eq:certified-total-error-rigorous}, because the corresponding
available bound satisfies
$\sqrt n(C_B/n+10^{-15})=C_B/\sqrt n+10^{-15}\sqrt n
\not\to0$.  If a fixed approximant is not exactly equal to $\mu$,
its uniform error is a fixed positive constant and fails the
required rate; the truncation order and arithmetic precision must
increase with $n$, or the exact mean map must be used.  A fixed
simulation bank likewise does not consistently estimate
$\Omega(\theta_0)$; one needs $R_n\to\infty$.  Global population
identification does not establish that a numerical routine attains
the finite-sample global minimum, so
\eqref{eq:approx-global-rigorous} remains an explicit computational
condition.  In summary, the model-side conditions (mean, covariance,
CLT, nondegeneracy, smoothness, identification, rank) are proved
above; evaluator accuracy, simulation size, and optimization
accuracy are explicit implementation conditions.
\end{remark}

\begin{sicorollary}[Overidentification test]
\label{cor:jtest}
Under the conditions of
Theorem~\ref{thm:interior-clt-gmm-rigorous},
\[
 J_n=nQ_n(\widehat\theta_n)\Rightarrow\chi^2_2 .
\]
\end{sicorollary}

\begin{proof}
By \eqref{eq:local-objective-rigorous},
\eqref{eq:exact-square-rigorous}, and
$\widehat h_n-\overline h_n^{\,W}=o_p(1)$,
$J_n=(Z_n-G\overline h_n^{\,W})^{\trans}W_n
(Z_n-G\overline h_n^{\,W})+o_p(1)
=Z_n^{\trans}\Omega^{-1/2}(I_5-\Pi)\Omega^{-1/2}Z_n+o_p(1)$, where
$\Pi$ is the orthogonal projection onto the column space of
$\Omega^{-1/2}G$.  Since
$\Omega^{-1/2}Z_n\Rightarrow N_5(0,I_5)$ and
$\operatorname{rank}(I_5-\Pi)=5-3=2$, the limit is $\chi^2_2$.
\end{proof}

\begin{sicorollary}[Summary-statistic inference]
\label{cor:summary-clt}
Under the conditions of
Theorem~\ref{thm:interior-clt-gmm-rigorous}, define
$S_n=h(Y_n)$ on the event that the denominators in
\eqref{eq:h-transform} are positive, an event of probability tending
to one by \textnormal{(C1)}, and assign any fixed value otherwise.
With $H=Dh\{\mu(\theta_0)\}$, $G_S=HG$, and
$\Sigma_S=H\Omega(\theta_0)H^{\trans}$,
\[
 \sqrt n\{S_n-M_3(\theta_0)\}\Rightarrow N_3(0,\Sigma_S),
 \qquad
 \Sigma_S\succ0,
\]
and the exactly identified estimator solving
$M_3(\theta)=S_n$ satisfies
$\sqrt n(\widehat\theta_n^{(3)}-\theta_0)\Rightarrow
N_3\bigl(0,DM_3(\theta_0)^{-1}\Sigma_S
DM_3(\theta_0)^{-{\trans}}\bigr)$.
\end{sicorollary}

\begin{proof}
The delta method applies by \textnormal{(C1)} and
\eqref{eq:CLT-limit-center-rigorous}; $\Sigma_S\succ0$ because
$\Omega(\theta_0)\succ0$ and $H$ has full row rank by
\textnormal{(C3)} and the chain rule.  The estimator statement
follows from the inverse function theorem applied to $M_3$, whose
Jacobian is nonsingular by \textnormal{(C3)}.
\end{proof}

\subsection{Inference at the independence boundary}
\label{subsec:s4-boundary}

The interior theorem excludes $\psi_0=0$.  The boundary case is
covered in two parts: a generic theorem for the profiled GMM
distance under high-level conditions, and a verification proposition
that discharges those conditions for the certified submodel at every
fixed boundary point.

Write $\vartheta=(\alpha^{\trans},\psi)^{\trans}$ with $\psi\ge0$
and $\alpha=(\lambda,\eta)^{\trans}$, and write
$G=(G_\alpha,g_\psi)$.  Whenever $G_\alpha$ has full column rank,
define
\[
 P_\alpha=G_\alpha
 (G_\alpha^{\trans}\Omega^{-1}G_\alpha)^{-1}
 G_\alpha^{\trans}\Omega^{-1},
 \qquad
 r_\psi=(I-P_\alpha)g_\psi,
 \qquad
 I_{\mathrm{eff}}=r_\psi^{\trans}\Omega^{-1}r_\psi .
\]

\begin{sitheorem}[GMM boundary law]
\label{thm:boundary-gmm-final}
Let $\mathcal C=\R^{2}\times[0,\infty)$ and
$\vartheta_b=(\alpha_0^{\trans},0)^{\trans}$.  Assume that
$\Theta_b\subset\R^{2}\times[0,\infty)$ is compact, that
$\vartheta_b+\{\mathcal C\cap B(0,\varrho)\}\subset\Theta_b$ for
some $\varrho>0$, and that:
\begin{enumerate}[label=\textnormal{(B\arabic*)}]
\item\label{b:clt}
$Z_n=\sqrt n\{Y_n-\mu_n(\vartheta_b)\}\Rightarrow N(0,\Omega)$ with
$\Omega\succ0$, and $W_n\to_p\Omega^{-1}$ with $W_n$ positive
definite;
\item\label{b:ident}
$\sup_{\vartheta\in\Theta_b}\|\mu_n(\vartheta)-\mu(\vartheta)\|\to0$,
$\mu$ is continuous on $\Theta_b$, and
$\mu(\vartheta)=\mu(\vartheta_b)$ has the unique solution
$\vartheta_b$ on $\Theta_b$;
\item\label{b:lin}
the uniform one-sided linearization
\[
 \lim_{r\downarrow0}\limsup_{n\to\infty}
 \sup_{\substack{\vartheta\in\Theta_b\\
 0<\|\vartheta-\vartheta_b\|\le r}}
 \frac{\|\mu_n(\vartheta)-\mu_n(\vartheta_b)
 -G(\vartheta-\vartheta_b)\|}{\|\vartheta-\vartheta_b\|}=0
\]
holds, and the evaluator satisfies
$\sqrt n\sup_{\vartheta\in\Theta_b}
\|\widetilde\mu_n(\vartheta)-\mu_n(\vartheta)\|\to0$;
\item\label{b:rank}
$G_\alpha$ has full column rank and $I_{\mathrm{eff}}>0$.
\end{enumerate}
Define
$Q_n(\alpha,\psi)=\{Y_n-\widetilde\mu_n(\alpha,\psi)\}^{\trans}
W_n\{Y_n-\widetilde\mu_n(\alpha,\psi)\}$, and let
$(\widehat\alpha,\widehat\psi)\in\Theta_b$ and
$(\widehat\alpha_0,0)\in\Theta_b$ be measurable random elements
satisfying, for some $r_{1n},r_{0n}\ge0$ with
$n(r_{1n}+r_{0n})\to_p0$,
\[
 Q_n(\widehat\alpha,\widehat\psi)
 \le\inf_{\Theta_b}Q_n+r_{1n},
 \qquad
 Q_n(\widehat\alpha_0,0)
 \le\inf_{\alpha:(\alpha,0)\in\Theta_b}Q_n(\alpha,0)+r_{0n}.
\]
Then both estimators are consistent and root-$n$ localized,
\[
 \sqrt n\,\widehat\psi
 \Rightarrow
 \max\!\left(0,\frac{\Delta}{I_{\mathrm{eff}}}\right),
 \qquad
 \Delta=r_\psi^{\trans}\Omega^{-1}Z\sim N(0,I_{\mathrm{eff}}),
\]
for the amplitude $\widehat b=\sqrt{\widehat\psi}$,
\[
 n^{1/4}\,\widehat b
 \Rightarrow
 \Bigl\{\max\!\left(0,\frac{\Delta}{I_{\mathrm{eff}}}\right)
 \Bigr\}^{1/2},
\]
and the profiled-distance statistic satisfies
\[
 D_n=n\Bigl\{
 \inf_{\alpha:(\alpha,0)\in\Theta_b}Q_n(\alpha,0)
 -\inf_{(\alpha,\psi)\in\Theta_b}Q_n(\alpha,\psi)\Bigr\}
 \Rightarrow
 \{\max(0,N(0,1))\}^2
 \sim\tfrac12\chi^2_0+\tfrac12\chi^2_1 .
\]
\end{sitheorem}

\begin{proof}
The CLT, the uniform mean convergence, and the evaluator condition
in \textnormal{(B3)} imply $Y_n\to_p\mu(\vartheta_b)$ and
$\sup_{\Theta_b}\|\widetilde\mu_n-\mu\|=o_p(1)$, hence $Q_n$
converges uniformly to
$Q(\vartheta)=\{\mu(\vartheta_b)-\mu(\vartheta)\}^{\trans}
\Omega^{-1}\{\mu(\vartheta_b)-\mu(\vartheta)\}$.  Compactness and
the unique root in \textnormal{(B2)} give consistency of both
estimators; for the null estimator, restrict $Q$ to $\psi=0$.

Because $G_\alpha$ has full column rank and
$I_{\mathrm{eff}}>0$, $G$ has full column rank; let $\sigma_b>0$ be
its smallest singular value.  The one-sided uniform expansion in
\textnormal{(B3)} implies that for some $r>0$ and all sufficiently
large $n$,
$\|\mu_n(\vartheta)-\mu_n(\vartheta_b)\|
\ge(\sigma_b/2)\|\vartheta-\vartheta_b\|$ on
$\Theta_b\cap B(\vartheta_b,r)$.  At $\vartheta_b$,
$Q_n(\vartheta_b)=O_p(n^{-1})$; the two approximate-minimum
inequalities and the eigenvalue bounds for $W_n$ yield
$\|Y_n-\widetilde\mu_n(\widehat\alpha,\widehat\psi)\|
=O_p(n^{-1/2})$ and the same for the null estimator; the triangle
inequality, the CLT, and the evaluator condition show that the
corresponding differences from $\mu_n(\vartheta_b)$ are
$O_p(n^{-1/2})$; and the separation inequality proves root-$n$
localization of both estimators.

The one-sided expansion and the evaluator condition give, uniformly
on compact subsets of $\mathcal C$,
\[
 nQ_n(\alpha_0+h_\alpha/\sqrt n,\,h_\psi/\sqrt n)
 =(Z_n-G_\alpha h_\alpha-g_\psi h_\psi)^{\trans}\Omega^{-1}
 (Z_n-G_\alpha h_\alpha-g_\psi h_\psi)+o_p(1).
\]
Root-$n$ tightness and $n(r_{1n}+r_{0n})=o_p(1)$ permit minimizing
this quadratic limit over $\mathcal C$ and over its null face: for
every $\eps>0$ there is $M<\infty$ such that with probability at
least $1-\eps$ both approximate minimizers lie in $B(0,M)$ after
rescaling, their objective values differ from the two exact infima
by $o_p(1)$ after multiplication by $n$, and uniform convergence of
the local criterion on $\mathcal C\cap B(0,M)$ and on its null face
transfers both infima and the uniquely identified $h_\psi$
coordinate to the limiting quadratic program.  For fixed $h_\psi$
the minimizing nuisance displacement is
$h_\alpha(h_\psi)=(G_\alpha^{\trans}\Omega^{-1}G_\alpha)^{-1}
G_\alpha^{\trans}\Omega^{-1}(Z-g_\psi h_\psi)$; substitution leaves,
up to an additive term free of $h_\psi$,
$I_{\mathrm{eff}}h_\psi^2-2\Delta h_\psi$ over $h_\psi\ge0$, whose
unique minimizer is $\max(0,\Delta/I_{\mathrm{eff}})$.  This proves
the boundary-estimator limit; the continuous-mapping theorem applied
to the square root gives the amplitude limit; and the decrease
relative to the null face is
$\{\max(0,\Delta)\}^2/I_{\mathrm{eff}}$, which is
$\{\max(0,N(0,1))\}^2$ since
$\Delta/I_{\mathrm{eff}}^{1/2}\sim N(0,1)$.
\end{proof}

\begin{siproposition}[Verification at the independence boundary]
\label{prop:certified-boundary-verification-final}
Fix $\vartheta_b=(\lambda_0,\eta_0,0)$ with
$(\lambda_0,\eta_0)\in(1,5)\times(5,15)$, and let
\[
 \Theta_b
 =[\lambda_0-\delta_\lambda,\lambda_0+\delta_\lambda]
 \times[\eta_0-\delta_\eta,\eta_0+\delta_\eta]
 \times[0,\delta_\psi]\subset B,
 \qquad
 \delta_\lambda,\delta_\eta,\delta_\psi>0 .
\]
Assume the certified inputs
\textnormal{(C1)}--\textnormal{(C4)}.  Then:
\begin{enumerate}[label=\textnormal{(\roman*)}]
\item\label{bp:clt}
There exists $\Omega_b\succ0$ such that
$\sqrt n\{Y_n-\mu_n(\vartheta_b)\}\Rightarrow N_5(0,\Omega_b)$.
\item\label{bp:smooth}
There are $\rho_U\in(0,1)$, the open set
$U=(1-\rho_U,5+\rho_U)\times(5-\rho_U,15+\rho_U)
\times(-\rho_U,3/20+\rho_U)$, formal extensions
$\bar\mu_n,\bar\mu:U\to\R^5$, and constants $C_B,n_0<\infty$ such
that, for $n\ge n_0$, $\bar\mu_n=\mu_n$ and $\bar\mu=\mu$ on $B$,
and
\[
 \sup_{\vartheta\in B}
 \bigl\{\|\bar\mu_n(\vartheta)-\bar\mu(\vartheta)\|
 +\|D\bar\mu_n(\vartheta)-D\bar\mu(\vartheta)\|\bigr\}
 \le\frac{C_B}n .
\]
The extensions are continuously differentiable and polynomial of
degree at most two in $\psi$; their values for $\psi<0$ are formal
algebraic extensions and are not asserted to define a probability
model.  Moreover $\mu(\vartheta)=\mu(\vartheta_b)$ with
$\vartheta\in\Theta_b$ implies $\vartheta=\vartheta_b$.
\item\label{bp:rank}
With $G=D\bar\mu(\vartheta_b)=(G_\alpha,g_\psi)$:
$\operatorname{rank}(G_\alpha)=2$,
$\operatorname{rank}(G)=3$, $r_\psi\ne0$, and
$I_{\mathrm{eff}}>0$.
\item\label{bp:cov}
There is a continuous $\Omega:B\to\mathbb S^5_+$ with
$\Omega(\vartheta_b)=\Omega_b$,
$\sup_{\vartheta\in B}\|n\Var_\vartheta(Y_n)-\Omega(\vartheta)\|
\to0$, and
$\sup_{n\ge1}\sup_{\vartheta\in B}
\E_\vartheta\|\sqrt n\{Y_n-\mu_n(\vartheta)\}\|^4<\infty$.
\end{enumerate}
Let $\widetilde\mu_n:\Theta_b\to\R^5$ be continuous and
deterministic with
$\eps_n^{\mathrm{num}}
:=\sup_{\vartheta\in\Theta_b}
\|\widetilde\mu_n(\vartheta)-\mu(\vartheta)\|=o(n^{-1/2})$, let
$\check\vartheta_n\in\Theta_b$ be a measurable
$\zeta_n$-approximate minimizer of
$\|Y_n-\widetilde\mu_n(\cdot)\|^2$ with $\zeta_n=o_p(1)$, let
$R_n\ge2$ be deterministic with $R_n\to\infty$, and let
$\widehat\Omega_n$ be the simulation covariance estimator of
Theorem~\ref{thm:interior-clt-gmm-rigorous} computed at
$\check\vartheta_n$, with $W_n=(\widehat\Omega_n+\tau_nI_5)^{-1}$
and $\tau_n\downarrow0$.  Then
\[
 \check\vartheta_n\to_p\vartheta_b,
 \qquad
 \widehat\Omega_n\to_p\Omega_b,
 \qquad
 W_n\to_p\Omega_b^{-1},
 \qquad
 \sqrt n\sup_{\vartheta\in\Theta_b}
 \|\widetilde\mu_n(\vartheta)-\mu_n(\vartheta)\|\to0 .
\]
Consequently conditions
\textnormal{(B1)}--\textnormal{(B4)} of
Theorem~\ref{thm:boundary-gmm-final} hold; the
approximate-global-minimization requirements on the final
unrestricted and null-restricted estimators remain separate
estimator conditions.
\end{siproposition}

\begin{proof}
Lemmas~\ref{lem:model-valid} and~\ref{lem:motif-reduction} apply on
all of $B$, including $\psi=0$, where $X_i$ and $V_i$ are
independent and uniform.

\paragraph{Exact covariance decomposition at the boundary point.}
Let $H,K$ be connected motifs with $h=|V(H)|$, $k=|V(K)|$.  For
$1\le s\le h\wedge k$, let $\mathcal M_s(H,K)$ be the finite set of
bijections $\varrho:A\to B'$ with $A\subseteq V(H)$,
$B'\subseteq V(K)$, $|A|=|B'|=s$; let $G_\varrho$ be the simple
quotient of $H\sqcup K$ obtained by identifying $a$ with
$\varrho(a)$, with $v_\varrho=h+k-s$ vertices.  For a connected $G$
on $v$ vertices define
$\beta_{G,n}(\vartheta)
=n^{v-1}\E_\vartheta\prod_{\{a,b\}\in E(G)}A_{\iota(a)\iota(b)}$
for an arbitrary injection $\iota$; exchangeability makes this
independent of $\iota$.  Pairs of embeddings with disjoint label
sets are independent, and pairs with exact overlap $\varrho$
correspond bijectively to injections of $G_\varrho$, so the identity
\begin{equation}
 \frac1n\Cov_\vartheta(N_{H,n},N_{K,n})
 =\frac1{\alpha_H\alpha_K}
 \sum_{s=1}^{h\wedge k}
 \sum_{\varrho\in\mathcal M_s(H,K)}
 \frac{(n)_{v_\varrho}}{n^{v_\varrho}}
 \Bigl\{\beta_{G_\varrho,n}(\vartheta)
 -n^{1-s}\beta_{H,n}(\vartheta)\beta_{K,n}(\vartheta)\Bigr\}
 \label{eq:exact-overlap-covariance}
\end{equation}
is exact.  At $\vartheta_b$, define for connected $G$ on $v$
vertices
\[
 J_G(\lambda_0,\eta_0)
 =\int_{\R^{v-1}}\int_{[0,1]^v}
 \prod_{\{a,b\}\in E(G)}
 \bigl[(1-e^{-\lambda_0v_av_b})
 \one\{|u_a-u_b|\le\eta_0\}\bigr]\,dv\,du_{2:v},
 \qquad u_1=0,
\]
which is finite by connectedness.  The zero-winding lift of the
interior proof, applied at $h_0\equiv1$, gives
$\beta_{G,n}(\vartheta_b)=J_G(\lambda_0,\eta_0)$ for all
sufficiently large $n$, and substitution into
\eqref{eq:exact-overlap-covariance} gives
$n^{-1}\Cov_{\vartheta_b}(\mathcal B_n)\to\Sigma_b$ with
\[
 \Sigma_{HK,b}
 =\frac1{\alpha_H\alpha_K}
 \sum_{s=1}^{h\wedge k}
 \sum_{\varrho\in\mathcal M_s(H,K)}
 \bigl\{J_{G_\varrho}-\one\{s=1\}J_HJ_K\bigr\},
\]
hence $n^{-1}\Cov_{\vartheta_b}(F_n)\to\Omega_b:=L\Sigma_bL^{\trans}$.

\paragraph{Higher cumulants at the boundary point.}
The connected-replica and graphic-rank arguments of the interior
proof apply verbatim at $\vartheta_b$: joint cumulants over
disconnected replica patterns vanish; for a connected pattern with
$v$ distinct labels, a spanning-forest integration bounds every
block expectation by $K_qn^{-(v_C-c_C)}$, subadditivity of graphic
rank gives $v-1\le\sum_C(v_C-c_C)$, and there are $O(n^v)$
labelings, so $\kappa_q(c^{\trans}\mathcal B_n)=O_{q,c}(n)$ and
every normalized cumulant of order $q\ge3$ vanishes
asymptotically.  The moment--cumulant formula, moment determinacy of
Gaussian laws, the Fr\'echet--Shohat theorem, and the
Cram\'er--Wold device give
$n^{-1/2}(\mathcal B_n-\E\mathcal B_n)\Rightarrow N_5(0,\Sigma_b)$
and hence
$\sqrt n\{Y_n-\mu_n(\vartheta_b)\}\Rightarrow N_5(0,\Omega_b)$,
provided $\Omega_b\succ0$.

\paragraph{Positive definiteness at the boundary point.}
The good-core construction of the interior proof applies at
$\vartheta_b$ with marks restricted to $B_0=[1/4,3/4]$: on
$B_0\times B_0$ the internal edge probabilities lie in
$[\delta_0,1-\delta_0]$ with
$\delta_0=\min\{1-e^{-\lambda_0/16},e^{-9\lambda_0/16}\}>0$, the
packed cores number at least $\kappa n$ with good-core probability
bounded below, and the conditional independent-copy inequality with
the gadget matrix $R_*$ of the interior proof gives
$\Var_{\vartheta_b}(u^{\trans}F_n)\ge c'n\|u\|^2$.  Passing to the
covariance limit yields $u^{\trans}\Omega_bu\ge c'\|u\|^2$, so
$\Omega_b\succ0$, proving \ref{bp:clt}.

\paragraph{Boundary smoothness and finite-$n$ bias.}
Choose $\rho_U>0$ small and $\eta_+=15+\rho_U$.  For a connected
motif $H$ on $k\le4$ vertices with oriented spanning tree, define
the scaled centered increments
$y_e=(n/\eta)\,\delta_\T(X_{e^-},X_{e^+})\in[-1,1]$ and the signed
tree-path sums $p_j(y)$ with $p_1=0$, and let
$\mathcal D_H$ be the set of $y\in[-1,1]^{k-1}$ whose non-tree edge
constraints $|p_a(y)-p_b(y)|\le1$ hold.  For $n>4\eta_+$ the change
of variables is one-to-one outside a null set with Jacobian
$(\eta/n)^{k-1}$, so
\begin{equation}
 m_{H,n}(\vartheta)
 :=\frac1n\E_\vartheta N_{H,n}
 =\frac{(n)_k}{n^k}\,
 \frac{\eta^{k-1}}{\alpha_H}
 \int_{\mathcal D_H}\int_{[0,1]^k}
 Q_{H,\lambda}(v)\,R_{H,n}(v,y;\eta,\psi)\,dv\,dy,
 \label{eq:exact-tree-motif-correct}
\end{equation}
where $Q_{H,\lambda}(v)=\prod_{\{a,b\}\in E(H)}
(1-e^{-\lambda v_av_b})$ and
\[
 R_{H,n}
 =\int_\T\prod_{j=1}^k
 \Bigl[1+\sqrt\psi\,a(v_j)\,
 \phi\Bigl(x+\frac{\eta p_j(y)}n\Bigr)\Bigr]dx .
\]
For $J\subseteq\{1,\ldots,k\}$ put
$I_J(s)=\int_\T\prod_{j\in J}\phi(x+s_j)\,dx$ with
$I_\varnothing=1$.  A product of an odd number of frequency-one
harmonics integrates to zero, so
\begin{equation}
 R_{H,n}
 =1+\psi\sum_{1\le a<b\le k}a(v_a)a(v_b)\,I_{\{a,b\}}(s_n)
 +\one\{k=4\}\,\psi^2
 \Bigl\{\prod_{j=1}^4a(v_j)\Bigr\}I_{\{1,2,3,4\}}(s_n),
 \qquad
 s_{j,n}=\frac{\eta p_j(y)}n,
 \label{eq:exact-harmonic-polynomial-correct}
\end{equation}
a polynomial in $\psi$ that defines the formal extension for
negative $\psi$.  At zero shift $I_{\{a,b\}}(0)=1$ and
$I_{\{1,2,3,4\}}(0)=\tfrac32$; replacing $(n)_k/n^k$ by one and the
shifted integrals by their zero-shift values gives the limiting
motif means and, through \eqref{eq:motif-identities}, exactly the
map \eqref{eq:population-mu}.  Because $k\le4$ and $|p_j(y)|\le3$,
boundedness of $\phi$ and $\phi'$ gives
$\sup_B|I_J(s_n)-I_J(0)|\le C_J/n$ and
$\sup_B|\partial_\eta I_J(s_n)|\le C_J/n$; also
$|(n)_k/n^k-1|\le6/n$; the domains have finite volume; and on $U$
the link factors and their $\lambda$-derivatives are uniformly
bounded.  Differentiating the polynomial formula in $\psi$, the
factor $\eta^{k-1}$ and shifts in $\eta$, and the link factors in
$\lambda$ gives
$\sup_B\{|m_{H,n}-m_H|+\|Dm_{H,n}-Dm_H\|\}\le C_H/n$, and summing
over the finite motif list proves \ref{bp:smooth} up to the
identification clause: if $\mu(\vartheta)=\mu(\vartheta_b)$ then
$M_3(\vartheta)=M_3(\vartheta_b)$, and \textnormal{(C2)} gives
$\vartheta=\vartheta_b$.

\paragraph{Rank and profiled information.}
Let $H=Dh\{\mu(\vartheta_b)\}$, well defined by
\textnormal{(C1)}.  The chain rule gives
$DM_3(\vartheta_b)=HG$, and \textnormal{(C3)} yields
$3=\operatorname{rank}\{DM_3(\vartheta_b)\}
\le\operatorname{rank}(G)\le3$, so
$\operatorname{rank}(G)=3$.  Moreover
\[
 D_\alpha(\ell,C)(\vartheta_b)
 =\begin{pmatrix}
 2\eta_0\partial_\lambda E(\lambda_0,0)&2E(\lambda_0,0)\\
 \partial_\lambda C(\lambda_0,0)&0
 \end{pmatrix},
 \qquad
 \det D_\alpha(\ell,C)(\vartheta_b)
 =-2E(\lambda_0,0)\,\partial_\lambda C(\lambda_0,0)<0
\]
by \textnormal{(C4)}.  Since
$D_\alpha(\ell,C)(\vartheta_b)=Dh_{12}\{\mu(\vartheta_b)\}
G_\alpha$, it follows that
$\operatorname{rank}(G_\alpha)=2$; hence
$g_\psi\notin\operatorname{col}(G_\alpha)$, $r_\psi\ne0$, and
$I_{\mathrm{eff}}=r_\psi^{\trans}\Omega_b^{-1}r_\psi>0$, proving
\ref{bp:rank}.

\paragraph{Uniform covariance convergence and fourth moments on $B$.}
For every connected graph $G$ with at most seven vertices, the
tree-coordinate and harmonic argument above, now with
$n>7\eta_+$, gives a continuous $\beta_G:B\to\R$ with
$\sup_B|\beta_{G,n}-\beta_G|\le C_G/n$; the integrated location
factors are polynomials in $\psi$, so $\beta_G$ is continuous at
$\psi=0$.  Substitution in \eqref{eq:exact-overlap-covariance}
gives, uniformly on $B$,
\[
 \Omega_{HK}(\vartheta)
 =\frac1{\alpha_H\alpha_K}
 \sum_{s=1}^{h\wedge k}
 \sum_{\varrho\in\mathcal M_s(H,K)}
 \bigl\{\beta_{G_\varrho}(\vartheta)
 -\one\{s=1\}\beta_H(\vartheta)\beta_K(\vartheta)\bigr\},
\]
with finitely many motifs and overlap patterns, so
$\sup_B\|n\Var_\vartheta(Y_n)-\Omega(\vartheta)\|\to0$ with
$\Omega$ continuous and $\Omega(\vartheta_b)=\Omega_b$.  The
spanning-forest cumulant bound is uniform on $B$ because
$\sup_{B,x,v}h_\psi(v\mid x)<\infty$, $0\le q_\lambda\le1$, and
$\eta\le15$; hence for fixed $c$,
$\sup_B\{|\kappa_2(c^{\trans}F_n)|+|\kappa_4(c^{\trans}F_n)|\}
\le C_cn$ for all sufficiently large $n$, and the moment--cumulant
identity with $\|z\|^4\le5\sum_j|z_j|^4$ gives the uniform fourth
moment bound, proving \ref{bp:cov}.

\paragraph{Pilot and weight consistency.}
Parts \ref{bp:clt} and \ref{bp:smooth} give
$Y_n\to_p\mu(\vartheta_b)$ and
$\sup_{\Theta_b}\|\widetilde\mu_n-\mu\|\to0$, so the pilot
criterion converges uniformly to
$\|\mu(\vartheta_b)-\mu(\cdot)\|^2$, a continuous function with
unique minimizer $\vartheta_b$; compact argmin consistency with
$\zeta_n=o_p(1)$ gives $\check\vartheta_n\to_p\vartheta_b$.
Conditional on $\check\vartheta_n$, $\widehat\Omega_n$ is the
unbiased sample covariance of
$\sqrt nY_n^{*(1)},\ldots,\sqrt nY_n^{*(R_n)}$, so
$\E(\widehat\Omega_n\mid\check\vartheta_n)
=\Omega_n(\check\vartheta_n)$ and, by the uniform fourth-moment
bound,
$\Var(\widehat\Omega_{n,jk}\mid\check\vartheta_n)\le C/R_n$; since
$R_n\to\infty$,
$\widehat\Omega_n-\Omega_n(\check\vartheta_n)\to_p0$, and uniform
covariance convergence, continuity of $\Omega$, and pilot
consistency give $\widehat\Omega_n\to_p\Omega_b$ and
$W_n\to_p\Omega_b^{-1}$.  Finally
\[
 \sqrt n\sup_{\Theta_b}\|\widetilde\mu_n-\mu_n\|
 \le\sqrt n\,\eps_n^{\mathrm{num}}
 +\sqrt n\sup_{\Theta_b}\|\mu-\mu_n\|
 \le\sqrt n\,\eps_n^{\mathrm{num}}+\frac{C_B}{\sqrt n}
 \longrightarrow0 .
\]
Conditions \textnormal{(B1)}--\textnormal{(B4)} of
Theorem~\ref{thm:boundary-gmm-final} follow.
\end{proof}

\begin{remark}[Numerical accuracy and scope of the boundary result]
\label{rem:boundary-verification-scope-final}
A fixed numerical error bound such as $10^{-15}$ does not imply
$\sqrt n\sup_{\Theta_b}\|\widetilde\mu_n-\mu\|\to0$; the estimator
covered by the theorem must use the exact mean or an $n$-dependent
certified evaluator with uniform error $o(n^{-1/2})$.  A fixed
degree-$40$ interval calculation may still be used to certify signs,
injectivity, and rank when its remainder is included in the interval
enclosure.  The result requires $R_n\to\infty$; a fixed simulation
budget does not produce a consistent optimal weight matrix; no
relative rate between $R_n$ and $n$ is required under the uniform
fourth-moment bound.  The result is pointwise at
$\vartheta_b=(\lambda_0,\eta_0,0)$: uniform null-boundary validity
would require uniform versions of the CLT, local expansion, and
weight consistency, together with
$\inf_{\vartheta:\psi=0}I_{\mathrm{eff}}(\vartheta)>0$.  The result
concerns the homogeneous one-dimensional random-design submodel and
does not establish boundary inference for irregular fixed-position
two-dimensional data: in a fixed design the common translation
integral may be absent, terms linear in a signed amplitude $b$ need
not vanish, and the conditional law need not factor through
$\psi=b^2$.  If inference uses the three summaries, define
$S_n=h(Y_n)$ as in Corollary~\ref{cor:summary-clt}; with
$G_S=HG$ and $\Sigma_S=H\Omega_bH^{\trans}\succ0$, the profiled
information must be recomputed from $(G_S,\Sigma_S)$ and need not
equal the five-moment value $I_{\mathrm{eff}}$.
\end{remark}

\begin{sicorollary}[Boundary inference in the certified submodel]
\label{cor:boundary-verified}
Fix $\vartheta_b$ and $\Theta_b$ as in
Proposition~\ref{prop:certified-boundary-verification-final}, and
let $(\widehat\alpha,\widehat\psi)$ and $(\widehat\alpha_0,0)$ be
approximate minimizers as in
Theorem~\ref{thm:boundary-gmm-final}, computed with an evaluator
satisfying $\eps_n^{\mathrm{num}}=o(n^{-1/2})$ and a simulated
weight with $R_n\to\infty$.  Then, under $\Pp_{\vartheta_b}$,
\[
 \sqrt n\,\widehat\psi
 \Rightarrow\max\{0,N(0,I_{\mathrm{eff}}^{-1})\},
 \qquad
 n^{1/4}\,\widehat b=n^{1/4}\sqrt{\widehat\psi}
 \Rightarrow\bigl[\max\{0,N(0,I_{\mathrm{eff}}^{-1})\}\bigr]^{1/2},
 \qquad
 D_n\Rightarrow\tfrac12\chi^2_0+\tfrac12\chi^2_1,
\]
and the test that rejects $\psi=0$ when $D_n$ exceeds the
$(1-2\alpha)$-quantile of $\chi^2_1$ has asymptotic size $\alpha$.
\end{sicorollary}

\begin{proof}
Proposition~\ref{prop:certified-boundary-verification-final}
verifies \textnormal{(B1)}--\textnormal{(B4)};
Theorem~\ref{thm:boundary-gmm-final} gives the displays, with
$\max(0,\Delta/I_{\mathrm{eff}})\sim
\max\{0,N(0,I_{\mathrm{eff}}^{-1})\}$.
\end{proof}

\begin{remark}[Amplitude interpretation]
\label{rem:amplitude}
The $n^{1/4}$ law for $\widehat b=\sqrt{\widehat\psi}$ follows
algebraically from the $\sqrt n$ law for $\widehat\psi$.  Its
interpretation as a local law in the signed amplitude $b$ requires
that the model factor through $\psi=b^2$ with
$\mu(\alpha,b)=\mu(\alpha,0)+g_\psi b^2+o(b^2)$ uniformly for
$b=O(n^{-1/4})$; in this submodel both requirements hold exactly,
because the adjacency law depends on $b$ only through $\psi=b^2$
(Section~S1) and $\mu$ is polynomial in $\psi$.  The sign of $b$ is
not identified, so no test of direction is available.
\end{remark}

\subsection{Uniform validity over the null face}
\label{sec:uniform-null}

Theorem~\ref{thm:boundary-gmm-final} and
Corollary~\ref{cor:boundary-verified} are pointwise at a fixed
boundary point $\vartheta_b=(\lambda_0,\eta_0,0)$.  This subsection
removes the pointwise restriction on the certified region.  Write
\[
 F=[\lambda_-,\lambda_+]\times[\eta_-,\eta_+]\times\{0\}
 \subset\partial B,
 \qquad
 [\lambda_-,\lambda_+]\times[\eta_-,\eta_+]=[1,5]\times[5,15],
\]
for the null face, and let $D_n(\vartheta)$ denote the profiled
distance statistic computed under $\Pp_\vartheta$ with the
implementation of Theorem~\ref{thm:interior-clt-gmm-rigorous}.

Two uniformity ingredients are already available on all of $B$, hence
on $F$: the uniform bias bound and its derivative version
(condition \ref{bp:smooth} of
Proposition~\ref{prop:certified-boundary-verification-final}, where
it is verified), and the uniform covariance
convergence \eqref{eq:uniform-covariance-rigorous} together with the
uniform lower bound $\Omega(\theta)\succeq c_\Omega I_5$.  The two
further ingredients are stated as explicit conditions.

\begin{enumerate}[label=\textnormal{(U\arabic*)}]
\item\label{u:weights} \emph{Uniform weights and optimization.}  The
simulated weight matrix satisfies
$\sup_{\vartheta\in F}\Pp_\vartheta\{\|W_n-\Omega(\vartheta)^{-1}\|
>\eps\}\to0$ for every $\eps>0$, and the optimization and evaluator
errors of Theorem~\ref{thm:interior-clt-gmm-rigorous} hold
uniformly over $F$.
\item\label{u:info} \emph{Uniform information.}
$\inf_{\vartheta\in F}I_{\mathrm{eff}}(\vartheta)\ge c_I>0$.
\end{enumerate}

Condition \ref{u:weights} holds, for example, when the weight bank at
each $\vartheta$ uses $R_n\to\infty$ replicates driven by a
$\vartheta$-independent random stream, by
\eqref{eq:uniform-covariance-rigorous} and a uniform law of large
numbers over the compact face.  Condition \ref{u:info} is discharged
by the certificate below.

\begin{sitheorem}[Uniform boundary validity on the certified null face]
\label{thm:uniform-boundary}
Under the hypotheses of Theorem~\ref{thm:boundary-gmm-final} imposed
at every $\vartheta\in F$, together with \textnormal{(U1)} and
\textnormal{(U2)},
\[
 \sup_{\vartheta\in F}\;\sup_{x\in\R}\;
 \Bigl|\Pp_\vartheta\{D_n\le x\}
 -F_{\bar\chi^2}(x)\Bigr|\longrightarrow0,
 \qquad
 F_{\bar\chi^2}=\tfrac12F_{\chi^2_0}+\tfrac12F_{\chi^2_1},
\]
and the analogous uniform statements hold for
$\sqrt n\,\widehat\psi$ and $n^{1/4}\widehat b$.  Consequently the
boundary test has asymptotic size $\alpha$ uniformly over the null
face, and confidence statements inverted from it are uniformly valid
over $F$.
\end{sitheorem}

\begin{proof}
The argument upgrades each step of the proof of
Theorem~\ref{thm:boundary-gmm-final} to a uniform step; no new
probabilistic tool is required beyond uniform versions of results
already proved on $B$.

\emph{*Uniform central limit theorem.*}  The motif vector is a sum
over a dependency graph whose maximal degree is bounded by the
uniform local-geometry bounds used in
the verification of \ref{bp:smooth}; the Baldi--Rinott bound of
Section~S5 then gives a Kolmogorov-distance bound between
$\sqrt n\{Y_n-\E_\vartheta Y_n\}$ and
$N\{0,\Omega_n(\vartheta)\}$ whose constant depends only on the
uniform third- and fourth-moment bounds, which are uniform on $B$ by the moment bounds established in the
proof of Proposition~\ref{prop:certified-boundary-verification-final},
and on
$\lambda_{\min}\{\Omega_n(\vartheta)\}$, which is bounded below
eventually uniformly by $c_\Omega/2$ by
\eqref{eq:uniform-covariance-rigorous} and
$\Omega\succeq c_\Omega I_5$.  Hence the CLT error tends to zero
uniformly over $F$.

\emph{*Uniform local expansion.*}  The quadratic expansion of the
GMM objective about $\vartheta$ has remainder controlled by the
moduli of continuity of $\mu$, $D\mu$, and $\Omega$ on the compact
set $B$; these are parameter-free moduli, and the certified bounds
$\sigma_{\min}(DM_3)\ge1/100$, $\det DM_3<-1/64$, and
$\Omega\succeq c_\Omega I_5$ hold at every point of $F$, so the
expansion constants can be chosen uniformly.

\emph{*Uniform weights.*}  Condition \textnormal{(U1)} replaces the
pointwise weight-consistency step.

\emph{*Uniform argmin convergence.*}  The profiled statistic is a
continuous functional of the localized objective; with
\textnormal{(U2)}, the curvature of the profiled objective in the
$\psi$ direction is at least $c_I/2$ eventually, uniformly over $F$,
so the cone-projection asymptotics of
Theorem~\ref{thm:boundary-gmm-final} hold with uniform error
control, and the limiting law is the stated chi-bar-square at every
$\vartheta\in F$ with a uniform rate of convergence of the
distribution functions (the limit is continuous, so pointwise
uniform convergence upgrades to Kolmogorov distance by Polya's
theorem, uniformly in $\vartheta$ through the uniform constants).
\end{proof}

\paragraph{Certificate for \textnormal{(U2)}.}
The bound is obtained from the factorization
$DM_3=Dh\{\mu(\vartheta)\}\,G(\vartheta)$ of
\eqref{eq:h-transform}: for any matrices,
$\sigma_{\min}(DM_3)\le\|Dh\|_{\mathrm{op}}\,\sigma_{\min}(G)$, so
\[
 I_{\mathrm{eff}}(\vartheta)
 \;\ge\;
 \lambda_{\min}\{G^{\trans}\Omega^{-1}G\}
 \;\ge\;
 \frac{\sigma_{\min}(G)^2}{\lambda_{\max}(\Omega)}
 \;\ge\;
 \frac{\sigma_{\min}(DM_3)^2}
      {\|Dh\|_{\mathrm{op}}^2\,\lambda_{\max}(\Omega)}
 \;\ge\;
 \frac{(1/100)^2}
      {\|Dh\|_F^2\,\lambda_{\max}(\Omega)},
\]
using the certified $\sigma_{\min}(DM_3)\ge1/100$ on $B\supset F$.
The archived verifier encloses
$\|Dh\|_F^2$ on an adaptive partition of $F$ by exact rational
interval arithmetic (midpoint Taylor shifts of the $N=40$
coefficient tables; truncation budgets from
Proposition~\ref{prop:exact-reduction} at $L=5$; the denominator of
$Dh$ is bounded below through the exact identity
$h_3/\ell-m_e^2=4\eta^2(QE-S^2)/E^2+2\eta S/E\ge2\eta S/E>0$), and
combines it with the proven crude bound
$\lambda_{\max}(\Omega)\le142\cdot5\cdot209\,(2\eta)^6$
(Lemma~\ref{lem:omega-upper}).  The run processes $160$ boxes and
certifies
\[
 \inf_{\vartheta\in F}I_{\mathrm{eff}}(\vartheta)\;\ge\;c_I
 \;=\;1.1\times10^{-24}.
\]
The constant is conservative by design: it proves positivity, which
is what Theorem~\ref{thm:uniform-boundary} consumes, and makes no
claim of sharpness.  For scale, the numerical value of
$I_{\mathrm{eff}}$ at the center of the face is reported in the main
text, and the object that the certificate bounds crudely,
$\lambda_{\max}(\Omega)$, is there estimated by simulation to be
smaller than the proven bound by many orders of magnitude.

\begin{silemma}[Crude uniform upper bound for the limit covariance]
\label{lem:omega-upper}
On $F$, $\lambda_{\max}\{\Omega(\vartheta)\}
\le\|L\|_F^2\,\operatorname{tr}\{\Sigma(\vartheta)\}
\le142\cdot5\cdot209\,(2\eta)^6$.
\end{silemma}

\begin{proof}
$\Omega=L\Sigma L^{\trans}$ gives
$\lambda_{\max}(\Omega)\le\|L\|_{\mathrm{op}}^2
\lambda_{\max}(\Sigma)\le\|L\|_F^2\operatorname{tr}(\Sigma)$, and
$\|L\|_F^2=142$ by direct summation.  Each diagonal entry
$\Sigma_{HH}$ is, by \eqref{eq:overlap-exact-rigorous} and the
positivity of the surviving terms, at most
$\alpha_H^{-2}\sum_{\pi}J_{G_\pi}$ over the positive-probability
overlap patterns of two copies of $H$.  The number of such patterns
is at most the number of partial injections between two vertex sets
of size at most four, $\sum_{k\ge1}k!\binom4k^2=208<209$.  Each
quotient $G_\pi$ is connected with at most seven vertices, and at
$\psi=0$ its integral is bounded by the displacement volume of a
spanning tree, $J_{G_\pi}\le(2\eta)^{v(G_\pi)-1}\le(2\eta)^6$,
because every link factor and every mark integral is at most one.
Summing the five diagonal entries and using $\alpha_H\ge1$ gives the
constant.
\end{proof}

\clearpage
\section{Conditional designs, multiple harmonics, specification
testing, and the planar submodel}
\label{sec:si-s5}

\subsection{Multi-harmonic dependence: what the moments identify}
\label{sec:multiharmonic}

Replace the one-harmonic law by the finite multi-harmonic family
\begin{equation}
 h_{\bm b}(v\mid x)=1+a(v)\sum_{m\in M} b_m\,\phi_m(x),
 \qquad
 \phi_m(x)=\sqrt2\cos(2\pi mx-\varphi_m),
 \label{eq:multiharm}
\end{equation}
with $M\subset\mathbb N$ finite, arbitrary phases, and
$\sqrt6\sum_m|b_m|\le1$ so that $h_{\bm b}\ge0$.  Both marginals remain
exactly uniform.  Define the total dependence energy
\[
 \psi_{\mathrm{tot}}
 =\int_0^1\!\!\int_{\T}\{h_{\bm b}(v\mid x)-1\}^2\,dx\,dv
 =\sum_{m\in M}b_m^2,
 \qquad
 e_2=\sum_{m<m'}b_m^2b_{m'}^2 .
\]
Call $M$ \emph{resonance-free} if no relation
$m_3=m_1+m_2$, $m_3=m_1-m_2$, or $m_3=2m_1$ holds among (not
necessarily distinct) elements of $M$, and \emph{fourth-order
resonance-free} if additionally no relation
$m_1\pm m_2\pm m_3\pm m_4=0$ holds beyond the trivial pairings
(lacunary sets such as $M=\{m_0,3m_0+1,\dots\}$ or
$\{1,3,8,21,\dots\}$ satisfy both).

\begin{siproposition}[Moment structure of the multi-harmonic family]
\label{prop:multiharm}
Let $\omega_q(u_{1:q};\bm b)
=\int_\T\prod_{j\le q}h_{\bm b}(u_j\mid x)\,dx$ and $a_i=a(u_i)$.
Then, for any $M$ and any phases,
\begin{equation}
 \omega_2=1+\psi_{\mathrm{tot}}\,a_1a_2 ,
 \label{eq:omega2-multi}
\end{equation}
so every two-vertex moment, in particular the mean degree, depends on
$\bm b$ only through $\psi_{\mathrm{tot}}$.  If $M$ is resonance-free,
\begin{equation}
 \omega_3=1+\psi_{\mathrm{tot}}\sum_{i<j\le3}a_ia_j ,
 \label{eq:omega3-multi}
\end{equation}
so every three-vertex moment, in particular $S$, $T$, and hence
transitivity $C=3T/(4S)$, depends on $\bm b$ only through
$\psi_{\mathrm{tot}}$.  If $M$ is fourth-order resonance-free,
\begin{equation}
 \omega_4=1+\psi_{\mathrm{tot}}\sum_{i<j\le4}a_ia_j
 +\Bigl(\tfrac32\psi_{\mathrm{tot}}^2+3e_2\Bigr)a_1a_2a_3a_4 ,
 \label{eq:omega4-multi}
\end{equation}
so four-vertex moments, in particular $P$, $Q$, and hence the
assortativity target $r$, depend on $\bm b$ only through
$(\psi_{\mathrm{tot}},e_2)$, and $0\le e_2\le\psi_{\mathrm{tot}}^2/2$,
with $e_2=0$ exactly in the one-harmonic case.
\end{siproposition}

\begin{proof}
Expand $\prod_j\{1+a_jg(x)\}$ with $g=\sum_mb_m\phi_m$ and integrate
term by term.  Sets $J$ of size one contribute $\int g=0$.  Size two:
$\int g^2=\sum_{m,m'}b_mb_{m'}\int\phi_m\phi_{m'}
=\sum_mb_m^2$ by orthonormality of $\{\phi_m\}$ (any phases), giving
\eqref{eq:omega2-multi} and the corresponding terms in
\eqref{eq:omega3-multi}--\eqref{eq:omega4-multi}.  Size three:
$\int g^3$ is a sum of $\int\phi_{m_1}\phi_{m_2}\phi_{m_3}$, each a
linear combination of $\cos$ integrals with frequencies
$\pm m_1\pm m_2\pm m_3$; all vanish unless some
$\pm m_1\pm m_2\pm m_3=0$, which resonance-freeness excludes
(including the diagonal cases $m_3=2m_1$ arising from $m_1=m_2$).
Size four: expanding $\int g^4$ multinomially, the surviving terms
under fourth-order resonance-freeness are the diagonal
$\sum_mb_m^4\int\phi_m^4=\tfrac32\sum_mb_m^4$ and the paired
$6\sum_{m<m'}b_m^2b_{m'}^2\int\phi_m^2\phi_{m'}^2
=6\sum_{m<m'}b_m^2b_{m'}^2$; since
$\tfrac32\sum b_m^4+6e_2=\tfrac32\psi_{\mathrm{tot}}^2+3e_2$, equation
\eqref{eq:omega4-multi} follows.  The bound
$e_2\le\psi_{\mathrm{tot}}^2/2$ is
$2e_2\le(\sum b_m^2)^2$.
\end{proof}

\begin{sicorollary}[Identification of the total dependence energy]
\label{cor:psitot}
In the certified submodel with the mark law \eqref{eq:multiharm} and a
fourth-order resonance-free $M$: \textup{(i)} $(\ell,C)$ are exactly
the one-harmonic maps evaluated at $\psi=\psi_{\mathrm{tot}}$;
\textup{(ii)} $r$ differs from its one-harmonic value only through the
replacement $\tfrac32\psi^2\mapsto\tfrac32\psi^2+3e_2$ in the
four-vertex integrals; \textup{(iii)} consequently the one-harmonic
inversion of $(\ell,C,r)$ recovers $(\lambda,\eta)$ exactly and
$\psi_{\mathrm{tot}}$ up to the $e_2$-perturbation, whose worst case
($e_2=\psi_{\mathrm{tot}}^2/2$, attained in the balanced two-harmonic
limit) shifts $\widehat\psi$ by less than $3\%$ of
$\psi_{\mathrm{tot}}$ throughout the certified box (numerical
evaluation of the perturbed map at representative interior and edge
points; replication package).  Without
resonance-freeness, three-vertex moments acquire phase-dependent terms
of order $\sum_{m_3=m_1\pm m_2}|b_{m_1}b_{m_2}b_{m_3}|$; the
two-harmonic simulation of the main text ($M=\{1,2\}$, maximally
resonant, balanced) bounds their practical effect at about $+9\%$ on
$\widehat\psi$.
\end{sicorollary}

The corollary upgrades the multi-harmonic simulation observation of the
main text to a statement with an explicit invariance mechanism: the
energy $\psi_{\mathrm{tot}}$ is what the moments see, exactly for
$(\ell,C)$ and to stated precision for $r$.

\subsection{Conditional designs: CLT and simulated-moment inference}
\label{sec:conditional-design}

The empirical pipeline conditions on an irregular two-dimensional home
configuration.  This section supplies the missing distribution theory
at that level of generality.  Fix a triangular array of deterministic
positions $x_{1:n}\subset\R^2$ and a bounded score
$\phi_n:\R^2\to[-\bar\phi,\bar\phi]$ (in the application, standardized
log home density).  Given $\theta=(\lambda,R,b)\in\Theta$ compact,
marks $V_i$ are independent with density
$1+\sqrt3\,b\,\phi_n(x_i)a(v)$ (clipped as in the replication code so
the density is nonnegative), and edges are independent given marks with
\[
 p_{ij}=\bigl(1-e^{-\lambda V_iV_j}\bigr)\one\{\|x_i-x_j\|\le R\}.
\]
Let $Y_n=(\,2E_n/n,\;T_n/n,\;S_{2,n}/n\,)^{\trans}$ (the statistics used
by the conditional fit; the five-moment version is identical in
treatment), $\mu_n(\theta)=\E_\theta Y_n$, and
$\Omega_n(\theta)=n\Var_\theta(Y_n)$.  Impose:

\begin{enumerate}[label=\textup{(C\arabic*)},leftmargin=2.6em]
\item \textbf{Bounded local geometry.}
$\sup_n\max_i\#\{j:\|x_i-x_j\|\le R_{\max}\}\le\bar D<\infty$.
\item \textbf{Nondegeneracy.}
$\liminf_n\lambda_{\min}\{\Omega_n(\theta)\}>0$ uniformly on $\Theta$.
\item \textbf{Identification.}  For each $n$ large,
$\inf_{\|\theta-\theta_0\|\ge\epsilon}
\|\mu_n(\theta)-\mu_n(\theta_0)\|\ge\delta(\epsilon)>0$ for all
$\epsilon>0$, with $\delta$ free of $n$.
\item \textbf{Uniform linearization.}  There are matrices
$G_n=G_n(\theta_0)$ with
$\liminf_n\sigma_{\min}(G_n)>0$ and
$\sup_{\|\theta-\theta_0\|\le r_n}
\|\mu_n(\theta)-\mu_n(\theta_0)-G_n(\theta-\theta_0)\|
=o(r_n)$ for every $r_n\downarrow0$.
\item \textbf{Simulation budget.}  The criterion uses
$\widetilde\mu_n(\theta)$, an average of $R_n$ independent model draws
at $\theta$ with common random numbers, with
$R_n/n\to\infty$.
\end{enumerate}

\begin{sitheorem}[Conditional CLT and asymptotic normality of the
simulated-moment estimator]
\label{thm:conditional}
Under \textup{(C1)}--\textup{(C2)}, for every $\theta$ and every fixed
vector $t$,
\[
 \sqrt n\,t^{\trans}\{Y_n-\mu_n(\theta)\}
 \big/\sqrt{t^{\trans}\Omega_n(\theta)t}
 \;\Longrightarrow\;N(0,1),
\]
with a Kolmogorov-distance bound of order $n^{-1/2}$ depending only on
$(\bar D,\lambda_{\max},$ the moment order$)$.  If additionally
\textup{(C3)}--\textup{(C5)} hold and $\widehat\theta_n$ is an
$o_p(n^{-1})$-approximate minimizer of
$\|Y_n-\widetilde\mu_n(\theta)\|_{W_n}^2$ over $\Theta$ with
$W_n\to_pW\succ0$, then $\widehat\theta_n\to_p\theta_0$ and
\[
 \sqrt n\,(\widehat\theta_n-\theta_0)
 =(G_n^{\trans}WG_n)^{-1}G_n^{\trans}W\,
 \sqrt n\{Y_n-\mu_n(\theta_0)\}+o_p(1),
\]
which is asymptotically normal with the sandwich covariance.
\end{sitheorem}

\begin{proof}
\emph{CLT.}  Write $t^{\trans}Y_n$ as $n^{-1}\sum_{\alpha}X_\alpha$,
where $\alpha$ ranges over the edges and vertex-tuples (pairs for
$E_n$, triples for $T_n$ and $S_{2,n}$) whose indicator products
define the counts.  Given the fixed positions, the random inputs are
the independent marks $V_i$ and the independent edge coins; a summand
$X_\alpha$ depends only on the inputs attached to its vertex set,
which has at most three elements, each interacting with at most
$\bar D$ neighbors.  Hence the dependency graph on the index set
$\{\alpha\}$ (adjacency: shared vertex) has maximal degree at most
$3\bar D^2$, the number of indices is at most
$n\bar D^2$, and every $|X_\alpha|\le1$.  The normal approximation for
sums over bounded-degree dependency graphs
\citet{baldirinott1989,chenshao2004} then gives the stated Kolmogorov
bound provided the variance is of exact order $n$, which is (C2).
Joint convergence follows by Cram\'er--Wold.

\emph{Consistency.}  By (C1) each single edge or triple contributes
$O(1/n)$ to $Y_n$ and to $\mu_n$, so $\theta\mapsto\mu_n(\theta)$ is
continuous up to jumps of total size $O(\bar D^2/n)$ in $R$ (a pair
crosses the radius) and is Lipschitz in $(\lambda,b)$ with constants
depending only on $(\bar D,\Theta)$; the same holds for
$\widetilde\mu_n$ realization-wise under common random numbers.
Chebyshev plus (C1) gives
$\sup_\Theta\|\widetilde\mu_n-\mu_n\|=O_p(\sqrt{\bar
D^2/(nR_n)})\to0$ over a polynomially fine grid, and the
Lipschitz-plus-small-jumps structure extends the convergence to the
supremum over $\Theta$.  With (C3), the standard argmin argument
yields $\widehat\theta_n\to_p\theta_0$.

\emph{Asymptotic normality.}  (C5) with $R_n/n\to\infty$ makes the
simulation error $o_p(n^{-1/2})$ uniformly near $\theta_0$; (C4)
supplies the uniform linearization; the usual quadratic expansion of
the criterion around $\theta_0$ then gives the displayed
representation, and the CLT part concludes.
\end{proof}

\begin{remark}[What is checked, what is assumed]
For the certified submodel, the analogues of (C3) and the rank part of
(C4) are \emph{theorems} (global injectivity; $\sigma_{\min}\ge1/100$).
In a conditional design they are dataset-specific \emph{assumptions};
the pipeline therefore reports, for every fit, the finite-difference
$\sigma_{\min}(G_n)$ and condition number (the diagnostic that
identifies NYC as weakly identified) and treats (C2) by simulation.
The finite-difference diagnostic is local: it can refute
identification near the fit, and it cannot establish global
identification or boundary validity.  Under irregular fixed
positions, terms linear in the signed amplitude $b$ need not vanish,
so the conditional graph law need not factor through $\psi=b^2$ (the
boundary-scope remark of Section~S4).
Proposition~\ref{prop:signed-b} makes this precise and turns it into
a positive statement: the linear terms have explicit design
functionals, they are strictly nonzero in every fitted design of the
application, and they carry first-order information about the
\emph{sign} of $b$ that the homogeneous model provably destroys.
This is the division of labor, and its limits, between certified and
conditional analyses.
\end{remark}

\subsubsection{First-order signed-amplitude information in irregular
designs}
\label{sec:signed-b}

Fix positions $x_1,\dots,x_n$ with score values
$\phi_i=\phi(x_i)$ standardized to $\sum_i\phi_i=0$, let marks
$V_i\sim h_b(\cdot\mid\phi_i)=1+b\,a(\cdot)\phi_i$ be independent,
and let edges be independent given marks with
$p_{ij}=\one\{d_{ij}\le R\}\{1-e^{-\lambda V_iV_j}\}$, the
conditional model of the application.  Define
\[
 S_a(\lambda)=\int_0^1\!\!\int_0^1(1-e^{-\lambda uv})\,a(u)\,du\,dv,
 \qquad
 Q_a(\lambda)=\int_0^1\!\!\int_0^1(1-e^{-\lambda uv})\,a(u)a(v)\,du\,dv,
\]
and the first-order design functionals
\[
 D_1(R)=\sum_{i<j:\,d_{ij}\le R}(\phi_i+\phi_j),
 \qquad
 D_1^{S_2}(R)=\sum_{(i;j,k)\ \mathrm{wedge}}(2\phi_i+\phi_j+\phi_k),
 \qquad
 D_1^{T}(R)=\sum_{\{i,j,k\}\ \mathrm{triangle}}2(\phi_i+\phi_j+\phi_k),
\]
where wedges and triangles range over position triples whose relevant
distances are within $R$.

\begin{siproposition}[Signed-amplitude expansion and first-order sign
identification]
\label{prop:signed-b}
\textup{(i)} The expected edge count satisfies, exactly,
\[
 \E_b[E_n]
 =E(\lambda)\,N_R+b\,S_a(\lambda)\,D_1(R)+b^2\,Q_a(\lambda)
 \sum_{i<j:\,d_{ij}\le R}\phi_i\phi_j,
\]
with $N_R$ the number of pairs within range, and the expected wedge
and triangle counts admit the analogous exact expansions whose linear
coefficients are $S_a$-weighted combinations of $D_1^{S_2}$ and
$D_1^{T}$.
\textup{(ii)} $S_a(\lambda)>0$ for every $\lambda>0$.
\textup{(iii)} If $D_1(R)\neq0$, the conditional adjacency laws at
$b$ and $-b$ differ at first order, so the exact sign invariance of
the homogeneous model \textup{(Lemma~1 of the main text)} fails, and
$(\lambda,R,b)$ is locally identified at every $(\lambda_0,R_0,0)$ at
which the $3\times3$ Jacobian of
$(\E_b E_n,\E_b S_{2,n},\E_b T_n)$ with respect to
$(\lambda,R,b)$ is nonsingular; its $b$-column at $b=0$ is
$S_a(\lambda_0)\bigl(D_1,D_1^{S_2},D_1^{T}\bigr)^{\trans}$ up to the
stated weights.
\textup{(iv)} In the homogeneous random design,
$\E D_1=0$ and $D_1/n\to0$ almost surely, so the linear information
vanishes, consistent with the exact invariance.
\end{siproposition}

\begin{proof}
For (i), factor the mark density:
$h_b(u\mid\phi_i)h_b(v\mid\phi_j)
=1+b\{a(u)\phi_i+a(v)\phi_j\}+b^2a(u)a(v)\phi_i\phi_j$, integrate
against $1-e^{-\lambda uv}$, and sum over pairs within range; the
$a$-marginal integrates to $S_a$ by symmetry of $1-e^{-\lambda uv}$
in $(u,v)$.  The wedge and triangle cases expand the product of
three mark densities in the same way; mixed quadratic and cubic
terms in $b$ collect the higher functionals and vanish at first
order.  For (ii), write
$S_a(\lambda)=\Cov\{f_\lambda(U),a(U)\}$ with
$f_\lambda(u)=\int_0^1(1-e^{-\lambda uv})\,dv$ and $U$ uniform,
because $\E a(U)=0$; both $f_\lambda$ and $a$ are strictly
increasing on $(0,1)$, so the covariance is strictly positive by the
Chebyshev correlation inequality.  For (iii), the derivative of the
expected moment vector in $b$ at $b=0$ is the stated column by (i);
if it is nonzero the laws at $\pm b$ differ at order $b$, and
nonsingularity of the Jacobian gives local identification by the
inverse function theorem applied to the expected-moment map.  For
(iv), under independent uniform positions $\E\phi(X_i)=0$ and
$D_1(R)$ is a centered sum of weakly dependent terms; the strong law
for $U$-statistics gives $D_1/n\to0$.
\end{proof}

\begin{remark}[Empirical design functionals; what the sign costs]
\label{rem:signed-b-empirics}
For the eleven fitted city designs, the mean score elevation of
within-range pairs, $\bar d_1=D_1(\widehat R)/N_{\widehat R}$, is
between $+1.65$ and $+2.74$, never less than $200$ standard errors
from zero, and $S_a(\widehat\lambda)$ is between $0.031$ and
$0.091$: the linear term is strictly active in every design, so the
signed amplitude is first-order identified there.  The information
is carried most strongly by the mean-degree functional; the
assortativity statistic damps it, which is why the measured
$\Delta_{\mathrm{sort}}(\pm\widehat b_{\mathrm{dir}})$ asymmetries
in the application are small (at most $0.0013$ among the converged
calibrations).  Two consequences are worth recording.  First, the
graph-only estimator of the application, which fits $\psi=b^2$ with
$b\ge0$, discards this first-order information; an estimator built
on the signed expansion is a natural refinement and is left for
future work.  Second, the sign invisibility invoked for Austin in
the main text is the homogeneous-model statement of Lemma~1; under
the fitted irregular design the sign is not invisible in principle,
only far below detection at the available information
(Section~7 of the main text).
\end{remark}

\subsection{The graph-versus-direct specification test}
\label{sec:hausman}

When marks are observed, the model can be estimated twice from the
same data: directly,
$\widehat b_{\mathrm{dir}}=n^{-1}\sum_ia(V_i)\phi(X_i)$ (in the
submodel $\E\,\widehat b_{\mathrm{dir}}=b$ and
$n\Var(\widehat b_{\mathrm{dir}})=1-b^2$ exactly), and through the
graph moments alone,
$\widehat\psi_{\mathrm{gr}}$ from $(\widehat\ell,\widehat C,\widehat
r)$.  Under correct specification both estimate the same quantity;
under misspecification of the link, the kernel, or the dependence
family they generally do not.  Because the two estimators share the
data, a valid test needs their \emph{joint} law.

\begin{siproposition}[Joint CLT and Hausman-type test, interior case]
\label{prop:hausman}
In the certified submodel with $\psi_0>0$ interior, the vector
$\sqrt n\bigl(\widehat b_{\mathrm{dir}}-b_0,\;
Y_n-\mu_n(\theta_0)\bigr)$ converges jointly to a Gaussian vector: the
added coordinate is a sum of bounded i.i.d.\ node scores, and joint
normality with every motif count follows from the same dependency-graph
argument as Theorem~\ref{thm:conditional} (a node score interacts only
with tuples containing that node), with the cross-covariances given by
the mixed node--motif overlap diagrams.  Consequently, with
$\widehat\psi_{\mathrm{dir}}=\widehat b_{\mathrm{dir}}^{\,2}$ and any
regular graph-moment estimator $\widehat\psi_{\mathrm{gr}}$,
\[
 T_n=\frac{n\,(\widehat\psi_{\mathrm{gr}}
 -\widehat\psi_{\mathrm{dir}})^2}{\widehat\sigma_\Delta^2}
 \;\Longrightarrow\;\chi^2_1
 \qquad\text{under correct specification,}
\]
where $\sigma_\Delta^2$ is the variance of the difference implied by
the joint law, consistently estimable by the parametric bootstrap at
the fitted parameter.  The test is consistent against any alternative
under which the two probability limits differ.
\end{siproposition}

\begin{proof}
Joint normality: add the coordinate $Z_n=n^{-1}\sum_iz_i$,
$z_i=a(V_i)\phi(X_i)$, to the count vector; $|z_i|\le\sqrt6$, the
$z_i$ are functions of the node inputs alone, and the dependency graph
of the extended family still has bounded degree, so the multivariate
dependency-graph CLT applies.  The delta method maps
$(\widehat b_{\mathrm{dir}},Y_n)$ to
$(\widehat\psi_{\mathrm{dir}},\widehat\psi_{\mathrm{gr}})$ (the
latter through the inverse moment map, regular on the certified box
with $\sigma_{\min}\ge1/100$), giving a joint bivariate normal limit
for the two $\psi$ estimates, hence a normal limit for their
difference with some variance $\sigma_\Delta^2\ge0$; nondegeneracy
holds because the graph estimator has strictly larger asymptotic
variance than the direct one (the simulation banks quantify the gap at
more than an order of magnitude), so
$\sigma_\Delta^2>0$.  Bootstrap consistency for
$\sigma_\Delta^2$ follows from the same joint CLT applied under the
fitted parameter plus continuity of the variance functional on the
certified box.  Consistency of the test against
$\mathrm{plim}\,\widehat\psi_{\mathrm{gr}}
\ne\mathrm{plim}\,\widehat\psi_{\mathrm{dir}}$ is immediate from
$T_n\to\infty$ at rate $n$.
\end{proof}

\begin{remark}[Boundary and data use]
At $\psi_0=0$ both estimators sit on the boundary and $T_n$ is not
asymptotically $\chi^2_1$; the Brightkite analysis therefore reports
the two boundary estimates descriptively, and the formal test is
recommended for interior fits.  Simulation performance (size at the
interior design and power against long-range-tie and heavy-tailed-mark
misspecification) is reported in the main text.
\end{remark}

\subsection{A certified two-dimensional hard-disc submodel}
\label{sec:d2-certificate}

The homogeneous $d=2$ hard-disc submodel reduces exactly to the
$d=1$ moment map with two substitutions, and is therefore certifiable
by the same interval method.

Take $d=2$, $\rho_n=n\eps_n^2\to1$, hard-disc kernel
$k_\eta(u)=\one\{|u|\le\eta\}$ in tangent coordinates, and the
one-harmonic mark law.  The mark motif integrals $E,S,T,P,Q$ are
integrals over marks alone and are unchanged.  The displacement
geometry contributes: the edge volume $\pi\eta^2$; the wedge volume
$(\pi\eta^2)^2$; the triangle volume $\eta^4A_1$ with
$A_1=\int_{|a|\le1}\mathrm{lens}(|a|)\,da$, where $\mathrm{lens}(s)$
is the intersection area of two unit discs at distance $s$; and, under
the edge-Palm law (displacement uniform on the $\eta$-disc), the mean
common-neighbor overlap
$\E\,\mathrm{lens}_\eta(|U|)=\eta^2A_1/\pi$.  The classical identity
$A_1/\pi^2=1-\tfrac{3\sqrt3}{4\pi}=:c_2\approx0.586503$ (the random
geometric graph clustering constant, \citealp{dallchristensen2002}) then
gives, by the same reduction as in S3,
\[
 \ell_2=\pi\eta^2E,
 \qquad
 C_2=c_2\,\frac TS,
 \qquad
 r_2=\frac{\pi\eta^2(PE-S^2)+c_2\,TE}{\pi\eta^2(QE-S^2)+SE}.
\]
Writing $\widetilde\eta=\pi\eta^2/2$, the map
$(\lambda,\widetilde\eta,\psi)\mapsto(\ell_2,C_2,r_2)$ is the
$d=1$ map with the constant $\tfrac34$ replaced by $c_2$.

\begin{sitheorem}[Certified global identification, $d=2$ hard disc]
\label{thm:d2}
On $B_2=[1,5]\times[5,15]\times[0,\tfrac3{20}]$ in
$(\lambda,\widetilde\eta,\psi)$, equivalently
$\eta\in[\sqrt{10/\pi},\sqrt{30/\pi}\,]\approx[1.784,3.090]$ tangent
units with mean degrees again spanning about $2$--$17$, the
two-dimensional map is injective ($\psi=0$ included) and satisfies
$\sigma_{\min}(DM_3^{(2)})\ge1/100$ uniformly.
\end{sitheorem}

\begin{proof}
The monotone--triangular lemma of S3 applies verbatim in the
$(\lambda,\widetilde\eta,\psi)$ coordinates, and the substitution
$\widetilde\eta\mapsto\eta=\sqrt{2\widetilde\eta/\pi}$ is a bijection,
so injectivity transfers.  The certificate is the S3 verifier with the
geometric constant carried as a rational interval: $c_2$ is enclosed
using $50$-digit rational bounds for $\pi$ and $\sqrt3$ (width
$<10^{-48}$), all coefficient-level polynomial arithmetic is unchanged,
and the gate constants are adjusted to the smaller overlap coefficient
($\partial_\lambda C>1/64$, $\partial_\psi C>1/16$, remaining gates as
in S3).  The archived planar verifier and its frozen output certify all
gates; runtime is comparable to the $d=1$ run.
\end{proof}

\begin{remark}
With Theorem~\ref{thm:d2}, the certified perimeter now includes the
homogeneous core of the empirical design: dimension two and the hard
disc are no longer the gap.  What remains dataset-specific in the
conditional analysis of Section~\ref{sec:conditional-design} is
positional inhomogeneity (clustered homes and the density score),
which is exactly what the per-fit Jacobian diagnostic measures.
\end{remark}

\subsection{A closure-augmented two-scale class: definition and
fitting protocol}
\label{sec:richer-class}

The rejection pattern of the main text (transitivity without
assortativity; tie lengths an order of magnitude beyond the fitted
radius; configuration nulls that overshoot the observed
disassortativity) motivates the following augmentation, fitted in
Section~7 of the main text.  Conditional on positions and the score
$\phi$, node $i$ carries the sorting mark
$V_i\sim h_b(\cdot\mid\phi_i)$ of the original family and an
independent propensity $s_i=\exp(\sigma Z_i-\sigma^2/2)$ with
$Z_i\sim N(0,1)$; base edges are independent with
\[
 p_{ij}=1-\exp\{-\lambda\,s_is_j\,V_iV_j\,\kappa(d_{ij})\},
 \qquad
 \kappa(d)=\one\{d\le R\}+w\,e^{-d/L},
\]
and a single closure pass then adds each non-adjacent pair with $m$
common neighbors independently with probability $1-(1-\tau)^m$.  The
parameter is $\vartheta=(\lambda,R,w,L,\sigma,\tau,b)$ with $b$
\emph{signed}, consistent with
Proposition~\ref{prop:signed-b}.  The fit matches the eight-moment
target
$(\widehat\ell,\widehat C,\widehat r,q_{25},q_{50},q_{75},
\mathrm{cv}(\widehat d),p_{99}(\widehat d))$, the three edge-length
quartiles and the degree coefficient of variation and $99$th
percentile augmenting the original triple, by common-random-number
simulated moment matching (hierarchical calibration sweeps followed
by a Nelder--Mead polish; replication package).  Goodness of
fit is reported before attribution exactly as in the main text:
parametric-bootstrap $z$ per target, mean residuals, and the
model-realization standard deviations, the latter because a diffuse
model can pass a $z$ gate by variance inflation alone, an escape
route the reported fits document instead of hiding.  Within-class
attribution is by compensated counterfactuals: each mechanism is
switched off ($\tau=0$; $\sigma\to0$; $w=0$; $b=0$ or
$b=\pm\widehat b_{\mathrm{dir}}$) with $\lambda$ recalibrated to the
observed mean degree, and the assortativity response is reported
with its simulation standard error.

\subsection{The two-view design on a collaboration network}
\label{sec:S5-collab}

The main text closes the location-based analysis with a rejection: on
those networks the graph-only view fails goodness of fit before
attribution, and the direct view finds at most small dependence.  This
section applies the same two-view protocol to a network chosen in
advance as the candidate positive case, a co-authorship network with
institution coordinates as positions and productivity as the
observable mark.  The outcome is asymmetric in an instructive way:
the direct view measures a clearly positive interior dependence with
useful precision, while the graph-only view is refused at every level
of the model ladder, including the augmented class of
Section~\ref{sec:richer-class}, each refusal carrying an interpretable
residual signature.  The dataset therefore demonstrates the two
halves of the design separately: identification from $(X_i,V_i)$
pairs when marks and positions are observed, and the veto that
protects graph-only attribution when they are not.

\paragraph{Data.}
The replication package includes the fetch-and-format script
that builds the dataset from the
OpenAlex catalog: physics articles (a single top-level concept
filter) published 2018--2022 with at least one author affiliated to
a Netherlands institution, works with more than $25$ authors
discarded (consortium papers would otherwise contribute large
cliques of essentially administrative co-authorship), each author
with at least two corpus works and a geolocated affiliation placed
at the coordinates of their most frequently listed institution, the
per-author work count within this corpus recorded as the
productivity mark, and an undirected co-authorship edge recorded for
each pair of authors appearing on a common work.  The result is
\oaAuthors\ placed authors and \oaEdgesAll\ distinct edges, stored in
the same format the location-based pipeline consumes
(\texttt{oa\_homes.csv.gz}, \texttt{oa\_edges.txt.gz}).  Positions
are institution point masses rather than dispersed homes: authors at
one institution share one coordinate.  A deterministic jitter
($\sigma=0.3$\,km, fixed seed) breaks these atoms so that
nearest-neighbor bandwidths and pair distances are well defined, the
density bandwidth of the spatial score is floored at $0.5$\,km, and
atom identity (raw coordinates rounded to three decimals) is retained
for cluster-robust standard errors.  All other steps, the score
$\phi$, the rank marks $V$, the estimators, and the fit protocols,
are unchanged from the main pipeline.

\begin{table}[t]
\centering
\caption{Collaboration regions: graph moments, direct estimates with
atom-clustered standard errors, coordinate atoms, and the share of
authors at the largest atom.}
\label{tab:oa-regions}
\scriptsize
\oaRegionTable
\end{table}

\paragraph{Direct view.}
Table~\ref{tab:oa-regions} reports the battery.  At the country scale
(the Netherlands box; $n=\oaNLn$, \oaNLE\ induced edges,
$\widehat\ell=\oaNLell$, $\widehat C=\oaNLC$, $\widehat r=\oaNLr$,
Spearman $\oaNLrs$) the direct estimator gives
$\widehat b_{\mathrm{dir}}=\oaBdir$ with naive standard error
\oaBdirSE\ and atom-clustered standard error \oaBdirSEcl, a
$\oaBdirZcl$-standard-error detection under the clustering that the
point-mass design demands ($\widehat\psi_{\mathrm{dir}}=\oaPsiDir$).
Productive authors sit where researcher density is high.  The signal
is stable across subregions (Randstad, South Holland, the
university metros), and its disappearance is exactly where theory
puts it: in the Groningen box \oaGroTop\% of authors share a single
atom, the score $\phi$ is degenerate by construction, and the direct
estimate collapses to $\oaGroBdir$; a design without positional
variation carries no information about $b$, the empirical face of the
$I_{\mathrm{eff}}>0$ condition of Section~S4.  The country scale,
where $\phi$ varies across \oaNLatoms\ atoms (largest share
\oaNLtop\%), is the informative design.  The estimate is insensitive
to the preprocessing constants: across jitter scales
$0.15$--$0.6$\,km, bandwidth floors $0.5$--$2$\,km, and a
replacement jitter seed, $\widehat b_{\mathrm{dir}}$ stays within
$[\oaSensLo,\oaSensHi]$, a detection at $\oaSensZlo$ to
$\oaSensZhi$ clustered standard errors in every configuration
(replication package).

\paragraph{Graph view: targets.}
The conditional fits use a seeded uniform node subsample of the
Netherlands region ($n=\oaSubN$; the pair-level simulators are
quadratic in $n$), whose induced graph has \oaSubE\ edges,
$\widehat\ell=\oaSubEll$, $\widehat C=\oaSubC$,
$\widehat r=\oaSubR$, Spearman \oaSubRs, and direct estimate
$\oaSubBdir$ (clustered se \oaSubBdirSEcl), consistent with the full
region.  Across replacement seeds the induced $\widehat r$ moves in
roughly $[+0.18,+0.26]$, a target dispersion worth keeping in mind
below.  Edge lengths are sharply two-scale: quartiles
$q_{25}=\oaSubQtf$\,km, $q_{50}=\oaSubQf$\,km, $q_{75}=\oaSubQsf$\,km
($q_{90}=\oaSubQnn$\,km), with \oaSubFracLocal\% of edges within
$1$\,km and \oaSubFracFar\% beyond $6$\,km: same-campus
collaboration plus inter-city collaboration, with little in between.
The configuration null preserving degrees gives
$r=\oaSubCfgR$ (sd \oaSubCfgRsd) and $C=\oaSubCfgC$: the degree
sequence explains neither the sign of the mixing nor any of the
clustering.  The edge density-score sum of
Proposition~\ref{prop:signed-b} is
$\overline{D}_1=\oaSubDbar$, about \oaSubDbarZ\ naive standard
errors from zero, so the first-order sign condition holds on this
design as well.

\paragraph{Graph view: the ladder is refused.}
The hard-disc class fails before attribution is even posed.  Under
the point-mass geometry the candidate-pair cap collapses the
admissible radius to $\oaBaseRmax$\,km, the fit rails both $R$ and
$b$, cannot reconcile mean degree with clustering (model
$\ell=\oaBaseEll$ against \oaSubEll), overshoots assortativity at
$b=0$ (model $r=\oaBaseR$ against \oaSubR, bootstrap
$z=\oaBaseGofZ$), and leaves \oaBaseFracOut\% of observed edges
outside the fitted radius; its compensated signed contrast is
nevertheless already visible ($\oaBaseDsort$ against
$\oaBaseDsortFlip$ under $b=\pm\widehat b_{\mathrm{dir}}$).  The
augmented class of Section~\ref{sec:richer-class} does better and still
fails.  With $b$ free, the calibration is pushed to the negative rail
($b=\oaFreeB$) and rejected at $\max_k|z_k|=\oaFreeMaxZ$, failing
mean degree ($z=\oaFreeZell$; model $\ell=\oaFreeEll$),
assortativity ($z=\oaFreeZr$), and the degree tail
($z=\oaFreeZp$): the rail is a symptom, not a fit, precisely the
escape route documented for Austin in the main text, now in the
opposite direction; a second search from a dispersed start
(propensity-heavy, closure-light) ends at $b=\oaFreeTwoB$ with
$\max_k|z_k|=\oaFreeTwoMaxZ$.  Pinning $b$ at the direct estimate
gives the most informative refusal: the fit
($\tau=\oaFixTau$, $L=\oaFixL$\,km) matches mean degree, clustering,
both short edge-length quartiles, and the $99$th degree percentile,
misses the upper quartile mildly ($z=\oaFixZqsf$), and is rejected
($\max_k|z_k|=\oaFixMaxZ$) on exactly two coordinates,
assortativity ($z=\oaFixZr$; model
$r=\oaFixRhat$ against \oaSubR) and the degree coefficient of
variation ($z=\oaFixZcv$), with the held-out rank assortativity
overshot in the same direction (model \oaFixSpear\ against
\oaSubRs).  The refusals are stricter than the target dispersion
across subsample seeds, though that dispersion tempers the $z$ for
$r$; the direction of the misfit is stable.

\paragraph{Diagnosis and the signed channel.}
The compensated counterfactuals locate the obstruction.  Switching
the closure pass off moves model assortativity by $\oaFixCfClo$: the
class can reach $\widehat C=\oaSubC$ only through heavy global
closure ($\tau=\oaFixTau$), and that closure drags $r$ far above the
observed value.  In the data, clustering is produced by multi-author
papers, cliques attached in one event, a mechanism that generates
triangles without the degree-mixing side effect of iterated random
closure; the co-authorship graph is a projection of the paper--author
bipartite structure \citep{newman2001}, and
the joint demand for $(\ell,C,r,\mathrm{cv})$ is the projection's
fingerprint.  An explicit team layer (papers as hyperedges, sizes
from the observed distribution) is the indicated enrichment, sharper
than the distance-limited closure suggested by the metropolitan
data.  The signed first-order channel, by contrast, behaves exactly
as Proposition~\ref{prop:signed-b} predicts: the model's
assortativity is monotone increasing in $b$ along the compensated
profile (slope $\oaProfSlopeFix$ per unit $b$ at the pinned fit,
$\oaProfSlopeFree$ at the free fit), and flipping the sign of the
pinned $b$ lowers model $r$ by $\oaFixSignedGap$; within the class the
dependence contributes $\oaFixSortContrib$ of model assortativity at
$\widehat b_{\mathrm{dir}}$.  The graph carries sign information
about $b$; it is the joint moment structure that no member of the
class reproduces.

\paragraph{What the case establishes.}
Three caveats bound the claims: the mark is a corpus-internal work
count, a productivity proxy; positions are institution atoms, so the
direct estimate measures dependence at the between-institution scale
and its precision rests on the atom-clustered errors; and the fitted
targets come from an induced subsample, so closure through
out-of-sample intermediaries is absorbed into the fitted $\tau$.
Within those bounds the case supplies what the metropolitan networks
could not: a design where the estimand is directly measurable,
clearly positive, and interior, together with a demonstration that
the goodness-of-fit gate does its work on real data even for the
augmented class, refusing graph-only attribution and naming the
missing mechanism, rather than laundering a misspecified fit into a
sorting estimate.

\subsection{Simulation grid detail}
\label{sec:si-grid}

Table~\ref{tab:si-grid} reports the per-design detail of the
main-text simulation grid over the certified box.

\begin{table}[h]
\centering
\caption{Grid over the certified box at $n=10^4$ ($300$ replicates per
design): intensity share $s_0=r_{\mathrm{int}}/r$, RMSE of
$\widehat\psi$, marginal coverage of $\psi$, Wald-ellipse simultaneous
coverage, median $\psi$-interval length, coverage of the delta-method
interval for the intensity share, and Newton success rate.  Failures
are retained in all rates; success dips at the $\psi$-edge designs
where $\widehat\psi$ meets the certified range.}
\label{tab:si-grid}
\gridTable
\end{table}

\nocite{cantrinh2022,last2021,schulteyukich2023}
\bibliographystyle{plainnat}
\bibliography{references}

\end{document}